\documentclass[a4paper,10pt]{book}
 \usepackage{amsfonts, amssymb, amscd, amsmath}
 \usepackage{calrsfs, pifont}
 \usepackage{latexsym}
 \usepackage{bbm}
\DeclareMathAlphabet{\mathbb}{U}{bbold}{m}{n}

\begin{document}
 \makeatletter%{@}
\renewcommand{\@cite}[2]{{[#1\if@tempswa, #2\fi]}}
\renewcommand{\@dotsep}{3}
\renewcommand{\@evenhead}{\raisebox{3mm}[\headheight][0pt]%
             {\vbox{\hbox to\textwidth{\thepage \hfill\strut {\text{\small A. E. Gutman, A. G. Kusraev, and S. S. Kutateladze}}\hfill}\hrule}}}
\renewcommand{\@oddhead}{\raisebox{3mm}[\headheight][0pt]%
             {\vbox{\hbox to\textwidth{\strut \hfill{\text{\small The Wickstead Problem}\hfill}\thepage}\hrule}}}

\renewcommand{\section}[1]{\medskip\begin{center}\textbf{#1}\end{center}}

\def\subsec#1{\smallskip{\rm\bf #1}}
\newcommand{\subsubsec}[1]{\hspace{17pt}{\rm #1}}
\newcommand{\Ref}{\medskip\begin{center}\textbf{References}\end{center}}
\newcommand{\bib}[2]{{\leftskip-20pt\baselineskip=11pt\footnotesize\item{}\textsl{#1}\,\ #2}}
\makeatother%@
%%%%%%%%%%%%%%%%%%%%%%%%%%%%%%%%%%%%%%%%%%
\newcommand{\Endproc}{\rm}
\newcommand{\Proclaim}[1]{\smallskip{\textbf{#1\/~}}\sl}
\newcommand{\proclaim}[1]{{\textbf{#1\/~}}\sl}
%%%%%%%%%%%%%%     English      %%%%%%%%%%%%%%%%%%%%%%%%%%
\newcommand{\Lemma}[1]{\smallskip\textbf{{Lemma#1}.\/~}\sl}
\newcommand{\lemma}[1]{\textbf{{Lemma#1}.\/~}\sl}
\newcommand{\Theorem}[1]{\smallskip\textbf{{Theorem#1}.\/ }\sl}
\newcommand{\theorem}[1]{{\textbf{Theorem#1.\/~}}\sl}
\newcommand{\Corollary}[1]{\smallskip\textbf{{Corollary#1}.\/~}\sl}
\newcommand{\corollary}[1]{\textbf{{Corollary#1}.\/~}\sl}
\newcommand{\Proposition}[1]{\smallskip\textbf{{Proposition#1}.\/~}\sl}
\newcommand{\proposition}[1]{\textbf{{Proposition#1}.\/~}\sl}
\newcommand{\Crit}[1]{\smallskip\textbf{{Criterion#1}.\/ }\sl}
%-------------
\newcommand{\defin}[1]{{\sc Definition#1.}}
\newcommand{\Definition}[1]{\smallskip{\sc Definition#1.}}
\newcommand{\Remark}[1]{\smallskip{\sc Remark#1.}~}
\newcommand{\Ex}[1]{\smallskip{\sc Example#1.}}
\newcommand{\Problem}[1]{{\sc Problem#1:}}
\newcommand{\Note}[1]{{\sc Note#1.}}
%%%%%%%%%%%%%%%%%%%%%%%%%%%%%%%%%%%%%%%%%%%%%%%%%%%%%%%%%%%%%%

\def\beginproof{\smallskip\par\text{$\vartriangleleft$}}
\def\endproof{\text{$\vartriangleright$}}
\def\gather{\gathered}
\def\endgather{\endgathered}
\def\align{\aligned}
\def\endalign{\endaligned}
\def\nomathbreak{}
\def\allowmathbreak{}
\def\Sb#1\endSb{_{\substack{#1}}}
\def\section#1{\bigskip\begin{center}{\scshape#1}\end{center}\medskip}
 \def\dom{\mathop{\fam0 dom}}
 \def\On{\mathop{\fam0 On}}
 \def\ZFC{\mathop{\fam0 ZFC}}
 \def\endproc{\rm}
 \def\Fin{\mathop{\fam0 fin}\nolimits}
 \def\Id{\mathop{\fam0 Id}}
\def\Rightharpoonup{\mathrel{\vcenter{\offinterlineskip\halign{\hfil##\hfil\cr
$\rightharpoonup$\cr\hskip.6pt\hrulefill\hskip.6pt\cr\noalign{\vskip.5pt}}}}}

\def\upwardarrow{\mathord{
  \hbox to 5pt{\hss$\vcenter{\hbox to 2.4pt{\hss$\mathchar"222$\hss}\hrule}\hss$}
}}
\def\downwardarrow{\mathord{
  \hbox to 5pt{\hss$\vcenter{\hrule\hbox to 2.4pt{\hss$\mathchar"223$\hss}}\hss$}
}}

\let\mathbbt\mathbb
\let\mathbbm\mathbb
\def\assign{:=} % You may redefine it as {\!:=}, but I suggest leaving it plain.
\let\iff\leftrightarrow % A workaround for proper coloring in the FAR editor
\let\la\leftarrow % A workaround for proper coloring in the FAR editor
\let\ra\rightarrow % A workaround for proper coloring in the FAR editor
%%%%%%%%%%%%%%%%%%

%%%%%%%%%%%%%%%%%%%%%%%%%%%%

\newpage\thispagestyle{empty}

 \normalsize
 \noindent
 \rm

\newcommand{\shortpage}{\enlargethispage{-\baselineskip}}
\def\goth#1{\mathfrak{#1}}
\def\osum{\mathop{o\text{-}\!\sum}}

\def\End{\mathop{\fam0 End}}
\def\Ext{\mathop{\fam0 Ext}}
\def\mix{\mathop{\fam0 mix}}
\def\sker{\mathop{\fam0 sker}}
\def\im{\mathop{\fam0 im}}
\def\uc{\mathop{\fam0 uc}\nolimits}
\def\sgn{\mathop{\fam0 sgn}\nolimits}
\def\Orth{\mathop{\fam0 Orth}\nolimits}
\def\Otimes{\overset{\underline{\phantom{\ \,}}}\otimes}
 \def\Ord{\mathop{\fam0 Ord}}
 \def\Cn{\mathop{\fam0 Cn}}
 \def\Card{\mathop{\fam0 Card}}

 \newpage\thispagestyle{empty}
 \

 \begin{center}
 {\LARGE
  A.\,E.\,Gutman\\
  A.\,G.\,Kusraev\\
  S.\,S.\,Kutateladze

  \vskip50mm

 {\bf THE WICKSTEAD PROBLEM}}

 \end{center}

 \newpage\thispagestyle{empty}

 \normalsize
 \bf Gutman~A.\,E.,\!\footnote{\ The work of the first author is supported by the Russian Science Support Foundation.} Kusraev~A.\,G., Kutateladze~S.\,S.~
 \rm The Wickstead Problem. 2007. 44~p.

 \bigskip
 {
 In 1977 Anthony Wickstead raised the question of the
 conditions for all band preserving linear operators to be
 order bounded in a vector lattice. This article overviews
 the main ideas and results on the Wickstead problem and
 its variations, focusing primarily on the case of band
 preserving operators in a~universally complete vector
 lattice.
}
 \bigskip

 \medskip

 
 \vfill

 {\parindent=10cm

 \copyright\ A.~E.~\,Gutman,~ 2007

 \copyright\ A.~G.~Kusraev,~ 2007

 \copyright\ \,S.\,~S.\,~Kutateladze,~ 2007

 \bigskip

 }

 \newpage\thispagestyle{empty}

 \large
 \begin{center}
 \bf  THE WICKSTEAD PROBLEM%\footnote[1]{Supported
 %by  the Russian Foundation for Basic Research,
 %Grant No.~06-01-00622.}
 \\[5mm]
 \sc  A. E. Gutman, \ A. G. Kusraev, \ S. S. Kutateladze
 \end{center}

\vskip 1truecm \rightline{\small\it To Anthony Wickstead on his
sixtieth birthday}

 \bigskip
 \bigskip
 \def\contline#1#2#3#4{\hangafter=1\hangindent=40pt #2\ \rm #3\ \dotfill\ #4\medskip}

  \contline{contpart}{Introduction}{}{5}

 \contline{contpart}{Part~1. Locally One-Dimensional Vector Lattices}{}{7}

 \contline{contpart}{\quad 1.1. Band Preserving Operators}{}{7}

 \contline{contpart}{\quad 1.2. A Local Hamel Basis}{}{8}

 \contline{contpart}{\quad 1.3. $\sigma$-Distributive Boolean Algebras}{}{10}

 \contline{contpart}{Part~2.  Boolean  Approach}{}{13}

 \contline{contpart}{\quad 2.1. Representation of a Band Preserving Operator}{}{14}

 \contline{contpart}{\quad 2.2. Representation of a~Locally One-Dimensional
 Vector Lattice}{}{16}

 \contline{contpart}{\quad 2.3. Dedekind Cuts and Continued Fractions in a Boolean Valued Model}{}{19}

 \contline{contpart}{Part~3. Automorphisms and Derivations}{}{20}

 \contline{contpart}{\quad 3.1. Band Preserving Operators in Complex Vector Lattices}{}{20}

 \contline{contpart}{\quad 3.2. Automorphisms and Derivations on the Complexes}{}{22}

 \contline{contpart}{\quad 3.3. Automorphisms and Derivations on Complex $f$-Algebras}{}{24}

 \contline{contpart}{Part~4. Variations on the Theme}{}{26}

 \contline{contpart}{\quad 4.1. The Wickstead Problem in Lattice Ordered Modules}{}{27}

 \contline{contpart}{\quad 4.2. The Wickstead Problem for Bilinear Operators}{}{28}

 \contline{contpart}{\quad 4.3. The Noncommutative Wickstead Problem}{}{30}

 \contline{contpart}{Part~5. Comments}{}{31}

 \contline{contpart}{\quad 5.1. Comments on Part 1}{}{31}

 \contline{contpart}{\quad 5.2. Comments on Part 2}{}{32}

 \contline{contpart}{\quad 5.3. Comments on Part 3}{}{33}

 \contline{contpart}{Appendix. Boolean Valued Analysis}{}{34}

 \contline{contpart}{\quad A1. Boolean Valued Universes}{}{34}

 \contline{contpart}{\quad A2. Escher Rules}{}{35}

 \contline{contpart}{\quad A3. Boolean Valued Numbers, Ordinals, and Cardinals}{}{39}

 \contline{contpart}{References}{}{42}

 \newpage\thispagestyle{empty}

 \vspace*{50pt}

 \begin{center}
 \Large{\emph{December 13, 2007\\ is the date of the sixtieth birthday of\\
 Professor Anthony Wickstead\\ who is associated with Queens'
 University at Belfast.\\
 Professor Wickstead works in positivity,\\
 a vast and attractive area of functional analysis\\
 which borders  many powerful theories\\ of modern mathematics.\\
 His contributions to positivity\\ has brought him world fame and recognition.\\
 He was the first Editor-in-Chief\\ of the international journal ``Positivity.''\\
 His efforts and contributions in this capacity\\
 have made this journal a natural epicenter\\ of research into
 positivity.\\
 %Russian mathematicians are grateful\\
 %to Professor Wickstead\\
 %not only for his scientific contributions\\
 %but also for the charm and kindness of his personality.\\
 Russian mathematicians appreciate not only \\
 the scientific contribution of
 Professor Wickstead\\ but also his charm and personality.\\
 On behalf of his colleagues in Russia\\ and the  Editorial Board\\
 of {\it Vladikavkaz Mathematical Journal}\\
 we heartily congratulate Professor Wickstead\\
 and wish him many happy returns of the day.}
 }
 \end{center}

\begin{flushright}
 \Large{\emph{A. Kusraev, S. Kutateladze}}
\end{flushright}

 \newpage

 \section{Introduction}

 %\noindent
 \textbf{WP:}~
 \textsl{When are we so happy in a vector lattice that all band preserving
 linear operators turn out to be order bounded?}
 %\smallskip
 This question was raised by~Wickstead in~\cite{Wic1}. The
 answer depends on the~vector lattice in which the~operator in
 question acts. There are several results that guarantee
 automatic order boundedness for a~band preserving operator
 acting in concrete classes of vector lattices
(cp.~\cite[Theorem~2]{AVK}, \cite[Theorems~3.2 and
 3.3]{AVK1}, and \cite[Corollary~2.3]{MW}). However, in this
 article we focus our attention on the case of universally
 complete vector lattices.

 Abramovich, Veksler, and~Koldunov
 were the first to announce an example of an~order unbounded
 band preserving operator in~\cite[Theorem~1]{AVK}.
 Later these authors \cite[Theorem~2.1]{AVK1} as well as
 McPolin and~Wickstead \cite[Theorem~3.2]{MW} showed
 that all band preserving operators in a universally complete
 vector lattice~$E$ are bounded automatically if and only if
 $E$ is locally one-dimensional. The Wickstead
 problem in the class of universally complete vector lattices
 was thus reduced to the characterization of locally one-dimensional
 vector lattices.

 This led to another problem posed by Wickstead \cite{AW}:
 {\sl Is the class of locally one-dimensional vector lattices
 coincident with the class of discrete vector lattices?}
 Gutman gave the~negative answer in~\cite{Gut6}:
 {\sl There is a continuous (purely nonatomic) locally
 one-dimensional universally complete   vector lattice}
 (cp.~\cite{Gut5, Gut1}). Also, Gutman described
 the bases of locally one-dimensional  universally
 complete vector lattices: these are exactly
 $\sigma$-distributive complete Boolean algebras.

 Furthermore, it is well known in Boolean valued analysis that the
 condition for a~universally complete vector lattice to be locally
 one-dimensional is related to the structure of the reals $\mathcal
 R$ inside an appropriate Boolean valued model ${\mathbb
 V}^{(\mathbb B)}$. In more detail  the situation is as follows
 (cp.~\cite{BVA}): By the Gordon Theorem,  each universally
 complete vector lattice may be represented as the descent
 $\mathcal R{\downarrow}$ of the Boolean valued reals
 $\mathcal R$,  while the image of the standard reals
 $\mathbbm{R}$
 (under the canonical embedding of the standard universe
 $\mathbbm{V}$ into the Boolean valued universe
 ${\mathbbm{V}}^{(\mathbbm{B})}$) is the subfield
 ${\mathbbm{R}}^{\scriptscriptstyle\wedge}$
 of~$\mathcal{R}$ inside~${\mathbb V}^{(\mathbb B)}$.
 It is easy and well-known in other terms that
 $\mathcal R{\downarrow}$ is locally one-dimensional
 if and only if ${\mathbbm{R}}^{\scriptscriptstyle\wedge}=\mathcal R$.
 The same is true for Boolean valued complexes $\mathcal C$
 and the image ${\mathbbm{C}}^{\scriptscriptstyle\wedge}$
 of the standard reals $\mathbbm{C}$.

 The Boolean  approach to band preserving operators
 as developed by~Kusraev in \cite{Kus1} reveals
 new interconnections. For example, the construction of
 an~order unbounded band preserving operator can be
 carried out inside an~appropriate Boolean valued universe
 by using a~Hamel basis of the reals $\mathcal R$ considered
 as a~vector space over its subfield
 ${\mathbbm{R}}^{\scriptscriptstyle\wedge}$ (cp.~\cite{DOP, IBA}).
 Of course, some important properties of~$\mathcal{R}{\downarrow}$
 are connected with the structure of the reals~$\mathcal{R}$
 as a~vector space over~${\mathbbm{R}}^{\scriptscriptstyle\wedge}$.
 In~particular, using a~Hamel basis, we can construct a~discontinuous
 $\mathbbm{R}^{\scriptscriptstyle\wedge}$-linear function
 in~$\mathcal{R}$ which gives an~order unbounded band preserving
 linear operator in the universally complete vector lattice
 $\mathcal{R}{\downarrow}$.

 As was demonstrated by~Kusraev in~\cite{Kus2},
 similar constructions can be carried out on using
 a~{\it transcendence basis\/} instead of a~Hamel basis.
 This approach yielded the new characterizations
 of universally complete  vector lattices
 with~$\sigma$-distributive base in terms
 of narrower classes of band preserving linear
 operators, namely, of derivations and automorphisms.
 In particular, working with a~transcendence basis, we
 can construct a~discontinuous
 $\mathbbm{C}^{\scriptscriptstyle\wedge}$-derivation and
 $\mathbbm{C}^{\scriptscriptstyle\wedge}$-automorphism
 in~$\mathcal{C}$ which gives an~order unbounded band
 preserving derivation or automorphism in~$\mathcal{C}{\downarrow}$.

 Summarizing the results of \cite{AVK, Gut6, Kus1, Kus2, MW}
 on the Wickstead problem, we can state the following

 \Theorem{ WP}Assume that $G$ is a~universally complete
 vector lattice with a~fixed order unity   $\mathbbm{1}$,
 while $G_{\mathbbm{C}}$ is the complexification~of~$G$, and\/
  $\mathbbm{B}\assign \mathfrak{E}(G)\assign \mathfrak{E}({\mathbbm{1}})$ is
 the Boolean algebra of all components of\/ $\mathbbm{1}$.
 Assume further that ${\mathcal R}$ and ${\mathcal C}$ stand for
 the reals and the complexes inside the Boolean valued
 universe ${{\mathbbm V}}^{({\mathbbm B})}$. Then the
 following are equivalent:
 \smallskip

 \subsubsec{WP(1)}~${\mathbbm B}$ is $\sigma$-distributive;
 \smallskip

 \subsubsec{WP(2)}~${\mathcal R}={{\mathbbm
 R}}^{\scriptscriptstyle\wedge}$ inside~${{\mathbbm V}}^{({\mathbbm
 B})}$;
 \smallskip

 \subsubsec{WP(2$^\prime$)}~${\mathcal C}= {{\mathbbm
 C}}^{\scriptscriptstyle\wedge}$ inside~${{\mathbbm V}}^{({\mathbbm B})}$;
 \smallskip

 \subsubsec{WP(3)}~$G$ is locally one-dimensional;
 \smallskip

 \subsubsec{WP(3$^\prime$)}~$G_{\mathbbm{C}}$ is locally one-dimensional;
 \smallskip

 \subsubsec{WP(4)}~Every band preserving linear operator
 in $G$ is order bounded;
 \smallskip

 \subsubsec{WP(4$^\prime$)}~Every band preserving linear operator
 in $G_{\mathbbm{C}}$ is order bounded;
 \smallskip

 \subsubsec{WP(5)}~There is no nontrivial\/ ${\mathbbm
 R}$-derivation~in the $f$-algebra $G$;
 \smallskip

 \subsubsec{WP(5$^\prime$)}~There is no nontrivial\/
 ${\mathbbm C}$-derivation~in the complex $f$-algebra $G_{\mathbbm{C}}$;
 \smallskip

 \subsubsec{WP(6)}~Each band preserving endomorphism of
 the complex $f$-algebra $G_{\mathbbm{C}}$ is~a~band
 projection;
 \smallskip

 \subsubsec{WP(7)}~There is no band preserving automorphism
 of\/~$G_{\mathbbm{C}}$ other than the identity.
 \smallskip

 %\subsubsec{WP(8)}~There is no nonzero separately band
 %preserving antisymmetric bilinear operator from
 %$G\times G$ to $G$;
 %\smallskip

 %\subsubsec{WP(9)}~Every separately band preserving bilinear
 %operator from $G\times G$ to $G$ is symmetric;
 %\smallskip

 %\subsubsec{WP(10)}~Every separately band preserving bilinear
 %operator from $G\times G$ to $G$ is order bounded.
 \Endproc

 %\smallskip
 The goal of this article is to examine the Wickstead
 problem for universally complete vector lattices and to
 prove the above theorem. The reader can find the necessary
 information on the theory of vector lattices in~\cite{AB, DOP, Z};
 Boolean valued analysis, in~\cite{Bell, BVA, IBA};
 and field theory, in~\cite{Bou1, Waer,  ZS}. Some aspects
 of the Wickstead problem are also presented
 in~\cite[Chapter~5]{DOP},
 \cite[Section 10.7]{IBA}, and \cite{Kus4}.

 By a vector lattice  throughout the sequel we will mean
 a real Archimedean vector lattice,
 unless specified otherwise. We let $\assign $
 denote  the assignment by definition, while
 $\mathbbm{N}$, $\mathbb{Z}$, $\mathbbm{Q}$, $\mathbbm{R}$, and
 $\mathbbm{C}$ symbolize the naturals, the integers,
 the rationals, the reals, and the complexes.
 We denote  the
 Boolean algebras of bands and band projections in a vector
 lattice~$E$  by  $\mathfrak{B}(E)$ and $\mathfrak{P}(E)$;
 and we let $\mathfrak{E}(\mathbbm{1})$ stand for
 the Boolean algebra of all components of $\mathbbm{1}$.
 \newpage

 \section{PART~1. LOCALLY ONE-DIMENSIONAL VECTOR LATTICES}

 In this part we introduce locally one-dimensional  vector lattices and $\sigma$-distributive Boolean algebras and
 prove that the following  are equivalent for each universally
 complete vector lattice
 $G$ with  base $\mathbbm{B}\assign \mathfrak{B}(G)$, the
 complete Boolean algebra of  bands in $G$:
 \smallskip

 \subsubsec{WP(1)}~${\mathbbm B}$ is $\sigma$-distributive;
 \smallskip

 \subsubsec{WP(3)}~$G$ is locally
 one-dimensional;

 \smallskip
 \subsubsec{WP(4)}~Every band preserving linear operator
 in $G$ is order bounded.

 \section{1.1. Band Preserving Operators}

 In this section we introduce the class of band preserving
 operators and  briefly overview some properties of
 orthomorphisms.

 \subsec{1.1.1.}~Consider a vector lattice $E$ and let
 $D$ be a~sublattice of~$E$. A~linear operator $T$ from $D$ into $E$ is
 \textit{band preserving} %\footnote{In Russian  literature the term \textit{nonextending} is also in use.}
 provided that one (and hence all) of the following  holds:
 %$$
 %\gather
 %Te\in \{e\}^{\perp\perp}\quad (e\in D),\\
 %e\perp f\Rightarrow Te\perp f\quad (e\in D,\,f\in E),\\
 %T(K\cap D)\subset K\quad (K\in\goth B(E)).
 %\endgather
 %$$

 \subsubsec{(1)}~$e\perp f$ implies $Te\perp f\ \ (e\in D,\,f\in E)$,

 \subsubsec{(2)}~$Te\in \{e\}^{\perp\perp}\ \ (e\in D)$ \
 (the disjoint complements are taken in $E$),

 \subsubsec{(3)}~$T(K\cap D)\subset K\ \ \bigl(K\in\goth B(E)\bigr)$.\\
 If $E$ is a vector lattice with the principal projection
 property and $D\subset E$ is an order dense ideal, then
 a linear operator $T:D\to E$ is band preserving if and
 only if $T$ commutes with  band projections; i.e.,

 \subsubsec{(4)}~$\pi Tx=T\pi x\ \ (\pi\in\mathfrak{P}(E),\,x\in D)$.

 \subsec{1.1.2.}~A~band preserving operator $T$ in $E$ need
 not be order bounded (cp.~Sections 1.2 and 1.3 below). However,
 the greatest order ideal $A_T$ in $E$ such that $T$ is order
 bounded on $A_T$ is a band (cp.~\cite{Pag1}). Now, if $A_T$ is a
 projection band then $A_T^\perp$ does not include any nonzero
 order ideal on which $T$ is order bounded. Thus, if $E$
 has the projection property then to each band preserving
 operator $T$ in $E$ there is a band projection $\pi$ such
 that $\pi T$ is order bounded and $\pi^\perp T$ has no
 order bounded components; i.e., $\rho T$ is not order bounded
 for any nonzero $\rho\leq\pi^\perp$.

 \subsec{1.1.3.}~An~order  bounded band preserving operator
 $\pi:D\to E$  on an~order  dense ideal $D\subset E$
 is an~\textit{extended orthomorphism\/} of~$E$
 (cp.~\cite{LS}). Since an extended orthomorphism is
 disjointness preserving, it is also regular according to
 the Meyer Theorem \cite{Mey1, DuM}. Let $\Orth(D,E)$
 signify the set of all extended orthomorphisms of $E$ %/?
 that are defined on a fixed order  dense ideal $D$. An extended %/?
 orthomorphism $\alpha\in\Orth(E,E)$  on the whole
 space $E$ is  an~\textit{orthomorphism}. The
 collection of all orthomorphisms $\Orth(E)$ of $E$ is a
 vector lattice under the pointwise algebraic and lattice
 operations. Let $\mathcal{Z}(E)$ stand for the order
 ideal generated by the identity operator $I_E$ in
 $\Orth(E)$. The space $\mathcal{Z}(E)$ is often called
 the \textit{ideal center} of~$E$.

 \Proclaim{1.1.4.}Every extended orthomorphism in
 a~vector lattice is order continuous. All extended
 orthomorphisms commute with one another.
 \Endproc

 \subsec{1.1.5.} The space of extended orthomorphisms
 $\Orth^{\infty}\,(E)$ is defined as follows: Denote by
 $\mathfrak{M}$ the collection of all pairs $(D,\pi)$, where
 $D$ is an order  dense ideal in~$E$ and $\pi\in\Orth(D,E)$.
 Elements $(D,\pi)$ and $(D',\pi')$ in $\goth M$ are
 announced {\it  equivalent\/} \big(in writing
 $(D,\pi)\sim (D',\pi')$\big) provided that the orthomorphisms
 $\pi$ and $\pi'$ coincide on~$D\cap D'$.
 The factor set  $\mathfrak{M}/{\sim}$ of $\mathfrak{M}$ by $\sim$
 is denoted by  $\Orth^\infty (E)$. The set $\Orth^\infty (E)$
 becomes a~vector lattice under the
 pointwise addition, scalar multiplication, and lattice
 operations. Moreover, $\Orth^\infty (E)$ is
 an~ordered algebra under composition. We will identify each
 orthomorphism $\pi\in\Orth(E)$ with the corresponding
 coset in $\Orth^\infty(E)$.

 \subsec{1.1.6.}~We now  list some useful results on
 orthomorphisms that can be found in~\cite{AB, Lux2, LS, Z}.
 \smallskip

 \subsubsec{(1)}\proclaim{}The ordered algebra
 $\Orth^\infty (E)$ is  a~laterally complete
 semiprime $f$-algebra with unity $I_E$. Moreover,
 $\Orth(E)$ is an $f$-subalgebra of $\Orth^\infty(E)$ and
 $\mathcal{Z}(E)$ is an $f$-subalgebra of $\Orth(E)$.
 \Endproc
 \smallskip

 \subsubsec{(2)}\proclaim{}Every Archimedean $f$-algebra $E$ with unity $\mathbbm{1}$
 is algebraically and latticially isomorphic to the
 $f$-algebra of orthomorphisms of $E$. Moreover, the ideal
 in~$E$ generated by $\mathbbm{1}$ is mapped onto $\mathcal{Z}(E)$.
 \Endproc
 \smallskip

 \subsubsec{(3)}\proclaim{}If $E$ is an~order complete vector
 lattice then $\Orth^\infty(E)$ is a~universally
 complete vector lattice and $\Orth(E)$ and $\mathcal{Z}(E)$ are order  dense ideals.
 \Endproc
 \smallskip

 \subsubsec{(4)}\proclaim{}Let $G$ be a~universally complete
 vector lattice equipped
 with the $f$-algebra multiplication uniquely determined by
 a choice of an order unity   in $G$. Also, let $E$ and $F$ be
 order dense ideals in $G$. Then, for every orthomorphism
 $\pi\in\Orth (E,F)$ there exists a~unique $g\in G$ such
 that $\pi x=g\cdot x$ for all~$x\in E$.
 \Endproc

 \subsec{1.1.7.}~An order bounded band preserving operator
 $\pi:D\to E$ is  a~{\it weak orthomorphism\/} of
 $E$ provided that $D$ is an order dense sublattice of $E$. In general,
 the weak orthomorphisms of $E$ do not comprise
 a good algebraic structure, while they do in the case of
 semiprime $f$-algebra. Denote by $\Orth^w(A)$ the set of
 all weak orthomorphisms with maximal domain. The set
 $\Orth^w(A)$ endowed with pointwise operations and ordering
 is an $f$-algebra (cp.~\cite{Wic2} for details).

 Denote by $Q(A)$ the \textit{maximal\/}
 (or \textit{complete}) ring of quotients of an $f$\hbox{-}algebra~$A$
 (cp.~\cite{Lam} for the definition). As was shown in
 \cite{Pag2}, $\Orth^\infty(A)$ and $Q(A)$ are not isomorphic.
 Nevertheless, $\Orth^\infty(A)$ can be embedded in $Q(A)$
 as an $f$\hbox{-}subalgebra \hbox{\cite{Pag2, Wic2}}. The following
 description of the maximal ring of quotients
 for an~(Archimedean) semiprime $f$-algebra is due to
 Wickstead \cite{Wic2}.

 \subsec{1.1.8.} \Theorem{}Let $A$ be a semiprime $f$-algebra. Then

 \subsubsec{(1)}~$\Orth^w(A)$ is a von Neumann regular
 $f$-algebra with unity $I_A$;

 \subsubsec{(2)}~$\Orth^\infty(A)$ is an $f$-subalgebra
 of~$\Orth^w(A)$;

 \subsubsec{(3)}~The maximal ring of quotients $Q(A)$
 coincides with $\Orth^w(A)$.

 \smallskip
 If, in addition, $A$ is relatively uniformly complete then
 \smallskip

 \subsubsec{(4)}~$Q(A)=\Orth^\infty(A)=\Orth^w(A)$.
 \Endproc

 \section{1.2. A Local Hamel Basis}

  Following \cite{MW}, we show  in this section that a
 universally complete vector lattice is locally one-dimensional
 if and only if all band preserving operators in it are
 automatically order bounded.

 \subsec{1.2.1.}~Let $G$ be an~arbitrary universally
 complete vector lattice with a fixed order unity    $\mathbbm{1}$.
 We~introduce some multiplication in $G$ that makes~$G$
 into a~commutative ordered algebra with unity~$\mathbbm{1}$.
 A~subset $\mathcal E\subset G$ is said to be {\it locally
 linearly independent\/} if whenever $e_1,\ldots,e_n\in\mathcal E$,
 $\lambda_1,\ldots,\lambda_n\in\mathbb R$, and $\pi$ is a band
 projection in~$G$ with $\pi (\lambda_1e_1+\cdots+\lambda_ne_n)=0$
 and $\pi e_1,\dots,\pi e_n$ nonzero and pairwise distinct
 we have $\lambda_k=0$ for all $k\assign 1,\ldots,n$.
 In other words, $\mathcal E$ is locally linearly independent
 if $\pi(\mathcal{E})\setminus\{0\}$, the set of all nonzero
 projections $\pi e$ of the elements $e\in\mathcal E$,
 is linearly independent for each nonzero $\pi\in\mathfrak{P}(G)$. A~maximal locally linearly independent set
 in~$G$ is  a~{\it local Hamel basis} for~$G$.

 \Proclaim{}There exists a~local Hamel basis for each universally %/?
 complete vector lattice.
 \Endproc

 \beginproof~Apply the Kuratowski--Zorn Lemma to the inclusion-ordered
 set of all locally linearly independent sets  in~$G$.~\endproof

 \Proclaim{1.2.2.}A~locally linearly independent set $\mathcal E$ in $G$ is
 a~local Hamel basis for $G$ if and only if for every
 $x\in G$ there exists a~partition of unity
 $(\pi_\xi)_{\xi\in\Xi}$ in $\goth P(G)$ such that
 $\pi_\xi x$ is a~finite linear combination of nonzero elements of
 $\pi_\xi\mathcal E$ for each $\xi\in\Xi$. Such representation
 of $\pi_\xi x$ is unique in the band $\pi_\xi (G)$.
 \Endproc

 \beginproof~$\la$:~Assume that  %/?
 $\mathcal E\subset G$ is locally linearly independent but is
 not a Hamel basis. Then we may find $x\in E$ such that
 $\mathcal E\cup\{x\}$ is locally linearly independent. Therefore,
 there is no nonzero band projection $\pi$ for which
 $\pi x$ is a linear combination of nonzero elements from
 $\pi\mathcal E$. This contradicts the existence of a~partition of
 unity with the mentioned properties.

 $\ra$:~If $\mathcal E$ is a local Hamel basis
 for~$G$ then $\mathcal E\cup\{x\}$ is not locally
 linearly independent for an arbitrary $x\in G$. Thus,
 there exist a nonzero band projection~$\pi$, reals
 $\lambda_0, \lambda_1,\ldots,\lambda_n\in\mathbb R$, and
 elements $e_1,\ldots,e_n\in\mathcal E$ such that
 $\pi(\lambda_0x+\lambda_1e_1+\cdots+\lambda_ne_n)=\nobreak 0$, while
 $\pi e_1,\ldots,\pi e_n$ are nonzero and pairwise distinct
 and not all $\lambda_0, \lambda_1,\ldots,\lambda_n$
 are equal to zero. Since the equality $\lambda_0=0$
 contradicts the~local linear independence of $\mathcal E$, it
 should be $\lambda_0\ne 0$, so that $\pi x$ is representable
 as a linear combination of $\pi e_1,\ldots,\pi e_n$. Now,
 the existence of the required partition of unity follows
 from the exhaustion principle.

 \subsec{1.2.3.}~Proposition 1.2.2 admits the following
 reformulation: A~locally linearly independent set $\mathcal E$
 in~$G$ is a
 local Hamel basis if and only if for every $x\in G$ there
 exist a partition of unity $(\pi_\xi)_{\xi\in\Xi}$ in
 $\goth P(G)$ and a family of reals
 $(\lambda_{\xi,e})_{\xi\in\Xi,e\in\mathcal E}$ such that
 $$
 x=\osum_{\quad\xi\in\Xi}\left(\sum_{e\in\mathcal E}\lambda _{\xi ,e}\pi_\xi
 e\right)\!,
 $$
 where  $\{e\in\mathcal E:\,\lambda_{\xi,e}\ne 0 \}$ is
 finite for every $\xi \in\Xi$. Moreover, the representation is
 unique in the sense that if $x$ admits one more representation
 $$
 x=\osum_{\quad\omega\in\Omega}\left(\sum_{e\in\mathcal E}
 \varkappa_{\omega,e}\rho_\omega e\right)\!,
 $$
 then for all $\xi\in\Xi $, $\omega\in\Omega$,  and $e\in\mathcal E$
 the relation $\pi_\xi\rho_\omega e\ne 0$ implies
 $\lambda_{\xi,e}=\varkappa_{\omega,e}$.

 \subsec{1.2.4.}~An~element $e\in G_{+}$ is  {\it locally
 constant\/}
 %\subject{locally constant element}
 with respect to $f\in\nomathbreak G_{+}$ if
 $e=\sup_{\xi \in \Xi }\lambda _{\xi}\pi _{\xi }f$ for
 some numeric family $(\lambda _{\xi })_{\xi\in\Xi}$ and
 a~family $(\pi _{\xi })_{\xi \in \Xi }$ of pairwise
 disjoint band projections.

 \Proclaim{} For each universally complete vector lattice $G$ the
 following are equivalent:

 \subsubsec{(1)}~All elements of~$G_{+}$ are locally constant with respect to $\mathbbm{1}$;

 \subsubsec{(2)}~All elements of~$G_{+}$ are locally constant with respect
 to an arbitrary order unity    $e\in G$;

 \subsubsec{(3)}~$\{\mathbbm{1}\}$ is a~local Hamel basis for $G$;

 \subsubsec{(4)}~Every local Hamel basis for $G$ consists of pairwise disjoint members.
 \Endproc

 \beginproof~Obviously, (2)~$\ra$~(1). To prove the converse,
 note that, given $x\in G$, we may choose a~partition of
 unity $(\pi _\xi )_{\xi \in\Xi}$ such that for each $\xi \in\Xi$
 both $\pi_\xi x$ and $\pi _\xi e$ are multiples of
 $\pi_\xi\mathbbm{1}$. So, $\pi_\xi x$ is a~multiple of $\pi _\xi e$. A~similar
 argument shows that $\{\mathbbm{1}\}$ is a~local Hamel basis if and
 only if so is $\{f\}$ for every order unity    $f\in G$. Thus, if (4)
 holds and $\mathcal E$ is a~local Hamel basis for $G$ then
 $f\assign \sup\{e:\,e\in\mathcal E\}$ exists and $\{f\}$ is a~local
 Hamel basis for~$G$. It follows that (4)~$\ra$~(3).
 Clearly, (3)~$\ra$~(1) by 1.2.3. To complete the
 proof, we had to show (1)~$\ra$~(4). If (4) fails
 then we may choose a nonzero band projection $\pi $ and
 a~local Hamel basis containing two members $e_1$ and $e_2$
 such that both $\pi e_1$ and $\pi e_2$ are nonzero
 multiples of $\pi\mathbbm{1}$. Consequently,
 $\pi(\lambda_1 e_1+\lambda _2 e_2)=0$ for some $\lambda_1,\lambda
 _2\in\mathbbm{R}$ and we arrive at the contradictory conclusion that
 $ \{e_1,e_2\}$ is not locally linearly independent.~\endproof

 A~universally complete vector lattice $G$ is  {\it locally one-dimensional\/}
 if $G$ satisfies the equivalent conditions (1)--(4) of the above
 proposition.

 \subsec{1.2.5.} \theorem{}Let~$G$ be a~universally complete vector lattice.
 Then the following are equivalent:

 \subsubsec{(1)}~$G$ is locally one-dimensional;

 \subsubsec{(2)}~Every band preserving operator $T:G\to G$ is
 order  bounded.
 \Endproc

 \beginproof~(1)~$\ra$~(2):~Recall that
 a~linear operator $T:G\to G$ is band preserving if and
 only if $\pi T=T\pi$ for every band projection $\pi$ in
 $G$ \big(cp.~1.1.1\,(4)\big). Assume that $T$ is band preserving and put
 $\rho\assign T\mathbbm{1}$. Since an~arbitrary $e\in G_+$ can
 be expressed as
 $e=\sup_{\xi\in\Xi}\lambda_\xi\pi_\xi\mathbbm{1}$, we deduce
 $$
 \pi_\xi Te=
 T(\pi_\xi e)=T(\lambda _\xi \pi_\xi\mathbbm{1})=
 \lambda_\xi\pi_\xi T(\mathbbm{1})=\pi_\xi(e)T(\mathbbm{1})=
 \pi_\xi e\rho,
 $$
 whence $Te=\rho e$. It follows that $T$ is a~multiplication
 operator in $G$ which is obviously order  bounded.

 (2)~$\ra$~(1):~Assume that (1) is false.
 According to 1.2.4\,(4) there is a~local Hamel
 basis $\mathcal E$ for $G$ containing two members $e_1$ and $e_2$
 that are not disjoint. Then the band projection
 $\pi\assign [e_1]\wedge [e_2]$ is nonzero.
 (Here and below $[e]$ is the band projection onto $\{e\}^{\perp\perp}$.)
 For an arbitrary $x\in G$
 there exists a partition of unity $(\pi_\xi)_{\xi \in\Xi }$ such
 that $\pi_\xi x$ is a finite linear combination of elements of~$\mathcal E$. Assume the elements of $\mathcal E$ have been labelled so
 that $\pi_\xi x=\lambda_1 \pi _\xi e_1+\lambda_2 \pi _\xi
 e_2+\cdots$. Define $Tx$ to be a~unique element in $G$ with
 $\pi_\xi Tx\assign \lambda_1 \pi \pi _\xi e_2$. It is easy to check
 that $T$ is a~well defined linear operator from $G$ into itself.

 Take $x,y\in G$ with $x\perp y$ and let $(\pi_\xi)_{\xi \in\Xi }$
 be a~partition of unity such that both $\pi_\xi x$ and $\pi_\xi y$
 are finite linear combination of elements from $\mathcal E$. Refining
 the partition of unity if need be, we may also require that at
 least one of the elements $\pi_\xi x$ and $\pi_\xi y$ equals
 zero for all $\xi \in\Xi $. If $\pi_\xi y\ne 0$ then $\pi_\xi
 x=0$,  and so the corresponding $\lambda_1e_1$ is equal to zero. If $\pi \pi _\xi \ne
 0$ then $\lambda_1=0$, and in any case $\pi_\xi Tx=0$. It follows
 that $Tx\perp y$ and $T$ is band preserving.
 If $T$ were order  bounded then $T$ would be presentable as $Tx=ax$
 $(x\in G)$ for some $a\in G$, see 1.1.6\,(4). In particular, $Te_2=ae_2$ and,
 since $Te_2=0$ by definition, we have $0=[e_2]|a|\ge\pi |a|$. Thus
 $\pi e_2=T(\pi e_1)=a\pi e_1=0$, contradicting the definition of~$\pi $.~\endproof

 \section{1.3. $\sigma$-Distributive Boolean Algebras}

 In this section we present the main result
 of~\cite{Gut6}: {\sl A~universally complete vector lattice~$G$ is
 locally one-dimensional if and only if the base of~$G$ is
 $\sigma$-distributive}.

 \subsec{1.3.1.} A~$\sigma $-complete Boolean algebra
 $\mathbbm{B}$ is said to be {\it $\sigma$-distributive\/}
 if $\mathbbm{B}$ satisfies one of the~following equivalent conditions
 (cp.~\cite[19.1]{Sik}):
 \smallskip

 \subsubsec{(1)} $\bigwedge_{n\in\mathbbm{N}}
 \bigvee_{m\in\mathbbm{N}} b_{m}^{n}=
 \bigvee _{m\in\mathbbm{N}^{\mathbbm{N}}}
 \bigwedge_{n\in{\mathbbm{N}}} b_{m(n)}^{n}$
 for all
 $(b_{m}^{n})_{n,m\in\mathbbm{N}}$ in~$\mathbbm{B}$;
 \smallskip

 \subsubsec{(2)} $\bigvee _{n\in\mathbbm{N}} \bigwedge _{m\in\mathbbm{N}} b_{m}^{n}=
 \bigwedge _{m\in\mathbbm{N}^{\mathbbm{N}}}
 \bigvee _{n\in\mathbbm{N}} b_{m(n)}^{n}$
 for all
 $(b_{m}^{n})_{n,m\in\mathbbm{N}}$ in~$\mathbbm{B}$;
 \smallskip

 \subsubsec{(3)} $\bigvee _{\varepsilon \in \{1,-1\}^{\mathbbm{N}}}
 \bigwedge _{n\in\mathbbm{N}} \varepsilon (n)b_{n}=\mathbbm{1}$
 for all $(b_{n})_{n\in\mathbbm{N}}$ in~$\mathbbm{B}$.
 \smallskip

 \noindent
 \big(Here $1b_{n}\assign b_{n}$ and~$(-1)b_{n}$ is the~complement
 of~$b_{n}$.\big)

 \subsec{1.3.2.}~Let~$\mathbb{B}$ be an~arbitrary Boolean algebra.
 A~subset of~$\mathbb{B}$ with supremum unity is called a~{\it cover\/} of~$\mathbb{B}$.
 Partitions of unity in~$\mathbb{B}$ are referred to as partitions of~$\mathbb{B}$ for brevity.
 Let~$C$ be a~cover of~$\mathbb{B}$. A~subset~$C_{0}$
 of~$\mathbb{B}$ is said to be {\it refined\/} from~$C$ if,
 for each $c_{0}\in C_{0}$, there exists $c\in C$ such
 that $c_{0}\leq c$. An~element $b\in\mathbb{B}$ is
 {\it refined\/} from~$C$ provided that  $\{b\}$ is refined
 from~$C$; i.e., $b\leq c$ for some  $c\in C$.
 If~$(C_{n})_{n\in {\mathbb N}}$ is a~sequence of covers
 of~$\mathbb{B}$ and $b\in\mathbb{B}$ is refined from each of
 the~covers~$C_{n}$ ($n\in {\mathbb N}$), then we say
 that~$b$ is refined from $(C_{n})_{n\in {\mathbb N}}$.
 We~also refer to a~cover whose all elements are refined
 from $(C_{n})_{n\in {\mathbbm N}}$ as refined from
 the~sequence.

 \subsec{1.3.3.}\proclaim{}Let~$\mathbbm B$ be
 a~$\sigma$-complete Boolean algebra.
 The~following  are equivalent:

 \subsubsec{(1)}~$\mathbbm B$ is $\sigma $-distributive;

 \subsubsec{(2)}~There is a~(possibly, uncountable) cover refined from
each sequence of countable covers of\/~$\mathbbm B$;

 \subsubsec{(3)}~There is a~(possibly, infinite) cover refined
 from each sequence of finite covers of\/~$\mathbbm B$;

 \subsubsec{(4)} There is a cover refined from each sequence of two-element partitions of\/~$B$.
 \Endproc

 \beginproof~A~proof of  (1)$\iff $(2)
 can be found in~\cite[19.3]{Sik}). Item~(4) is
 a~paraphrase of 1.3.1\,(3) in the~definition of
 $\sigma $-distributivity. The~implications
 (2)$\ra $(3)$\ra $(4) are
 obvious.~\endproof

 \subsec{1.3.4.} \proclaim{}Let~$\mathbbm B$ be a~complete Boolean algebra. The~following
 are equivalent:

 \subsubsec{(1)}~$\mathbbm B$ is $\sigma $-distributive;

 \subsubsec{(2)}
  There is a~(possibly, uncountable) partition refined
 from each sequence of countable partitions of\/~$\mathbbm B$;

 \subsubsec{(3)}
 There is a~(possibly, infinite) partition refined from each
 sequence of finite partitions of\/~$\mathbbm B$;

 \subsubsec{(4)}
 There is a~partition refined from each sequence
 of two-element partitions of\/~$\mathbbm B$.
 \Endproc

 \beginproof~The~claim follows from 1.3.3 in view of
 the~exhaustion principle.~\endproof

 \subsec{1.3.5.} Let $Q$ stand for the Stone space of
 $\mathbb{B}$ and denote by ${\rm Clop}(Q)$ the Boolean
 algebra of all clopen sets in $Q$. We~say that a~function
 $g\in C_{\infty }(Q)$ is {\it refined from a~cover}
 %\subject{refined function from a~cover}
 $C$ of the Boolean algebra ${\rm Clop}(Q)$ if, for every two
 points $q^{\prime },q^{\prime\prime }\in Q$ satisfying
 the~equality $g(q^{\prime })=g(q^{\prime\prime })$, there exists
 an~element $U\in C$ such that $q^{\prime },q^{\prime\prime }\in
 U$. If~$(C_{n})_{n\in {\mathbb N}}$ is a~sequence of covers of
  ${\rm Clop}(Q)$ and a~function~$g$ is refined from
 each of the~covers~$C_{n}$ ($n\in {\mathbb N}$), then we say that~$g$ is {\it refined from}
 %\subject{refined function from a~sequence of covers}
 $(C_{n})_{n\in {\mathbb N}}$.

 \subsec{1.3.6.} \lemma{}
 For every sequence of finite covers of\/~${\rm Clop}(Q)$,
 there is  a~function of~$C(Q)$ refined from the~sequence.
 \Endproc

 \beginproof~Let~$(C_{n})_{n\in {\mathbb N}}$ be a~sequence of finite
 covers of~${\rm Clop}(Q)$. By induction, it is easy to construct a~sequence of
 partitions $P_{m}=\{U_{1}^{m},U_{2}^{m},\ldots,U_{2^{m}}^{m}\}$ of
 ${\rm Clop}(Q)$ with the~following properties:

 (1) for every $n\in {\mathbb N}$, there is  $m\in {\mathbb N}$
 such that the~partition~$P_{m}$ is refined from~$C_{n}$;

 (2) $U_{j}^{m}=U_{2j-1}^{m+1}\lor U_{2j}^{m+1}$ for all
 $m\in {\mathbb N}$ and $j\in \{1,2,\ldots,2^{m}\}$.

 Given $m\in {\mathbb N}$, define the~two valued function
 $\chi _{m}\in C(Q)$ as follows:
 $$
 \chi _{m}\assign \sum _{i=1}^{2^{m-1}}\chi (U_{2i}^{m}),
 $$
 where~$\chi (U)$ is the~characteristic function of
 $U\subset Q$. Since the~series $\sum _{m=1}^{\infty }
 \frac1{3^{m}}\chi _{m}$ is uniformly convergent, its sum~$g$
 belongs to~$C(Q)$. We~will show that~$g$ is refined
 from $(C_{n})_{n\in {\mathbb N}}$. By property~(1) of
 the~sequence $(P_{m})_{m\in {\mathbb N}}$, it  suffices
 to establish that~$g$ is refined from $(P_{m})_{m\in
 {\mathbb N}}$.

 Assume the~contrary and consider the~least $m\in {\mathbb
 N}$ such that $g$ is not refined from~$P_{m}$. In~this case,
 there are two points
 $q^{\prime },q^{\prime\prime }\in Q$ satisfying the~equality
 $g(q^{\prime })=g(q^{\prime\prime })$ and belonging to distinct
 elements of~$P_{m}$. Since~$g$ is refined from
 $P_{m-1}$ (for $m>1$), from property~(2) of
 the~sequence $(P_{m})_{m\in {\mathbb N}}$ it follows that
 $q^{\prime }$ and $q^{\prime\prime }$ belong to some adjacent elements
 of~$P_{m}$; i.e., elements of the~form~$U_{j}^{m}$ and~$U_{j+1}^{m}$,
 with $j\in \{1,\ldots,2^{m}-1\}$. For~definiteness,
 suppose that~$q^{\prime }$ belongs to an~element with an even
 index and~$q^{\prime\prime }$, to that with an odd index; i.e., $\chi
 _{m}(q^{\prime })=1$ and $\chi _{m}(q^{\prime\prime })=0$.
 Since $\chi _{i}(q^{\prime
 })=\chi _{i}(q^{\prime\prime })$ for all $i\in \{1,\ldots,m-1\}$;
therefore, we have:
 $$
 g(q^{\prime })-g(q^{\prime\prime })=\frac1{3^{m}}+ \sum
 _{i=m+1}^{\infty }\frac1{3^{i}}\big(\chi _{i}(q^{\prime
 })-\chi _{i}(q^{\prime\prime })\big)\geq \frac1{3^{m}}-\sum _{i=m+1}^{\infty }
 \frac1{3^{i}}=\frac1{2\cdot 3^{m}}>0,
 $$
 which contradicts the~equality $g(q^{\prime })=g(q^{\prime\prime})$.~\endproof

 \subsec{1.3.7.} \theorem{}A~universally complete vector lattice $G$
 is locally one-dimensional if and only if the base of~$G$ is
 $\sigma$-distributive.
 \Endproc

 \beginproof~Let~$Q$ be the~Stone  space of the base of~$G$. Suppose that~$G$
 is locally one-dimensional and consider an~arbitrary sequence
 $(P_{n})_{n\in {\mathbb N}}$ of finite partitions of~${\rm Clop}(Q)$.
 By~1.3.4, to prove the $\sigma $-distributivity of~$G$, it
  suffices to refine a~cover of ${\rm Clop}(Q)$ from
 $(P_{n})_{n\in {\mathbb N}}$. By Lemma~1.3.6, we may
 refine  $g\in C_{\infty }(Q)$ from the~sequence
 $(P_{n})_{n\in {\mathbb N}}$. Since~$G$ is locally one-dimensional,
 there exists a~partition $(U_{\xi })_{\xi \in \Xi }$ of
 ${\rm Clop}(Q)$ such that~$g$ is
 constant on each of the~sets~$U_{\xi }$. Show that
 $(U_{\xi })_{\xi \in \Xi }$ is refined from
 $(P_{n})_{n\in {\mathbb N}}$. To~this end,  fix arbitrary indices
 $\xi \in \Xi $ and $n\in {\mathbb N}$ and establish that
 $U_{\xi }$ is refined from~$P_{n}$. We~may assume
 that $U_{\xi }\ne \varnothing $. Let~$q_{0}$ be an~element of~$U_{\xi }$.
 Finiteness of~$P_{n}$ allows us to find
 an~element~$U$ of $P_{n}$ such that $q_{0}\in U$. It~remains to
 observe that $U_{\xi }\subset U$. Indeed, if $q\in U_{\xi }$ then
 $g(q)=g(q_{0})$ and, since~$g$ is refined from
 $P_{n}$, the~points~$q$ and~$q_{0}$ belong to the~same element
 of~$P_{n}$; i.e., $q\in U$.

 Assume now that the~base of~$G$ is $\sigma $-distributive and
 consider an~arbitrary $g\in C_{\infty }(Q)$. By
 the~definition of a~locally one-dimensional
 vector lattice, it suffices to construct a~partition
 $(U_{\xi })_{\xi \in \Xi }$ of~${\rm Clop}(Q)$ such that
 $g$ is constant on each of the~sets~$U_{\xi }$.
 Given a~natural~$n$ and integer~$m$, denote by~$U_{m}^{n}$
 the~interior of the~closure of the~set of all points
 $q\in Q$ for which $\frac mn\leq g(q)<\frac{m+1}n$
 and put $P_{n}\assign \big\{U_{m}^{n}:m\in\nomathbreak {\mathbb
 Z}\big\}$. By~1.3.4, from the~sequence
 $(P_{n})_{n\in {\mathbb N}}$ of countable partitions of
 ${\rm Clop}(Q)$, we may refine some partition $(U_{\xi })_{\xi
 \in \Xi }$. It~is easy that  this is a~desired partition.~\endproof

 \subsec{1.3.8.}~\theorem{}There exists a~purely nonatomic
 locally one-dimensional universally complete vector lattice.
 \Endproc

 \beginproof~Theorem 1.3.7 reduces the problem to the existence of
 a~purely nonatomic $\sigma $\hbox{-}distributive complete Boolean
 algebra. An~algebra of this kind is constructed below in
 1.3.9--1.3.11.~\endproof

 \subsec{1.3.9.} A~Boolean algebra~$\mathbb{B}$ is
 {\it $\sigma $-inductive\/} provided that each decreasing
 sequence of nonzero elements of~$\mathbb{B}$ has
 a~nonzero lower bound. A~subalgebra~$\mathbb{B}_{0}$
 of~$\mathbb{B}$ is  {\it dense\/} if, for every nonzero
 element $b\in\mathbb{B}$, there exists a~nonzero
 element $b_{0}\in\mathbb{B}_{0}$ such that $b_{0}\leq b$.

 \Lemma{}If~a $\sigma $-complete Boolean algebra $\mathbb{B}$
 has a~$\sigma $-inductive dense subalgebra then $\mathbb{B}$ is $\sigma
 $-distributive. \Endproc

 \beginproof~Let~$\mathbb{B}_{0}$ be a~$\sigma $-inductive dense subalgebra of~$\mathbb{B}$.
 Consider an~arbitrary sequence $(C_{n})_{n\in {\mathbb N}}$ of
 countable covers of~$\mathbb{B}$, denote by~$C$ the~set of all elements
 in~$\mathbb{B}$ that are refined from $(C_{n})_{n\in {\mathbb N}}$, and
 assume by way of contradiction that~$C$ is not a~cover of~$\mathbb{B}$.
 Then there is a~nonzero element $b\in \mathbb{B}$ disjoint from
 all elements of~$C$.

 By~induction, we construct the sequences $(b_{n})_{n\in {\mathbb N}}$ and
 $(c_{n})_{n\in {\mathbb N}}$ as follows: Let~$c_{1}$ be an~element
 of~$C_{1}$ such that $b\land c_{1}\ne 0$. Since~$\mathbb{B}_{0}$ is
 dense, there is an~element $b_{1}\in \mathbb{B}_{0}$ such that
 $0<b_{1}\leq b\land c_{1}$. Suppose that~$b_{n}$
 and~$c_{n}$ are already constructed. Let~$c_{n+1}$ be
 an~element of~$C_{n+1}$ such that $b_{n}\land c_{n+1}\ne 0$. As
 $b_{n+1}$ we take an~arbitrary element of~$\mathbb{B}_{0}$ such that
 $0<b_{n+1}\leq b_{n}\land c_{n+1}$.

 Thus, we have constructed sequences $(b_{n})_{n\in {\mathbb N}}$ and
 $(c_{n})_{n\in {\mathbb N}}$ such that $b_{n}\in \mathbb{B}_{0}$,
 $b_{n}\leq c_{n}\in C_{n}$ and $0<b_{n+1}\leq
 b_{n}\leq b$ for all $n\in {\mathbb N}$. Since~$\mathbb{B}_{0}$
 is $\sigma $-inductive, $\mathbb{B}_{0}$ contains a~nonzero
 element~$b_{0}$
 that satisfies $b_{0}\leq b_{n}$ for all $n\in {\mathbb N}$.
 By the~inequalities $b_{0}\leq c_{n}$, we see that~$b_{0}$ is
 refined from $(C_{n})_{n\in {\mathbb N}}$; i.e., $b_{0}$ belongs
 to~$C$. On~the other hand, $b_{0}\leq b$, which contradicts
 the disjointness of~$b$ from all elements of~$C$.~\endproof

 \subsec{1.3.10.}~As~is well known, to every Boolean algebra~$\mathbb{B}$ there
is a~complete Boolean algebra~$\overline{\mathbb{B}}$ including%
~$\mathbb{B}$ as a~dense subalgebra (cp.~\cite[Section~35]{Sik}).
 This~$\overline{\mathbb{B}}$ is unique to within
 an~isomorphism and called a~{\it completion\/} of~$\mathbb{B}$.
 Obviously, a~completion of a~purely nonatomic Boolean
 algebra is purely nonatomic. Moreover, by
 Lemma~1.3.9, a~completion of a~$\sigma $-inductive
 algebra is $\sigma $-distributive. Therefore, in order to
 prove existence of a~purely nonatomic $\sigma $-distributive
 complete Boolean algebra, it  suffices to exhibit
 an~arbitrary purely nonatomic $\sigma $-inductive Boolean
 algebra. The examples of these algebras are readily available.
 For~the sake of completeness, we present here one of
 the~simplest constructions.

 \subsec{1.3.11.}\proclaim{}Let~$\mathbb{B}$ be the~boolean
 of\/~${\mathbb N}$ and let~$I$ be
 the~ideal of\/~$\mathbb{B}$ comprising all finite subsets of\/~${\mathbb N}$.
 Then the~quotient algebra~$\mathbb{B}/I$
 (cp.~\cite[Section~10]{Sik}) is purely nonatomic and
 $\sigma$-inductive.
 \endproc

 \beginproof~The pure nonatomicity of~$\mathbb{B}/I$ is
 obvious. In~order to prove that $\mathbb{B}/I$ is
 $\sigma $\hbox{-}inductive, it  suffices to consider an~arbitrary
 decreasing sequence $(b_{n})_{n\in {\mathbb N}}$ of infinite
 subsets of~${\mathbb N}$ and construct an~infinite subset~$b\subset {\mathbb N}$ such that the~difference $b\backslash b_{n}$
 is finite for each $n\in {\mathbb N}$. We~can easily obtain
 the~desired set $b\assign \{m_{n}:n\in {\mathbb N}\}$ by
 induction, letting $m_{1}\assign \min b_{1}$ and
 $m_{n+1}\assign \min\{m\in b_{n+1}: m>m_{n}\}$.~\endproof

 \section{PART 2. BOOLEAN APPROACH}

 The purpose of the this part is to present the approach
of  Boolean valued analysis to the Wickstead problem and
 prove that if $G$ is a~universally complete vector lattice and $\mathbbm{B}\assign \mathfrak{P}(G)$ is the base of~$G$
 then the following  are equivalent:
 \smallskip

 \subsubsec{WP(1)}~$\mathbbm{B}$ is $\sigma$-distributive;
 \smallskip

 \subsubsec{WP(2)}~$\mathcal{R}=\mathbbm{R}^{\scriptscriptstyle\wedge}$
 inside  $\mathbbm{V}^{(\mathbbm{B})}$;
 \smallskip

 \subsubsec{WP(3)}~$G$ is locally one-dimensional;
 \smallskip

 \subsubsec{WP(4)}~Every band preserving linear operator
 in $G$ is order bounded.
 \smallskip

 In Sections 2.1, 2.2, and 2.3 we will give purely
 Boolean valued proofs of the equivalences
 WP(2)~$\iff$~WP(4), WP(2)~$\iff$~WP(3),
 and WP(1)~$\iff$~WP(2), respectively. It turns
 out that all these equivalences reduce to
 some simple properties of reals and cardinals in
 an appropriate Boolean valued model (cp.~\cite{Kus1}).

 Throughout this part $\mathcal{R}$ denotes the Boolean valued
 reals and ${\mathbbm R}^{\scriptscriptstyle\wedge}$ is
 considered as a~dense subfield of~$\mathcal R$. More
 precisely, ${\mathcal R}\in{\mathbbm V}^{({\mathbbm B})}$
 and $[\![\,{\mathcal R}$ is a~field of reals$\,]\!]={\mathbbm 1}$,
 while $[\![\,{\mathbbm R}^{\scriptscriptstyle\wedge}$ is
 a~dense subfield
 of~$\mathcal R\,]\!]={\mathbbm 1}$ (cp.~A3.3 and A3.4).
 The Gordon Theorem A3.6 says that if $G$ is a universally
 complete vector lattice and $\mathbb{B}\assign \mathfrak{P}(G)$,
 then $\mathcal R{\downarrow}$ is a universally complete
 vector lattice isomorphic to $G$.

 \section{2.1. Representation of a Band Preserving Operator}

 In this section we show that the equivalence
 WP(2)~$\iff$~WP(4) is immediate from the
 Boolean valued representation of band preserving operators.

 \subsec{2.1.1.}~Throughout the~section
 we let $G$
 stand for  the~universally complete vector lattice
 $\mathcal R{\downarrow}$.
 Recall that $G$~is a faithful  $f$-ring with unity
 $\mathbbm{1}\assign 1^{\scriptscriptstyle\wedge}$.

 Let $\End_N(G)$~be the set of all band preserving endomorphisms
 of~$G$. Clearly, $\End_N(G)$ is a~vector space. Moreover,
 $\End_N(G)$ becomes a faithful unitary module over $G$ on letting
 $gT$ be equal to $gT:x\mapsto g\cdot Tx$ for all $x\in G$. This is
 immediate since the multiplication by an element of $G$ is band
 preserving and the composite of band preserving operators is band
 preserving too. By
 $\End_{\mathbb R^{\scriptscriptstyle\wedge}}(\mathcal R)$
 we denote the element of $\mathbb V^{(\mathbb B)}$ that represents
 the space of all $\mathbb R^{\scriptscriptstyle\wedge}$-linear
 operators from  $\mathcal R$ into $\mathcal R$. Then $\End_{\mathbb R^{\scriptscriptstyle\wedge}}(\mathcal R)$ is a~vector
  space over $\mathcal R$ inside $\mathbb V^{(\mathbb B)}$, and
 $\End_{\mathbb R^{\scriptscriptstyle\wedge}}(\mathcal R){\downarrow}$ is
 a~faithful unitary module over~$G$.

 \proclaim{2.1.2.}A linear operator~$T$ on a universally
 complete vector lattice $G$ is band pre\-ser\-ving if
 and only if $T$ is extensional.
 \Endproc

 \beginproof~By the Gordon Theorem A3.6 and A2.4\,(7),
  $T:G\to G$ is extensional if and only if, for all
  $x\in G$ and $\pi \in\goth P(G)$,
 from $\pi x=0$ it follows that $\pi Tx=0$.
 By taking $x\assign \pi^\perp y$ we conclude that
 $\pi T\pi^\perp=0$ or, in other words, $\pi T=\pi T\pi$.
 Substituting $\pi^\perp$ for~$\pi$, we see that $T\pi =\pi T\pi$,
 and so $\pi T=T\pi$.
 Hence, $T$ is band preserving by 1.1.1\,(4).
 Conversely, for a band preserving $T$ we see that
 $\pi x=0$ implies $\pi Tx=0$ by definition.~\endproof

 \subsec{2.1.3.}~If $\sigma\in\mathbb{V}^{(\mathbb{B})}$ and
 $[\![\,\sigma:\mathcal{R}\to\mathcal{R}\,]\!]=\mathbbm{1}$, then
 there exists a unique map
 $S:\mathcal{R}{\downarrow}\to\mathcal{R}{\downarrow}$ such that
 $$
 [\![S(x)=\sigma(x)]\!]=\mathbbm{1}\quad (x\in\mathcal{R}{\downarrow}).
 $$
 This map $S$ is called the \textit{descent} of $\sigma$
 and is denoted by $\sigma{\downarrow}$. It is of importance
 that the descent is \textit{extensional} (cp.~A2.5):
 $$
 [\![x=y]\!]\leq[\![S(x)=S(y)]\!]\quad(x,y\in\mathcal{R}{\downarrow}).
 $$
 It is immediate from A3.6 that $S$ is extensional if and only
 if $bx=by$ implies $bS(x)=bS(y)$
 for all $x,y\in\mathcal{R}{\downarrow}$ and $b\in\mathbb{B}=\mathfrak{P}(\mathcal{R}{\downarrow})$.

 Conversely, given an extensional map
 $S:\mathcal{R}{\downarrow}\to\mathcal{R}{\downarrow}$, there
 exists a unique function $\sigma:\mathcal{R}\to\mathcal{R}$
 inside $\mathbb{V}^{(\mathbb{B})}$, such that $S=\sigma{\downarrow}$.
 We say that $\sigma$ is the \textit{ascent} of $S$ and write
 $\sigma=S{\uparrow}$ (cp.~A2.4). Thus, the descent and
 ascent carry out a bijection between the sets of all extensional
 mappings from $\mathcal{R}{\downarrow}$ into $\mathcal{R}{\downarrow}$
 and all elements $\sigma\in\mathbb{V}^{(\mathbb{B})}$
 with $[\![\,\sigma:\mathcal{R}\to\mathcal{R}\,]\!]=\mathbbm{1}$
 (cp. the Escher rules for arrow cancellations in A2.6). Denote the
 latter set by $F(\mathcal{R}){\downarrow}$.

 \subsec{2.1.4.}~Let $\Ext(\mathcal{R}{\downarrow})$ be the~set
 of all extensional mappings from $\mathcal{R}{\downarrow}$ into
 $\mathcal{R}{\downarrow}$. The pointwise operations make this
 set into a unital module over the ring $\mathcal{R}{\downarrow}$.
 The set $F(\mathcal{R}){\downarrow}$ can be endowed with a
 module structure over $\mathcal{R}{\downarrow}$ by analogy  to
 A3.5.

 \Proclaim{}The bijection in 2.1.3 is an isomorphism of the
 modules $\Ext(\mathcal{R}{\downarrow})$ and $F(\mathcal{R}){\downarrow}$.
 \Endproc

 \smallskip
 \beginproof~This is immediate from the following identities:
 $$
 \gathered
 (S+T){\uparrow}x=(S+T)x=Sx+Tx=S{\uparrow}x\oplus T{\uparrow}x=
 (S{\uparrow}\oplus T{\uparrow})x\quad(x\in\mathcal{R}{\downarrow});
 \\
 (\alpha\cdot S){\uparrow}x=(\alpha\cdot S)x=\alpha\cdot(Sx)=
 \alpha\odot (S{\uparrow}x)=
 (\alpha\odot S{\uparrow})x\quad(\alpha,x\in\mathcal{R}{\downarrow}),
 \endgathered
 $$
 where $\oplus$ and $\odot$ stand for the operations in
 $\mathcal{R}$ and $F(\mathcal{R})$, while $+$ and $\cdot$
 symbolize the operations in $\mathcal{R}{\downarrow}$ and
 $\Ext(\mathcal{R}{\downarrow})$.~$\endproof$

 \Proclaim{2.1.5.}The modules $\End_N(G)$ and $\End_{\mathbb R^{\scriptscriptstyle\wedge}}(\mathcal R){\downarrow}$
 are isomorphic. The isomorphy may be established by sending a band preserving operator to its ascent.
 \Endproc

 \beginproof~Each $T\in \End_N(G)$ is extensional by 2.1.2,
 and so $T$ has the ascent $\tau\assign T{\uparrow}$ presenting
 the unique mapping from $\mathcal R$ into $\mathcal R$
 such that $[\![\tau (x)=Tx]\!]$ for all $x\in G$
 (cp.~2.1.3). Using this identity and the definition
 of the ring structure on $\mathcal R{\downarrow}$, we see
 $$
 \gathered
 \tau (x\oplus y)=T(x+y)=Tx+Ty=\tau (x)\oplus\tau (y)\quad(x,y\in G),
 \\
 \tau (\lambda^{\scriptscriptstyle\wedge}\odot x)=T(\lambda\cdot x)=
 \lambda\cdot Tx=\lambda^{\scriptscriptstyle\wedge}\odot\tau (x)
 \quad(x\in G,\ \lambda \in\mathbb R).
 \endgathered
 $$
 Hence, $[\![\,\tau:\mathcal R\to\mathcal R$ is an~$\mathbb
 R^{\scriptscriptstyle\wedge}$-linear function\,$]\!]=\mathbbm1$;
 i.e., $[\![\,\tau\in\End_{\mathbb R^{\scriptscriptstyle\wedge}}(\mathcal R)\,]\!]= \mathbbm1$.
 If~%
 $\tau\in\End_{\mathbb R^{\scriptscriptstyle\wedge}}(\mathcal
 R){\downarrow}$ then the descent
 $\tau{\downarrow}:G\to G$ is extensional
 (cp.~2.1.3). The same arguments as above convince us that
 if $\tau$ is $\mathbb R^{\scriptscriptstyle\wedge}$-linear inside
 ${\mathbb V}^{(\mathbb B)}$ then
 $\tau{\downarrow}$ is a linear operator.
 By 2.1.2 \ $\tau{\downarrow}$ is band preserving.
 The claim results now from~2.1.4.~\endproof

 \subsec{2.1.6.}~In 2.1.5 we encountered the
 following situation: There is some ordered subfield
 $\mathbb P$ of the reals $\mathbb R$ that includes
 $\mathbb Q$. Consequently, $\mathbb R$ is a vector space over
  $\mathbb P$ and has a Hamel basis, say $\mathcal E$.
  Denote the set of all $\mathbb P$-linear functions in
  $\mathbb R$ by $\End_{\mathbb P}(\mathbb R)$.
 For the sake of completeness, we  recall the two well-known facts:

 \smallskip
 \subsubsec{(1)} \proclaim{}Let $\mathbb P$ be a~subfield of\/ $\mathbbm{R}$.
 The general form of a~$\mathbbm{P}$-linear function
 $f:\mathbb R\to\mathbb R$ is given as

 $$
 f(x)=\sum_{e\in\mathcal E} x_e\phi(e)\quad \text{if }\ x=\sum_{e\in\mathcal E} x_e e,
 $$
 where $\phi:\mathcal E\to\mathbb R$ is an~arbitrary function
 and the second formula  is the expansion of $x\in\mathbbm{R}$
 with respect to the Hamel basis $\mathcal E$ and the
 coefficients $(x_e)_{e\in\mathcal{E}}$ are such that
 $\{e\in\mathcal{E}:\,x_e\ne0\}$ is a~finite set.
 \Endproc

 \beginproof~This is immediate from the definition and properties
 of a Hamel basis.~\endproof

 \smallskip
 \subsubsec{(2)} \proclaim{}An arbitrary $\mathbb P$-linear function
 $f:\mathbb R\to\mathbb R$ admits the representation
 $f(x)=cx$ $(x\in\mathbb R)$ with some
 $c\in\mathbb R$ if and only if $f$ is bounded above or below
 on some interval $\mathopen{]}a,b\mathclose{[}\subset\mathbb R$, with $a<b$.
 \Endproc

 \beginproof~Necessity is obvious. To prove sufficiency, assume that
 $f$ is bounded above by a real $M$ on $\mathopen{]}a,b\mathclose{[}$. Then the open set
 $\{(s,t)\in\mathbb R^2:\,a<s<b,\,M<t\}$ is disjoint from the graph of
 $f$, and so the graph of~$f$ cannot be dense in $\mathbb R^2$.
 However, if $f$ fails to admit the desired representation then
 the graph of~$f$ is dense in~$\mathbb R^2$. This is established
 in much the same way as in the case of the Cauchy functional equation
 (cp.~\cite[Chapter~2, Theorem~3]{AD}).~\endproof

 \subsec{2.1.7.}~We now exhibit the two corollaries for
 band preserving operators which are  the Boolean valued interpretations
 of~2.1.6\,(1),\,(2).

 \smallskip
 \subsubsec{(1)}\proclaim{}A band preserving operator
 $T\in\End_N(G)$ is order bounded if and only if
 $T$ may be presented as $Tx=g\cdot x$
 $(x\in G)$ for some fixed $g\assign g_T\in G$.
 \Endproc

 \beginproof~It suffices to observe that the ascent functor
 preserves the property of order boundedness in~2.1.5
 and apply
 2.1.6\,(2) inside ${\mathbb V}^{(\mathbb B)}$.~\endproof

 \smallskip
 \subsubsec{(2)}\proclaim{}For every band preserving
 endomorphism of~$G\assign \mathcal R{\downarrow}$ to be order bounded
 it is necessary and sufficient that  ${\mathbb V}^{(\mathbb B)}\models\mathcal R=
 \mathbb R^{\scriptscriptstyle\wedge}.$
 \Endproc

 \beginproof~$\la$:~If $\mathbb R^{\scriptscriptstyle\wedge}$
 coincides with the reals $\mathcal R$ inside
 ${\mathbb V}^{(\mathbb B)}$ then
 $\End_{\mathbb R^{\scriptscriptstyle\wedge}}(\mathcal R){\downarrow}$
 is the set of all $\mathcal R$-linear functions
 in~$\mathcal R$. However, each $\mathcal R$-linear function
 $\phi$ in~$\mathcal R$ admits the representation
 $\phi(x)=cx$  for all $x\in \mathcal R$. Hence, $\End_N(G)$ consists of order bounded operators by~(1).

 $\ra$:~If $\mathbb R^{\scriptscriptstyle\wedge}\ne\mathcal R$,
 then each Hamel basis $\mathcal E$ for the vector space $\mathcal R$ over
 $\mathbb R^{\scriptscriptstyle\wedge}$
 has at least two distinct elements $e_1\ne e_2$.
 Defining some function $\phi_0:\mathcal E\to\mathcal R$ so that
 $\phi_0(e_1)/e_1\ne \phi_0(e_2)/e_2$, we may extend $\phi_0$
  to an~$\mathbb R^{\scriptscriptstyle\wedge}$-linear function
  $\phi:\mathcal R\to\mathcal R$ as in 2.1.6\,(1) which
  cannot be bounded  by~2.1.6\,(2). Therefore, the descent
  of~$\phi$ would be a band preserving linear operator that
  fails to be order bounded \big(cp.~(1)\big).~\endproof

 \section{2.2. Representation of\\ a Locally One-Dimensional Vector Lattice}

 A proper delineation of the notion of local Hamel
 basis is simply a~Hamel basis in an appropriate
 Boolean valued model. As an easy consequence we get
 WP(2)~$\iff$~WP(3).

 \Proclaim{2.2.1.}The universally complete vector lattice
 $G\assign \mathcal R{\downarrow}$ is locally one-dimensional
 if and only if\/
 ${\mathbb V}^{(\mathbb B)}\models\mathcal{R}=\mathbbm{R}^{\scriptscriptstyle\wedge}$.
 \Endproc

 \beginproof~Clearly, $[\![\mathcal R=\mathbb R^{\scriptscriptstyle\wedge}]\!]=\mathbbm{1}$
 amounts to
 $\mathcal R{\downarrow}=\mathbb R^{\scriptscriptstyle\wedge}{\downarrow}$ (cp.~\cite[3.3.3]{BVA}).
 Therefore, it suffices to check that
 $G$ is locally one-dimensional if and only if
 $G=\mathbb R^{\scriptscriptstyle\wedge}{\downarrow}$.
 However, by  \cite[3.1.1]{BVA}
 $\mathbb R^{\scriptscriptstyle\wedge}{\downarrow}$
 consists of all mixings of the shape
 $\mix_{t\in\mathbb
 R}(b_tt^{\scriptscriptstyle\wedge})$, where
 $(b_t)_{t\in\mathbb R}$ is an~arbitrary partition of unity
 in~$\mathbb B$. Considering the  properties
 of the universally complete vector lattice
 $G$  (cp.~\cite[5.2.2 and 5.2.3]{BVA}), we see  that
 $G=\mathbb R^{\scriptscriptstyle\wedge}{\downarrow}$
means the possibility of presenting each
 $x\in G$ as $\osum_{t\in\mathbb R}t\chi (b_t)\mathbbm{1}$
 with a suitable partition of unity
 $(b_t)_{t\in\mathbb R}$ in~$\mathbb B$.
 The latter rephrases as
 $G$ is locally one-dimensional,
 since we may put  $\pi_t\assign \chi (b_t)$ and rewrite
 the above presentation as
 $$
 x=\osum_{\quad t\in\mathbb R,\,t>0}
 t\pi_t\mathbbm{1} + \osum_{\quad t\in\mathbb R,\,t<0}
 t\pi_t\mathbbm{1}
 = \sup_{t\in\mathbb R,\,t>0}
 t\pi_t\mathbbm{1} - \sup_{t\in\mathbb R,\,t<0}(-t)\pi_t\mathbbm{1};
 $$
 moreover, $x^+=\sup\{t\pi_t\mathbbm{1}:\,t\in\mathbb R,\,t>0\}$ and
 $x^- =\sup\{-t\pi_t\mathbbm{1}:\,t\in\mathbb R, \,t<0\}$.~\endproof

 \subsec{2.2.2.}~Thus, the universally complete vector lattice
 $G=\mathcal R{\downarrow}$ is locally one-dimensional if
 and only if $[\![\,\mathcal R$ as a vector space over
 $\mathbb R^{\scriptscriptstyle\wedge}$ has dimension~1\,$]\!]
 =\mathbbm{1}$. Consequently, it stands to reason
 to find out what construction
 in  $G=\mathcal R{\downarrow}$ corresponds to a Hamel basis
 for the vector space $\mathcal R$ over $\mathbb R^{\scriptscriptstyle\wedge}$.
 We will presume that  $G$ is furnished with the only
  multiplicative structure making $G$~%
 into an ordered commutative algebra with ring unity
 $\mathbbm{1}\assign 1^{\scriptscriptstyle\wedge}$.

 We will say that  $x,y\in G$ differ at
 $\pi\in\goth P(G)$ provided that
 from $\rho x=\rho y$ it follows that $\pi\rho=0$
 for all $\rho\in\goth P(G)$. This amounts clearly to
 the condition $\pi (G)\subset\{|x-y|\}^{\perp\perp}$.

 A subset $\mathcal E$ of $G$ is {\it locally linearly independent\/}
 provided that, for an arbitrary nonzero band projection
 $\pi$  in~$G$, each collection of elements
 $e_1,\ldots,e_n\in\mathcal E$ that differ pairwise
 at $\pi$ and reals
 $\lambda_1,\ldots,\lambda_n\in\mathbb R$, the condition
 $\pi(\lambda_1e_1+\cdots+\lambda_ne_n)=0$ implies that
 $\lambda_k=0$ for all $k\assign 1,\ldots,n$. An inclusion-maximal
 locally linearly independent subset of~$G$ is a {\it
 local Hamel basis\/} for~$G$.

 Observe that this definition of a local Hamel basis differs
 from that given in 1.2.1. The concept of a local Hamel
 basis in 1.2.1 (cp.~\cite{MW}) corresponds to the
 interpretation of the set $\mathcal E\cup\{0\}$, where
 $[\![\,\mathcal E$ is a~Hamel basis for the vector space
 $\mathcal R$ over
 $\mathbb R^{\scriptscriptstyle\wedge}\,]\!]=\mathbbm{1}$.

 \Proclaim{}There is a local Hamel basis for an arbitrary
 universally complete vector lattice.
 \Endproc

 \beginproof~It suffices to apply the Kuratowski--Zorn Lemma
 to the inclusion ordered set of all locally linearly independent
  subsets of~$G$.~\endproof

 \Proclaim{2.2.3.}Assume that $G\assign \mathcal R{\downarrow}$,
 $\mathcal E\in\mathbb V^{(\mathbb B)}$, and
 $[\![\,\mathcal E\subset\mathcal R\,]\!]=\mathbbm{1}$. Then
 $[\![\,\mathcal E$ is a linearly independent subset
 of the vector space $\mathcal R$ over $\mathbb R^{\scriptscriptstyle\wedge}\,]\!]=\mathbbm{1}$
 if and only if
 $\mathcal E{\downarrow}$ is a~locally
 linearly independent subset of~$G$.
 \Endproc

 \beginproof~$\la$:~Put $\mathcal E^\prime\assign
 \mathcal E{\downarrow}$ and assume that
 $\mathcal E^\prime$ is locally linearly independent.
 Given a~natural $n$, let the formula
 $\varphi(n,\tau,\sigma)$ express the following: $\tau$
 and $\sigma$ are maps from
 $n\assign \{0,1\dots,n-1\}$ into $\mathbb R^{\scriptscriptstyle\wedge}$ and $\mathcal E$
 respectively,
 $\sigma(k)\ne\sigma(l)$ for different $k$ and $l$ in~$n$,
 and
 $\sum_{k\in n}\tau(k)\sigma(k)=0$.
 Denote the formula
 $$
 (\forall\,\tau)(\forall\,\sigma)\bigl(\varphi(n,\tau,\sigma)
 \ra(\forall\,k\in n)\,\tau (k)=0\bigr)
 $$
 by $\psi(n)$.
 Then the linear independence of $\mathcal E$ inside
 $\mathbb V^{(\mathbb B)}$ amounts to the equality
 $$
 \mathbbm{1}=[\![\,(\forall\,n\in\mathbb N^{\scriptscriptstyle\wedge})\,\psi(n)\,]\!]
 = \bigwedge_{n\in\mathbb N}[\![\,\psi(n^{\scriptscriptstyle\wedge})\,]\!].
 $$
 Hence, we are left with proving that
 $[\![\,\psi(n^{\scriptscriptstyle\wedge})\,]\!]=\mathbbm{1}$
 for all $n\in\mathbb N$.
 Calculate the truth  values, using
 the construction of the formula $\psi$ and the rules of
 Boolean valued analysis (cp.~\cite[2.3.8]{BVA}).
 The result is as follows:
 $$
 \bigwedge\Big\{[\![\,(\forall\,k\in
 n^{\scriptscriptstyle\wedge})\,\tau(k)=0\,]\!]:\ \tau,\sigma\in\mathbb V^{(\mathbb B)};\, [\![\,\varphi
 (n^{\scriptscriptstyle\wedge},\tau,\sigma)\,]\!]=\mathbbm{1}\Big\}.
 $$
 Take some $\tau,\sigma\in\mathbb V^{(\mathbb B)}$
 and $n\in\mathbb N$ such that
 $[\![\varphi(n^{\scriptscriptstyle\wedge},\tau,\sigma)]\!]=\mathbbm{1}$.
 Then $[\![\tau:n^{\scriptscriptstyle\wedge}
 \to\mathbb R^{\scriptscriptstyle\wedge}]\!]=\mathbbm{1}$ and
 $[\![\,\sigma:n^{\scriptscriptstyle\wedge}\to\mathcal E\,]\!]=\mathbbm{1}$.
 Moreover,
 $[\![\,\sigma(k)\ne\sigma(l)$~for distinct $k$ and $l$
 in~$n^{\scriptscriptstyle\wedge}$, and
 $\sum_{k\in n^{\scriptscriptstyle\wedge}}\tau(k)\sigma(k)=0\,]\!]=\mathbbm{1}$.

 Let $t:n\to\mathbb R^{\scriptscriptstyle\wedge}{\downarrow}$ and
 $s:n\to\mathcal E^\prime$ stand for the modified descents of $\tau$
 and $\sigma$ (cp.~\cite[3.5.5]{BVA}). Then
 $$
 \gathered
 \mathbbm{1}=[\![(\forall\,k,l\in n^{\scriptscriptstyle\wedge})
 \bigl(k\ne l\ra\sigma(k)\ne\sigma(l)\bigr)]\!]=
 % \\
 \bigwedge_{\substack{ k,l\in n\\ k\ne l}}
 [\![\sigma(k^{\scriptscriptstyle\wedge})\ne
 \sigma(l^{\scriptscriptstyle\wedge})]\!]= \bigwedge_{\substack{k,l\in n \\ k\ne l}}
 [\![s(k)\ne s(l)]\!],
 \endgathered
 $$
 and so $s(k)$ and $s(l)$ differ at the identity projection
 for $k$ and $l$ distinct. Furthermore,
 $$
 \bigg[\!\!\bigg[\sum_{k=0}^{n-1}t(k)s(k)=0\bigg]\!\!\bigg]= \bigg[\!\!\bigg[\!\sum_{\ k\in
 n^{\scriptscriptstyle\wedge}}\!\tau(k)\sigma(k)=0\bigg]\!\!\bigg]=\mathbbm{1}.
 $$
 Hence, $\sum_{k=0}^{n-1}t(k)s(k)=0$. Since
 $t(k)\in\mathbb R^{\scriptscriptstyle\wedge}{\downarrow}$
 for all $k\in n$, there is a partition of unity
 $(b_\xi)_{\xi\in\Xi}$ in $\mathbb B$ and, to each
 $k\in n$, there is a numerical family
 $(\lambda_{\xi,k})_{\xi\in\Xi}$ such that
 $$
 t(k)=\osum_{\quad\xi\in\Xi}\lambda_{\xi,k}\chi(b_\xi) \mathbbm 1\quad(k\assign 0,1,\ldots,n-1).
 $$
 Inserting these expressions into the equality
 $\sum_{k=0}^{n-1}t(k)s(k)=0$, we obtain
 $$
 0=\sum_{k=0}^{n-1}\left(\osum_{\quad\xi\in\Xi}
 \lambda_{\xi,k}\chi(b_\xi)\mathbbm 1\right)\!\!s(k)=
 \osum_{\quad\xi\in\Xi}\chi(b_\xi) \sum_{k=0}^{n-1}\lambda_{\xi,k}s(k).
 $$
 Consequently, $\chi(b_\xi)\sum_{k=0}^{n-1}\lambda_{\xi,k}s(k)=0$ and,
 since $s(k)$ and $s(l)$ differ at $\chi(b_\xi)$ for
 distinct $k,l\in n$, by the definition of
 local linear independence we have $\lambda_{\xi,k}=0$ $(k\!=0,1,\ldots,n-1)$.
 Thus $t(k)=0$ $(k\!=0,1,\ldots,n-1)$, and so
 $$
 \mathbbm{1}=\bigwedge_{k\in n}[\![t(k)=0]\!]= \bigwedge_{k\in
 n}[\![\tau(k^{\scriptscriptstyle\wedge})=0]\!]= [\![(\forall\,k\in
 n^{\scriptscriptstyle\wedge})\,\tau(k)=0]\!],
 $$
 which was required.

 $\ra$:~Assume that $[\![\,\mathcal E$ is an~$\mathbb R^{\scriptscriptstyle\wedge}$-linearly
 independent set in~$\mathcal R\,]\!]=\mathbbm1$.
 Consider arbitrary $\pi\in\goth P(G)$, $n\in\mathbb N$, $t:n\to\mathbb R$ and $s:n\to\mathcal E'$
 such that $\pi\ne0$, $s(k)$ and $s(l)$ differ at $\pi $ for distinct $k,l\in n$,
 and $\pi\sum_{k=0}^{n-1}t(k)s(k)=0$.
 Our goal is now to prove that $t(k)=0$ $(k=0,\ldots,n-1)$.

 Let $\tau,\sigma\in\mathbb V^{(\mathbb B)}$ be the modified ascents of $t$ and $s$ (cp.~\cite[3.5.5]{BVA}).
 Then, inside~$\mathbb V^{(\mathbb B)}$, we have
 $\tau:n^{\scriptscriptstyle\wedge} \to\mathbb R^{\scriptscriptstyle\wedge}$,
 $\sigma:n^{\scriptscriptstyle\wedge}\to\mathcal E$, and
 $$
 \biggl((\forall\,k,l\in n^{\scriptscriptstyle\wedge})\,\bigl(k\ne l\ra\sigma(k)\ne\sigma(l)\bigr)\wedge
 \sum_{\ k\in n^{\scriptscriptstyle\wedge}}\!
 \tau(k^{\scriptscriptstyle\wedge})\sigma(k^{\scriptscriptstyle\wedge})=0\biggr)\!
 \ra\,(\forall\,k\in n^{\scriptscriptstyle\wedge})\,\tau(k)=0.
 $$
 Calculating the truth value of the latter formula, we obtain
 $$
 b\assign \bigwedge_{\substack{ k,l\in n\\k\ne l}}
 [\![s(k)\ne s(l)]\!] \wedge\bigg[\!\!\bigg[\sum_{k=0}^{n-1}t(k)s(k)=0\bigg]\!\!\bigg]
 \leq \bigwedge_{k=0}^{n-1}[\![t(k)^{\scriptscriptstyle\wedge}=0]\!].
 $$
 According to the initial properties of $\pi$, $s$, and $t$,
 by~virtue of A3.6 we have $\pi\leq\chi(b)$ implying that $\pi t(k)^{\scriptscriptstyle\wedge}=0$
 for all $k\in n$ again by~A3.6.
 Since $\pi\neq0$, we have $t(k)=0$ $(k=0,\ldots,n-1)$.~\endproof

 \Proclaim{2.2.4.}If $\mathcal E_0$ is a locally linearly independent
 subset of $G$ and $\mathcal E\assign \mathcal E_0{\uparrow}$ then
 $[\![\,\mathcal E$~is\/~$\mathbb R^{\scriptscriptstyle\wedge}$\hbox{-}linearly independent in~$\mathcal R\,]\!]=\mathbbm{1}$.
 \Endproc

 \beginproof~By 2.2.3 it suffices to show that
 $\mathcal E_0^\prime\assign \mix(\mathcal E_0)=\mathcal E{\downarrow}
 = \mathcal E_0{\uparrow}{\downarrow}$ is locally linearly independent.
 Take some nonzero band projection $\pi$ in~$G$,
   elements $e_1,\ldots,e_n\in\mathcal E_0^\prime$ that
 differ at~$\pi$, and  reals $\lambda_1,\ldots,\lambda_n\in\mathbb R$
 satisfying
 $\pi(\lambda_1e_1+\cdots+\lambda_ne_n)=0$. There are a~partition of unity
 $(b_\xi)$ in~$\mathbb B$ and families $(g_{\xi,k})\subset\mathcal E_0$
 such that $e_k=\osum_{\xi}\chi (b_\xi)g_{\xi,k}$. Clearly,
 $\rho\assign \pi \chi (b_\eta)\ne 0$ for some index
 $\eta$. The elements $g_{\eta,1},\dots,g_{\eta,n}$ differ pairwise
 at~$\rho$ and
 $\rho(\lambda_1g_{\eta,1}+\cdots+\lambda_ng_{\eta,n})=0$. Since
  $\mathcal E_0$ is locally linearly independent,
 $\lambda_1=\cdots=\lambda_n=0$.~\endproof

 \Proclaim{2.2.5.}Assume that $G\assign \mathcal R{\downarrow}$,
 $\mathcal E\in\mathbb V^{(\mathbb B)}$, and\/
 $[\![\,\mathcal E\subset\mathcal R\,]\!]=\mathbbm{1}$.
 Then $[\![\,\mathcal E$ is a~Hamel basis for the vector space
 $\mathcal R$ over $\mathbb R^{\scriptscriptstyle\wedge}]\!]=\mathbbm{1}$
  if and only if $\mathcal E{\downarrow}$ is a~local Hamel basis
 for~$G$.
 \Endproc

 \beginproof~Immediate from 2.2.3 and 2.2.4.~\endproof

 \Proclaim{2.2.6.}A universally complete vector lattice $G$ is locally
 one-dimensional if and only if\/ $\{\mathbbm{1}\}$ is a~local Hamel
 basis for $G$.
 \Endproc

 \beginproof~Immediate from 2.2.1 and 2.2.5.~\endproof

 \section{2.3. Dedekind Cuts and Continued Fractions\\
 in a Boolean Valued Model}

 The behavior of Dedekind cuts and continued fractions in
 a Boolean valued model clarifies the equivalence
 WP(1)~$\iff$~WP(2).

 \Proclaim{2.3.1.}For all $a\subset\mathbb Q$ and $\bar{a}\subset\mathbb Q$,
 the following holds:
 $$
 (a,\bar{a})\text{ is a~Dedekind cut}\,\iff\,
 [\![\,(a^{\scriptscriptstyle\wedge},\bar{a}^{\scriptscriptstyle\wedge})
 \text{ is a~Dedekind cut}\,]\!]=\mathbbm 1.
 $$
 \Endproc

 \beginproof~Indeed, the formula $\varphi (a,\bar{a},\mathbb Q)$
 stating that $a\subset\mathbb Q$ and
 $\bar{a}\subset\mathbb Q$ comprise a~Dedekind cut, is bounded. So we
 are done by restricted transfer (cp.~A2.2).~\endproof

 \Proclaim{2.3.2.}If\/ $\mathbb B$ is $\sigma$-distributive then\,
 ${\mathbb V}^{(\mathbb B)}\models
 \mathcal R\subset{\mathbb R}^{\scriptscriptstyle\wedge}$.
\Endproc

 \beginproof~Note that the claim means precisely
  WP(1)~$\ra$~WP(2).
 Assume that $\mathbb B$ is $\sigma$\hbox{-}distributive.
 By A3.9\,(3) $\mathcal P(\omega^{\scriptscriptstyle\wedge})=
 \mathcal P(\omega)^{\scriptscriptstyle\wedge}$.
 Let $\mathcal{Q}$ denote the rationals inside~%
 ${\mathbb V}^{(\mathbb B)}$. Since the set of rationals
 can be defined by a restricted set-theoretic formula, we have
 ${\mathbb V}^{(\mathbb B)}\models\mathcal{Q}=
 \mathbb{Q}^{\scriptscriptstyle\wedge}$ (cp.~A2.2). Thus,
 we also conclude that
 $\mathcal P(\mathbb Q^{\scriptscriptstyle\wedge})=
 \mathcal P(\mathbb Q)^{\scriptscriptstyle\wedge}$.
 To demonstrate the desired inclusion
 we are to show only that
 $[\![t\in\mathcal R]\!]=\mathbbm 1$ implies
  $[\![t\in{\mathbb R}^{\scriptscriptstyle\wedge}]\!]=\mathbbm 1$.
 Assume that $[\![t\in\mathcal R]\!]=\mathbbm 1$; i.e., $t$
 is a~Dedekind cut inside~${\mathbb V}^{(\mathbb B)}$.
 We~then see
 inside ${\mathbb V}^{(\mathbb B)}$ that
 $$
 \bigl(\exists\,a\in\mathcal P(\mathbb Q^{\scriptscriptstyle\wedge})\bigr)
 \bigl(\exists\,\bar{a}\in\mathcal P(\mathbb Q^{\scriptscriptstyle\wedge})\bigr)\,
 \varphi (a,\bar{a},\mathbb Q^{\scriptscriptstyle\wedge})\wedge
 t=(a,\bar{a}),
 $$
 where $\varphi$ is the same as in~2.3.1.
 Calculating the truth value of the above formula
 and considering that $\mathcal P(\mathbb Q^{\scriptscriptstyle\wedge})=
 \mathcal P(\mathbb Q)^{\scriptscriptstyle\wedge}$, we infer
 $$
 \mathbbm 1=\bigvee_{a\subset\mathbb Q}\,
 \bigvee_{\bar{a}\subset\mathbb Q}\,
 [\![\varphi (a^{\scriptscriptstyle\wedge},\bar{a}^{\scriptscriptstyle\wedge},
 \mathbb Q^{\scriptscriptstyle\wedge})]\!]\wedge
 [\![t=(a,\bar{a})^{\scriptscriptstyle\wedge}]\!].
 $$
Choose a partition of unity $(b_\xi)\subset \mathbb B$ and two
families $(a_\xi)$ and $(\bar{a}_\xi)$ in~$\mathcal P(\mathbb Q)$
so that
 $$
 b_\xi\leq [\![\varphi (a_\xi^{\scriptscriptstyle\wedge},
 \bar{a}_\xi^{\scriptscriptstyle\wedge},
 \mathbb Q^{\scriptscriptstyle\wedge})]\!]\wedge
 [\![t=(a_\xi,\bar{a}_\xi)^{\scriptscriptstyle\wedge}]\!].
 $$
 It follows that
 $t=\mix_{\xi}b_\xi
 (a_\xi,\bar{a}_\xi)^{\scriptscriptstyle\wedge}$, and
 $b_\xi\leq [\![\varphi (a_\xi^{\scriptscriptstyle\wedge},
 \bar{a}_\xi^{\scriptscriptstyle\wedge},
 \mathbb Q^{\scriptscriptstyle\wedge})]\!]$.
 If $b_\xi\ne\mathbb 0$
 then $[\![\varphi (a_\xi^{\scriptscriptstyle\wedge},
 \bar{a}_\xi^{\scriptscriptstyle\wedge},
 \mathbb Q^{\scriptscriptstyle\wedge})]\!]=\mathbbm 1$,
 since  $\varphi(x_1,x_2,x_3)$
 is a bounded formula and the truth value
 $[\![\varphi(x_1^{\scriptscriptstyle\wedge},x_2^{\scriptscriptstyle\wedge},x_3^{\scriptscriptstyle\wedge})]\!]$
 of a bounded formula may be either $\mathbb 0$ or $\mathbb 1$
 by the definitions and rules of transformation of
 truth values (cp.~\cite[2.2.3\,(2)]{BVA}).
 By restricted transfer (cp.~A2.2 and \cite[2.2.9]{BVA})
 $\varphi (a_\xi,\bar{a}_\xi,\mathbb Q)$; i.e.,
 $(a_\xi,\bar{a}_\xi)$ is a~Dedekind cut. It is evident now that
 $b_{\xi}\leq [\![t=(a_\xi,\bar a_\xi)^{\scriptscriptstyle\wedge}
 \in\mathbb R^{\scriptscriptstyle\wedge}]\!]$. Hence,
 $[\![t\in\mathbb R^{\scriptscriptstyle\wedge}]\!]=\mathbbm 1$.~\endproof

 \subsec{2.3.3.}~We now prove the implication WP(2)~$\ra$~WP(1).
 To this end we use  continued fractions. Put
 $$
 \align
 \mathbb I&\assign \{t\in\mathbb R\ :\ 0<t<1,\ t\text{ is\ irrational}\},\\
 \mathcal I&\assign \{t\in\mathcal R:\ 0<t<1,\ t\text{ is\ irrational\} \big(inside ${\mathbb V}^{(\mathbb B)}$\big)}.
 \endalign
 $$
 It is well known that there is a bijection
 $\lambda:\mathbb I\to{\mathbb N}^{\mathbb N}$
 sending a real $t$ to the sequence $\lambda(t)=a:\mathbb N\to\mathbb N$
 of partial continued fractions of the
 continued fraction expansion of~$t$:
 $$
 t= \frac1{ a(1)+\frac1{a(2)+\frac1{a(3)+\ldots}}}\,.
 $$
 Given sequences $a:\mathbb N\to\mathbb N$ and
 $s:\mathbb N\to\mathbb I$, consider the bounded formula
 $\varphi(a,s,t,\mathbb N)$ stating that $s(1)=t^{-1}$
 and
 $$
 a(n)=\bigg[\frac{1}{s(n)}\bigg],\quad s(n+1)=\frac{1}{s(n)}-a(n),
 $$
 for all $n\in\mathbb N$, where $[\alpha]$ is the~integer part
 of~$0<\alpha\in\mathbb R$ which is expressed by the
 bounded formula $\psi (\alpha,[\alpha],\mathbb N)$:
 $$
 [\alpha]\in\mathbb N\wedge[\alpha]\leq\alpha\wedge
 (\forall\, n\in\mathbb N)(n\leq\alpha\ra n\leq[\alpha]).
 $$
 The equality $\lambda (t)=a$ means the existence
 of a sequence $s:\mathbb N\to\mathbb I$ such that
 $\varphi(a,s,t,\mathbb N)$. Call the bijection $\lambda$
 the {\it continued fraction expansion}.
 By transfer (cp.~A1.2), the continued fraction expansion
 $\tilde{\lambda}:\mathcal I\to(\aleph_0)^{\aleph_0}=(\mathbb N^{\scriptscriptstyle\wedge})^{\mathbb N^{\scriptscriptstyle\wedge}}$
 exists inside~${\mathbb V}^{(\mathbb{B})}$.

 \Proclaim{2.3.4.}Inside~${\mathbb V}^{(\mathbb{B})}$, the restriction of\/ $\tilde{\lambda}$
 to\/~$\mathbb I^{\scriptscriptstyle\wedge}$ coincides with
 $\lambda^{\scriptscriptstyle\wedge}$; i.e.,
 $$
 {\mathbb V}^{(\mathbb B)}\models(\forall\,t\in\mathbb I^{\scriptscriptstyle\wedge})\,
 \tilde{\lambda}(t)=\lambda^{\scriptscriptstyle\wedge}(t).
 $$
 \endproc

 \beginproof~The desired is true only if
 $\tilde{\lambda}(t^{\scriptscriptstyle\wedge})=
 \lambda (t)^{\scriptscriptstyle\wedge}$ for all~$t\in \mathbb I$.
 By the above definition of the bijection $\tilde{\lambda}$
 we have to demonstrate the validity inside ${\mathbb V}^{(\mathbb B)}$
 of the following formula:
 $(\exists\,s\in\mathcal I^{\mathbb N^{\scriptscriptstyle\wedge}})\,
 \varphi(\lambda(t)^{\scriptscriptstyle\wedge},s,
 t^{\scriptscriptstyle\wedge},\mathbb N^{\scriptscriptstyle\wedge})$.
By the definition of~$\lambda$ there is a sequence
 $\sigma:\mathbb N\to\mathbb I$ satisfying
  $\varphi(\lambda (t),\sigma,t,\mathbb N)$. Since
  $\varphi$ is bounded, $\mathbbm 1=[\![\varphi(\lambda(t)^{\scriptscriptstyle\wedge},
 \sigma^{\scriptscriptstyle\wedge}, t^{\scriptscriptstyle\wedge},\
 \mathbb N^{\scriptscriptstyle\wedge})]\!]$.
 Note that $\sigma^{\scriptscriptstyle\wedge}:
 \mathbb N^{\scriptscriptstyle\wedge}\to\mathbb I^{\scriptscriptstyle\wedge}
 \subset\mathcal I$; i.e., $[\![\sigma^{\scriptscriptstyle\wedge}
 \in\mathcal I^{\mathbb N^{\scriptscriptstyle\wedge}}]\!]=\mathbbm{1}$.
 Summarizing the above, we may write
 $
 [\![\,(\exists\,s\in\mathcal I^{\mathbb N^{\scriptscriptstyle\wedge}})\,
 \varphi (\lambda(t)^{\scriptscriptstyle\wedge},s,
 t^{\scriptscriptstyle\wedge},\mathbb N^{\scriptscriptstyle\wedge})\,]\!]\geq
 [\![\varphi (\lambda(t)^{\scriptscriptstyle\wedge},
 \sigma^{\scriptscriptstyle\wedge},
 t^{\scriptscriptstyle\wedge},\mathbb N^{\scriptscriptstyle\wedge})]\!]=\mathbb
 1.~\endproof
 $

 \Proclaim{2.3.5.}If\/ ${\mathbb V}^{(\mathbb B)}\models\mathcal R=
 {\mathbb R}^{\scriptscriptstyle\wedge}$ then
 $\mathbb B$ is $\sigma$-distributive.
 \Endproc

 \beginproof~By hypothesis
 $\mathcal I=\mathbb I^{\scriptscriptstyle\wedge}$
 inside~${\mathbb V}^{(\mathbb B)}$. Hence,
 $\tilde{\lambda}$ and $\lambda^{\scriptscriptstyle\wedge}$ are
 bijections, $\tilde{\lambda}$ extends
 $\lambda^{\scriptscriptstyle\wedge}$, and their
 images coincide.
 Clearly,  the domains coincide in this event too
 (and, moreover,
 $\tilde{\lambda}=\lambda^{\scriptscriptstyle\wedge}$). Therefore,
  $(\mathbb N^{\mathbb N})^{\scriptscriptstyle\wedge}=
 (\mathbb N^{\scriptscriptstyle\wedge})^{\mathbb N^{\scriptscriptstyle\wedge}}$.
 By~A3.9\,(2) we infer that  $\mathbb B$
 is~$\sigma$\hbox{-}distributive.~\endproof

 \section{PART 3. AUTOMORPHISMS AND DERIVATIONS}

 The goal of this part is to prove that if $G_{\mathbbm{C}}$
 is the complexification of a universally complete
 vector lattice $G$ then the following  are equivalent:

 \subsubsec{WP(1)}~${\mathbbm B}$ is $\sigma$-distributive;

 \subsubsec{WP(2$^\prime$)}~${\mathcal C}= {{\mathbbm
 C}}^{\scriptscriptstyle\wedge}$ inside~${{\mathbbm V}}^{({\mathbbm
 B})}$;

 \subsubsec{WP(4$^\prime$)}~Every band preserving linear operator in
 $G_{\mathbbm{C}}$ is order bounded;

 \subsubsec{WP(5$^\prime$)}~There is no nontrivial
 ${\mathbbm C}$-derivation~in the complex $f$-algebra
 $G_{\mathbbm{C}}$;

 \subsubsec{WP(6)}~Each band preserving endomorphism of
 the complex $f$-algebra $G_{\mathbbm{C}}$ is~a~band
 projection;

 \subsubsec{WP(7)}~There is no band preserving automorphism
 of\/~$G_{\mathbbm{C}}$ other than the identity.

 \section{3.1. Band Preserving Operators
 in Complex Vector Lattices}

 Consider some properties of band preserving operators in
 a~complex vector lattice.

 \subsec{3.1.1.}~A vector lattice $E$ is called
 \textit{square-mean closed\/} if for all $x,y\in E$ the
 set $\{(\cos\theta)x+(\sin\theta)y:\,0\leq\theta<2\pi\}$
 has a supremum $\mathfrak{s}(x,y)$ in $E$. Every Banach lattice
 as well as every relatively uniformly complete vector lattice is square-mean
 closed. However, a square-mean closed Archimedean vector
 lattice need not be relatively uniformly complete. If $E$
 is a square-mean closed $f$-algebra, then
 $\mathfrak{s}(x,y)^2=x^2+y^2$ for all $x,y\in E$. (It is worth
 mentioning that in \cite{ABB, BBT2} the so-called \textit{geometric-mean
 closed\/} vector lattices were also considered: this class
 is defined by the property that for all $x,y\in E_+$ the
 set $\{\frac12x+\frac1{2t}y:\,0<t<+\infty\}$ has an infimum
 $\mathfrak{g}(x,y)$ in~$E$. More details on the theme
 see in \cite[Section 3]{BBT2}, \cite{ABB, BRe}.)

 Recall that a~{\it complex vector lattice\/} is the
 complexification $E_{\mathbb C}\assign E\oplus iE$ of a~real
 square-mean closed vector lattice~$E$. Thus, each element
 $z\in E_{\mathbb C}$ in a complex vector lattice
 has the absolute value~$|z|$ defined by the formula
 $$
 |z|\assign \mathfrak{s}(x,y)\quad(z\assign x+iy\in E_{\mathbb C}).
 $$
 As usual, the notion of disjointness of elements $z\assign x+iy$ and
 $z^\prime\assign x^\prime+iy^\prime$ in~$E_{\mathbb C}$ is defined by the
 formula $z\perp z^\prime\iff|z|\wedge|z^\prime|=0$ and
 is equivalent to the relation $\{x,y\}\perp\{x^\prime,y^\prime\}$.
 An~ideal~$J$ in~$E_{\mathbb C}$ is defined as the linear subspace
 which is solid: $|x|\leq|y|$ with $x\in E_{\mathbb C}$ and
 $y\in J$ implies $x\in J$. As in the real case, a band
 in~$E_{\mathbb C}$ can be defined as
 $\{z\in E_{\mathbb C}:(\forall\,v\in V)\,z\perp v\}$,
 where $V$ is a~nonempty subset of~$E_{\mathbb C}$. The
 ideals and bands of $E_{\mathbb C}$ are precisely the
 complexifications~of  ideals and bands of~$E$ (cp.~\cite[Chapter~II, \S\,11]{Sch} and
 \cite[Section 91]{Z} for more detail).

 \subsec{3.1.2.} Consider real vector lattices~$E$ and $F$. The
 space $L(E_{\mathbb C}, F_{\mathbb C})$ of $\mathbb C$-linear operators is
 isomorphic to the complexification of the real space~$L(E, F)$ of
 $\mathbb R$-linear operators. An~operator $T\in L(E_{\mathbb C}, F_{\mathbb
 C})$ is uniquely representable as $T=T_1+iT_2$, where $T_1,T_2\in
 L(E, F)$, and an~arbitrary operator $S\in L(E, F)$ is identified
 with the canonical extension
  $\widetilde S\in L(E_{\mathbb C},
 F_{\mathbb C})$ of $S$ defined by the formula $\widetilde Sz\assign Sx+iSy$,
 $z=x+iy$. In particular, if $E$ and $F$ are considered as real
 subspaces of~$E_{\mathbb C}$ and~$F_{\mathbb C}$ then the space $L(E,F)$
 can be considered as a~real subspace of~$L(E_{\mathbb C},F_{\mathbb
 C})$.

 An~operator $T=T_1+iT_2$ is {\it positive\/} provided that $T_1\geq 0$ and
 $T_2=0$ and \textit{order bounded} provided that
 for every $e\in E_+$
 there is $f\in F_+$ satisfying $|Tx|\leq f$ whenever $|x|\leq e$.
 The space of all order bounded linear operators from
 $E_{\mathbb C}$ into~$F_{\mathbb C}$ is the complexification
 of the space of all order bounded linear operators from
 $E$ into~$F$.

 If $E_{\mathbb C}=J\oplus J^\perp$ for some ideal $J\subset
 E_{\mathbb C}$ then there is a~projection $P:E_{\mathbb C}\to E_{\mathbb
 C}$ with kernel~$J^\perp$ and range~$J$. The restriction of~$P$
 to~$E$ is a~band projection in~$E$; in particular, $P$ is
 a~positive operator. More details can be found in
 \cite[Chapter~II]{Sch} and \cite[Section 92]{Z}.

 \subsec{3.1.3.} Suppose that $F$ is a~sublattice of a~vector
 lattice $E$. As in the real case \cite[3.3.2]{DOP},
 a~linear operator~$T$ from~$F_{\mathbb C}$ to~$E_{\mathbb C}$
 is {\it band preserving\/} provided that
 $$
 z\perp z^\prime\ra Tz\perp z^\prime\quad
 (z\in F_{\mathbb C},\ z^\prime\in E_{\mathbb C}),
 $$
 where the disjointness relations are understood in~$E_{\mathbb C}$.

 \Proclaim{}A~linear operator $T\assign T_1+iT_2$ from~$F_{\mathbb C}$
 to~$E_{\mathbb C}$ is band preserving if and only if such are the
 real linear operators~$T_1,T_2:F\to E$.
 \Endproc

 \beginproof~Assume that  $T_1$ and $T_2$ are band preserving. If
 $z\assign x+iy$ and $w\assign u+iv$ are disjoint then
 $\{x,y\}\perp\{u,v\}$. Therefore,
 $\{x,y\}\perp\{T_1u-T_2v,T_1v+T_2u\}$. Hence, $z\perp Tw$
 since $Tw=(T_1u-T_2v)+i(T_1v+T_2u)$.

 Conversely, if $T$ is band preserving and  $x\in E$ and
 $u\in F$ are disjoint then $x\perp Tu=T_1u+iT_2u; $ hence,
 $x\perp\{T_1u,T_2u\}$.~\endproof

 In particular, if $E$ is a~vector lattice
 enjoying the principal projection property and $F$ is an order
 dense ideal of~$E$ then a~linear operator
 $T=T_1+iT_2:F_{\mathbbm{C}}\to E_{\mathbbm{C}}$ is band
 preserving if and only if $\pi T_k z=T_k\pi z$
 $(z\in F_{\mathbb C}$, $k=1,2)$
 for every band projection $\pi\in\goth P(E)$. An order bounded
 band preserving operator in $E_{\mathbb C}$ is called an
 \textit{orthomorphism} and the set of all orthomorphisms in
 $E_{\mathbb C}$ is denoted by $\Orth(E_{\mathbb C})$. Clearly,
 $\Orth(E_{\mathbb C})$ is the complexification of $\Orth(E)$.

 \subsec{3.1.4.} Henceforth $\mathbb B$ is a~complete Boolean algebra.
 %and ${\mathbb V}^{(\mathbb B)}$ is the corresponding Boolean valued
 %universe in which the Boolean truth value of an~arbitrary
 %set-theoretic formula~$\varphi(x_1,\ldots,x_n)$ with elements
 %$x_1,\ldots,x_n\in{\mathbb V}^{(\mathbb B)}$ is denoted
 %by~$[\![\varphi(x_1,\ldots,x_n)]\!]$. Moreover,
 %$[\![\varphi(x_1,\ldots,x_n)]\!]\in\mathbb B$
 %and the validity of~$\varphi$ in ${\mathbb V}^{(\mathbb B)}$ means by
 %definition that $[\![\varphi(x_1,\ldots,x_n)]\!]=\mathbbm{1}$
 %(cp.~A1.1).

 By the maximum principle (cp.~A1.4 and \cite[Theorem~4.3.9]{IBA}),
 there is an~element $\mathcal C\in {\mathbb V}^{(\mathbb B)}$
 for which $[\![\,\mathcal C$ is the complexes$\,]\!]
 ={\mathbbm{1}}$. Since the equality
 ${\mathbb C}={\mathbb R}\oplus i{\mathbb R}$ is expressed
 by a~bounded set-theoretic formula, from the restricted
 transfer principle A2.2 (cp.~\cite[4.2.9\,(2)]{IBA}) we
 obtain $[\![\,{\mathbb C}^{\scriptscriptstyle\wedge} =
 {\mathbb R}^{\scriptscriptstyle\wedge}\oplus
 i^{\scriptscriptstyle\wedge}{\mathbb R}^{\scriptscriptstyle\wedge}\,]\!]=
 {\mathbbm{1}}$. Moreover, ${\mathbb R}^{\scriptscriptstyle\wedge}$
 is assumed to be a~dense subfield of~$\mathcal R$;
 therefore, we can also assume that
 ${\mathbb C}^{\scriptscriptstyle\wedge}$ is a~dense subfield
 of~$\mathcal C$. If $1$ is the unity of~$\mathbb C$ then
 $1^{\scriptscriptstyle\wedge}$ is the unity of~$\mathcal
 C$ inside~$\mathbb V^{(\mathbb B)}$. We write $i$ instead
 of~$i^{\scriptscriptstyle\wedge}$ and $\mathbbm{1}$ instead of
 $1^{\scriptscriptstyle\wedge}$.

 The {\it descent\/} of~$\mathcal C$ is the set $\mathcal
 C{\downarrow}\assign \{x\in{\mathbb V}^{(\mathbb B)}:\,[\![x\in\mathcal
 C]\!]=\mathbbm{1}\}$ endowed with the structure of a~commutative complex
 ordered ring by descending the operations (cp.~A2.4 and~\cite[Section 5.3]{IBA}). Moreover, $\mathcal C{\downarrow}=\mathcal
 R{\downarrow}\oplus i\mathcal R{\downarrow}$; consequently, by
 the Gordon Theorem (cp. A3.6 and~\cite[Theorem~10.3.4]{IBA}),
 $\mathcal C{\downarrow}$ is a~universally complete complex
 vector lattice and a~complex $f$-algebra simultaneously;
 moreover, $\mathbbm{1}\assign 1^{\scriptscriptstyle\wedge}$ is the order and
 ring unity in~$\mathcal C{\downarrow}$. The space
 $\mathcal C{\downarrow}$ depends only on~$\mathbb B$ and $\mathbb C$;
 therefore, we will also use the notation $\mathbb B(\mathbb C)\assign
 \mathcal C{\downarrow}$.

 \subsec{3.1.5.} Let $\operatorname{End}_N(G_{\mathbb C})$ be the set
 of all band preserving linear operators in $G_{\mathbb C}$, where
 $G\assign \mathcal R{\downarrow}$. It is clear that
 $\operatorname{End}_N(G_{\mathbb C})$ is a~complex vector space.
 Moreover, $\operatorname{End}_N(G_{\mathbb C})$ becomes a~faithful unitary
 module over $G_{\mathbb C}$ if the operator $gT$ is defined
 by the formula $gT:x\mapsto g\cdot Tx$ $(x\in G_{\mathbb C})$. This
 follows from the fact that multiplication by an~element
 of~$G_{\mathbb C}$ is a~band preserving operator and the composition
 of band preserving operators is a~band preserving operator.

 Denote by  $\operatorname{End}_{\mathbb
 C^{\scriptscriptstyle\wedge}}(\mathcal C)$ the element of~$\mathbb
 V^{(\mathbb B)}$ that depicts the space of all $\mathbb
 C^{\scriptscriptstyle\wedge}$-linear mappings from $\mathcal C$
 into~$\mathcal C$. Then $\operatorname{End}_{\mathbb
 C^{\scriptscriptstyle\wedge}}(\mathcal C)$ is a~vector space
 over~$\mathbb C^{\scriptscriptstyle\wedge}$ inside $\mathbb
 V^{(\mathbb B)}$ and $\operatorname{End}_{\mathbb
 C^{\scriptscriptstyle\wedge}}(\mathcal C){\downarrow}$ is
 a~faithful unitary module over~$G_{\mathbb C}$.

 \subsec{3.1.6.} As in~2.1.2, we can prove that
 a~linear operator in a~universally complete vector lattice $G_{\mathbb C}$ is band preserving
 if and only if it is extensional. Since extensional mappings admit
 ascent, each operator~$T\in\operatorname{End}_N(G_{\mathbb C})$ has
 the ascent $\tau\assign T{\uparrow}$ which is the unique function from
 $\mathcal C$ into~$\mathcal C$ (inside ${\mathbb V}^{(\mathbb B)}$) satisfying the condition
 $[\![\tau (x)=Tx]\!]=\mathbbm{1}$ for all $x\in G_{\mathbb C}$
 (cp.~\cite[Theorem~5.5.6]{IBA}).

 \Proclaim{}The modules~$\operatorname{End}_N(G_{\mathbb C})$ and
 $\operatorname{End}_{\mathbb C^{\scriptscriptstyle\wedge}}(\mathcal
 C){\downarrow}$ are put into isomorphy by sending a~band preserving
 operator to its ascent.
 \Endproc

 \beginproof~Repeat the arguments of~2.1.2 with~3.1.3 and
 3.1.4 taken into account.~\endproof

 \section{3.2. Automorphisms and Derivations on the Complexes}

 We start with introducing notions and notation needed for the current and next subsections.

 \subsec{3.2.1.}~Define \textit{a~complex $f$-algebra\/}
 to be the complexification $A_{\mathbb C}$ of a~real
 square-mean closed $f$-algebra~$A$ (cp.~%
 Definition~3.1.1). The multiplication in~$A$ extends
 naturally to~$A_{\mathbb C}$ by the formula
 $$
 (x+iy)(x^\prime+iy^\prime)=
 (xx^\prime-yy^\prime)+i(xy^\prime+x^\prime y),
 $$
 and so  $A_{\mathbb C}$ becomes a commutative complex algebra. Moreover,
 $|z_1 z_2|=|z_1||z_2|$ $(z_1,z_2\in A_{\mathbb C})$.
 In this situation $A_{\mathbb C}$ is called a \textit{complex
 $f$-algebra} (cp.~\cite{BHP, Z}). A~complex $f$-algebra $A_{\mathbb C}$ is
 semiprime whenever $x\perp y$ is equivalent to $xy=0$ for all $x,y\in A_{\mathbb C}$.

 If $G$ is a~universally complete vector lattice with
 a~fixed order unity $\mathbbm{1}\in G$ then there is
 a~unique multiplication in~$G$ which makes~$G$ into
 an~$f$-algebra and $\mathbbm{1}$ into the multiplicative
 unity. Thus, $G_{\mathbb C}$ is an~example of a~complex
 $f$-algebra. We~will always keep this circumstance in
 mind while considering a~universally complete vector
 lattice an~$f$-algebra.

 \subsec{3.2.2.} Given an~algebra~$A$ and a~subalgebra~$A_0$
 of~$A$, we call a~linear operator $D:A_0\to A$   a~{\it
 derivation\/} provided that
 $$
 D(uv)=D(u)v+uD(v)\quad (u,v\in A_0).
 $$
 The kernel of a~derivation is a~subalgebra. A~nonzero derivation
 is called {\it nontrivial}.

 An~{\it endomorphism of an~algebra\/} is a~linear multiplicative
 operator in it. A~bijective endomorphism is  an~{\it
 automorphism}. The identical automorphism is commonly referred to as the
 {\it trivial automorphism}.

 If the above definitions of an~automorphism and a~derivation
 relate to an~algebra over a~field $\mathbb P$ then we also speak of
 $\mathbb P$-automorphisms and $\mathbb P$-derivations.

 For completeness of exposition, we give some properties of
 the complexes which we need below. In the next section we
 will give the Boolean valued interpretation of these
 properties. As above, $\mathcal C$ is the complexes inside
 $\mathbbm{V}^{(\mathbbm{B})}$. Recall that $\mathcal{C}$ includes the
 subfield $\mathbbm{C}^{\scriptscriptstyle\wedge}$
 inside~$\mathbbm{V}^{(\mathbbm{B})}$.  The following  was  obtained in
 \cite{Kus2}:

 \subsec{3.2.3.} \theorem{}Inside ${\mathbb V}^{(\mathbb B)}$,
 the field $\mathbbm{C}^{\scriptscriptstyle\wedge}$ is algebraically
 closed in~$\mathcal{C}$. In particular, if\/
 $\mathbbm{V}^{(\mathbbm{B})}\models{\mathbbm{C}}^{\scriptscriptstyle\wedge}\ne\mathcal{C}$
 then
 $\mathbbm{V}^{(\mathbbm{B})}\models\text{``$\,\mathcal{C}$
 is a~transcendental extension of\/
 $\mathbbm{C}^{\scriptscriptstyle\wedge}$.''}$
 \Endproc

 Thus, under the canonical embedding of the complexes
 into the Boolean valued model, either ${\mathbb
 C}^{\scriptscriptstyle\wedge}=\mathcal C$ or the field of
 complexes is a~transcendental extension of some
 subfield of~$\mathcal C$. The same is true for the reals. To
 analyze this situation, we need the notion of an~algebraic
 or transcendence basis of a~field over some subfield.

 Let $\mathbb P$ be a~subfield of~$\mathbb C$ such that $\mathbb
 C$ is a~transcendental extension of~$\mathbb P$. By
 the Steinitz Theorem~\cite[Chapter~5, \S\,5, Theorem~1]{Bou1},
 there is
 a~transcendence basis~$\mathcal E\subset\mathbb C$. This means that
 the set $\mathcal E$ is algebraically independent over~$\mathbb P$
 and $\mathbb C$~is an~algebraic extension of the field~$\mathbb
 P(\mathcal E)$ obtained by addition of the elements of~$\mathcal
 E$ to~$\mathbb P$. The field $\mathbb P(\mathcal E)$ is  a~{\it
 pure extension\/} of~$\mathbb P$.

 \Proclaim{3.2.4.}Let $\mathbb C$ be a~transcendental
 extension of a~field\/~$\/\mathbb P$.~Then there is a~nontrivial
 $\mathbb P$-automorphism~of~$\mathbb C$.
 \Endproc

 \beginproof~Let $\mathcal E$ be a~transcendence basis for the
 extension~$\mathbb C$ over~$\mathbb P$. Since $\mathbb C$ is
 an~algebraically closed extension of~$\mathbb P(\mathcal E)$,
 every $\mathbb P$-automorphism~$\phi$ of the field $\mathbb P(\mathcal E)$
 extends to a~$\mathbb P$-automorphism~$\Phi$ of the field~$\mathbb C$
 (cp.~\cite[Chapter~5, \S\,4, the Corollary to Theorem~1]{Bou1}).
 It is clear that if $\phi$ is nontrivial then so is $\Phi$.

 To construct a~nontrivial $\mathbb P$-automorphism in $\mathbb P(\mathcal
 E)$, we firstly consider the case when $\mathcal E$ contains only one
 element~$e$; i.e., when $\mathbb C$ is an~algebraic extension of
 a~simple transcendental extension~$\mathbb P(e)$. Take
 $a,b,c,d\in\mathbb P$ such that $ad-bc\ne 0$. Then
 $e^\prime=(ae+b)/(ce+d)$ is a~generator of the field~$\mathbb P(e)$
 different from~$e$. The field $\mathbb P(e)=\mathbb P(e^\prime)$ is
 isomorphic to the field of rational fractions in one variable~$t$;
 consequently, the linear-fractional substitution
 $t\mapsto(at+b)/(ct+d)$ defines a~$\mathbb P$-automorphism~$\phi$ of
 the field $\mathbb P(e)$ which takes $e$ into~$e^\prime$
 (cp.~\cite[Section 39]{Waer}).

 Assume now that $\mathcal E$ contains at least two different
 elements~$e_1$ and~$e_2$ and take an~arbitrary bijective mapping
 $\phi_0:\mathcal E\to\mathcal E$ for which $\phi_0(e_1)= e_2$. Again,
 using the circumstance that $\mathbb C$ is an~algebraically closed
 extension of~$\mathbb P(\mathcal E)$, we can construct a~$\mathbb
 P$-automorphism~$\phi$ of~$\mathbb C$ such that
 $\phi_0(e)=\phi(e)$ for all $e\in\mathcal E$ (cp.~\cite[Chapter~5,
 Section 6, Proposition~1]{Bou1}). Clearly, $\phi$ is nontrivial.~\endproof

 \Proclaim{3.2.5.}Let $\mathbb C$ be a~transcendental
 extension of a~field\/~$\mathbb P$. Then there is a~nontrivial $\mathbb
 P$-derivation on~$\mathbb C$.
 \Endproc

 {\beginproof} We again use a~transcendence basis~$\mathcal E$ for the
 extension $\mathbb C$ over~$\mathbb P$. It is well known that every
 derivation of~$\mathbb P$ extends onto a~purely
 transcendental extension; moreover, this extension is defined
 uniquely by prescribing arbitrary values at the elements of a~transcendence basis
 (cp.~\cite[Chapter~V, Section 9, Proposition~4]{Bou1}). Thus, for every
 mapping $d:\mathcal E\to\mathbb C$, there is a~unique derivation
 $D:\mathbb P(\mathcal E)\to\mathbb C$ such that $D(e)=d(e)$ for
 all $e\in\mathcal E$ and
 $D(x)=0$ for $x\in\mathbb P$. Now, $\mathbb C$ is a~separable
 algebraic extension of~$\mathbb P(\mathcal E)$; consequently,
 $D$ admits a~unique extension to some derivation $\overline D:\mathbb
 C\to\mathbb C$ (cp.~\cite[Chapter~V, Section 9, Proposition~5]{Bou1}). It is
 obvious that the freedom in the choice of~$d$ guarantees that
 $\overline D$ is nontrivial.~\endproof

 \subsec{3.2.6.} Using the same arguments as above, we can
 show that some analogs of~3.2.3 and~3.2.5 are valid for the
 reals. More precisely, the following is valid:

 \subsubsec{(1)}~$[\![\,{\mathbb R}^{\scriptscriptstyle\wedge}$
 is algebraically closed in~$\mathcal R\,]\!]=\mathbbm{1}$;

 \subsubsec{(2)}~If
 $\mathbbm{V}^{(\mathbbm{B})}\models{\mathbbm{R}}^{\scriptscriptstyle\wedge}\ne\mathcal{R}$,
 then
 $\mathbbm{V}^{(\mathbbm{B})}\models\text{``$\,\mathcal{R}$
 is a~transcendental extension of
 $\mathbbm{R}^{\scriptscriptstyle\wedge}$;''}$

 \subsubsec{(3)}~If $\mathbb R$ is a~transcendental extension
 of a~field~$\mathbb P$ then there is a~nontrivial
 $\mathbb P$-derivation on~$\mathbb R$.

 However, 3.2.4 is not valid for the reals: there
 is no nontrivial automorphism on~$\mathbb R$. This is connected
 with the fact that $\mathbb R$ is not an~algebraically closed
 field.

 \subsec{3.2.7.}
 \proclaim{Theorem.} Let ${\mathbb C}$ be an extension of
 an algebraically closed subfield\/ ${\mathbb P}$. Then the
 following  are equivalent:

 \subsubsec{(1)}~${\mathbb P}={\mathbb C}$;

 \subsubsec{(2)}~Every ${\mathbb P}$-linear function in
 ${\mathbb C}$ is order bounded;

 \subsubsec{(3)}~There are no nontrivial\/ ${\mathbb P}$-derivations
 on ${\mathbb C}$;

 \subsubsec{(4)}~Each ${\mathbb P}$-linear endomorphism of\/
 ${\mathbb C}$ is the zero or identity function;

 \subsubsec{(5)}~There is no ${\mathbb P}$-linear automorphism
 of\/ ${\mathbb C}$ other than the identity.
 \Endproc

 \beginproof~The equivalence (1)~$\iff$~(2) is
 checked by using a~Hamel basis of the vector space
 ${\mathbb C}$ over~${\mathbb P}$. The remaining equivalences
 follow on replacing a~Hamel basis with a~transcendence basis
 from 3.2.4 and 3.2.5 (for details, cp.~\cite{Kus2}).~\endproof

 \section{3.3. Automorphisms and Derivations on Complex $f$-Algebras}

 Consider the question of existence of nontrivial automorphisms and
 derivations on a~universally complete complex $f$-algebra. In
 this section $G$ is a~universally complete vector lattice with a~fixed
 multiplicative structure, $E$ is a~subring and a~sublattice
 in~$G$, while $G_{\mathbb C}\assign G\oplus iG$ and $E_{\mathbb C}\assign E\oplus iE$.

\subsec{3.3.1.}  {\sl Let $D\in L(E_{\mathbb C},G_{\mathbb
 C})$ and $D=D_1+iD_2$. The~operator~$D$ is a~complex derivation if
 and only if $D_1$ and $D_2$ are real derivations from~$E$
 into~$G$.}

 \beginproof~We only have to insert $D\assign D_1+iD_2$ in the equality
 $D(uv)=D(u)v+uD(v)$, take  $u\assign x\in E$ and $v\assign y\in E$, and
 then equate the real and imaginary parts of the resulting
 relation.~\endproof

 \proclaim{3.3.2.}If $E^{\perp\perp}=G$ then each
 derivation from $E_{\mathbb C}$ into~$G_{\mathbb C}$ is a~band
 preserving operator.
 \Endproc

 {\beginproof}~By~3.1.3 and~3.3.1, we only have to establish that every
 real derivation is a~band preserving operator. Let $D:E\to G$ be
 a~real derivation. Take disjoint $x,y\in E$. Since the relation
 $x\perp y$ in an~$f$-algebra implies $xy=0$, we have
 $0=D(xy)=D(x)y+xD(y)$. But the elements $D(x)y$ and $xD(y)$ are
 disjoint as well by the definition of an~$f$-algebra; therefore,
 $D(x)y=0$ and $xD(y)=0$. Hence, since the
 $f$-algebra~$E$ is faithful, we obtain $D(x)\perp y$ and $x\perp D(y)$. Now,
 consider disjoint $x\in E$ and $g\in G$. By condition, the order ideal
 $I$ generated by  $\{x\}^\perp\cup\{x\}$ is order dense
 in~$G$; therefore, without loss of generality we may assume that
 $g\in I$. At the same time, $|g|\leq y$ for some $y\in E_+$;
 consequently, $D(x)\perp g$ by the above.~\endproof

 \subsec{3.3.3.} Let $\mathcal D(\mathcal C{\downarrow})$ be the
 set of all derivations on the $f$-algebra $\mathcal C{\downarrow}$ and let $\mathcal M_N(\mathcal C{\downarrow})$ be
 the set of all band preserving automorphisms of $\mathcal C{\downarrow}$. Let $\mathcal D_{\mathbb
 C^{\scriptscriptstyle\wedge}}(\mathcal C)$ and $\mathcal M_{\mathbb
 C^{\scriptscriptstyle\wedge}}(\mathcal C)$ be the elements of~${\mathbb
 V}^{(\mathbb B)}$ that depict the sets of all $\mathbb
 C^{\scriptscriptstyle\wedge}$-derivations and all $\mathbb
 C^{\scriptscriptstyle\wedge}$\hbox{-}automorphisms in~$\mathcal C$. Clearly,
 $\mathcal D(\mathcal C{\downarrow})$ is a~module over~$\mathcal
 C{\downarrow}$ and $[\![\,\mathcal D_{\mathbb
 C^{\scriptscriptstyle\wedge}}(\mathcal C)$ is a~complex vector
 space$]\!]=\mathbbm{1}$.

 \Proclaim{}The descent and ascent produce isomorphisms
 between the modules $\mathcal D_{\mathbb
 C^{\scriptscriptstyle\wedge}}(\mathcal C){\downarrow}$ and $\mathcal
 D(\mathcal C{\downarrow})$ as well as bijections between $\mathcal
 M_{\mathbb C^{\scriptscriptstyle\wedge}}(\mathcal C){\downarrow}$ and
 $\mathcal M_N(\mathcal C{\downarrow})$.
 \Endproc

  {\beginproof}~The proof follows from~3.1.6. We only have to note that
   $T\in\operatorname{End}_N(\mathcal C{\downarrow})$ is
 a~derivation (automorphism) if and only if $[\![\,\tau\assign T{\uparrow}$
 is a~derivation (auto\-morphism)$\,]\!]=\mathbbm{1}$.~\endproof

 \Proclaim{3.3.4.}An~order bounded derivation and
 an~order bounded band preserving automorphism of a~universally
 complete $f$-ring $G_{\mathbb C}$ are trivial.
 \Endproc

 {\beginproof} We may assume that $G_{\mathbb C}=\mathcal C{\downarrow}$.
 If $T$ is a~derivation (a~band preserving automorphism) of the
 $f$-ring $G_{\mathbb C}$ then $[\![\,\tau\assign T{\uparrow}$ is a~$\mathbb
 C^{\scriptscriptstyle\wedge}$-derivation ($\mathbb
 C^{\scriptscriptstyle\wedge}$-automorphism)
 of~$\mathcal
 C\,]\!]=\mathbbm{1}$. Moreover, $T$ is order bounded if and only if
 $[\![\,\tau$ is order bounded in~$\mathcal C\,]\!]=\mathbbm{1}$. However,
 every order bounded $\mathbb
 C^{\scriptscriptstyle\wedge}$-derivation on the field $\mathcal C$ is
 zero and every order bounded $\mathbb
 C^{\scriptscriptstyle\wedge}$-automorphism is the identity
 mapping. In the first case we have $T=0$ and in the second,
 $T=I$.~\endproof

 \Proclaim{3.3.5.}If\/ $\mathbb V^{(\mathbb B)}\models \mathbb
 C^{\scriptscriptstyle\wedge}\ne\mathcal C$ then there exist
 a~nontrivial derivation and a~nontrivial band preserving
 automorphism on the universally complete complex $f$-algebra $\mathbb
 B(\mathbb C)=\mathcal C{\downarrow}$.
 \Endproc

 {\beginproof} It follows from the condition $\mathbb
 C^{\scriptscriptstyle\wedge}\ne\mathcal C$ that $\mathcal C$ is
 a~transcendental extension of~$\mathbb
 C^{\scriptscriptstyle\wedge}$ inside~${\mathbb V}^{(\mathbb B)}$
 (cp.~3.2.3). By~3.2.4 and~3.2.5, there exist a~nontrivial
 $\mathbb C^{\scriptscriptstyle\wedge}$-derivation
 $\delta:\mathcal C\to\mathcal C$ and a~nontrivial
 $\mathbb C^{\scriptscriptstyle\wedge}$-automorphism
 $\alpha:\mathcal C\to\mathcal C$. If $D\assign \delta{\downarrow}$ and
 $A\assign \alpha{\downarrow}$ then, according to~3.3.3, $D$ is
 a~nontrivial derivation and~$A$ is a~nontrivial band preserving
 automorphism of the $f$-algebra~$\mathcal C{\downarrow}$.~\endproof

 \subsec{3.3.6.} \theorem{}For an~arbitrary complete Boolean
 algebra~$\mathbb B$ the following are equivalent:

 \subsubsec{(1)}~$\mathbb B$ is $\sigma$-distributive;

 \subsubsec{(2)}~${\mathbb V}^{(\mathbb B)}\models \mathcal C={\mathbb
 C}^{\scriptscriptstyle\wedge}$;

 \subsubsec{(3)}~All band preserving linear operators on the
 universally complete vector lattice $\mathbb B(\mathbb C)=\mathcal C{\downarrow}$ are order
 bounded;

 \subsubsec{(4)}~There are no nonzero derivations on the complex $f$-algebra
 $\mathbb B(\mathbb C)=\mathcal C{\downarrow}$;

 \subsubsec{(5)}~Each band preserving endomorphism of the complex
 $f$-algebra $\mathbb B(\mathbb C)=\mathcal C{\downarrow}$ is a~band
 projection;

 \subsubsec{(6)}~In the complex $f$-algebra $\mathbb B(\mathbb C)=\mathcal
 C{\downarrow}$ there are no nontrivial band preserving
 automorphisms.
 \Endproc

 \beginproof~(1) $\iff$ (2):~As is known (cp.~Section~2.3),
 a~Boolean algebra~$\mathbb B$ is $\sigma$-distributive if and only if
 $\mathbb V^{(\mathbb B)}\models\mathbb
 R^{\scriptscriptstyle\wedge}=\mathcal R$. Hence, using the
 restricted transfer principle~A2.2
 (\cite[4.2.9\,(2)]{IBA}), we conclude that
 $
 \mathbb V^{(\mathbb B)}\models\mathcal C=\mathcal R\oplus i\mathcal
 R=\mathbb R^{\scriptscriptstyle\wedge} \oplus i\mathbb
 R^{\scriptscriptstyle\wedge}= \mathbb C^{\scriptscriptstyle\wedge}.
 $
 The converse is proved similarly.

 (2) $\ra$ (3):~If $\mathbb V^{(\mathbb B)}\models\mathbb
 C^{\scriptscriptstyle\wedge}=\mathcal C$ then, inside $\mathbb
 V^{(\mathbb B)}$, the set $\operatorname{End}_{\mathbb
 C^{\scriptscriptstyle\wedge}}(\mathcal C)$ consists of the
 functions $\tau:\mathcal C\to\mathcal C$ of the form $\tau(z)=cz$,
 where $c\in \mathcal C$. But then the operator
 $T\assign \tau{\downarrow}$ from~$\mathcal C{\downarrow}$ into~$\mathcal
 C{\downarrow}$ has the form $T(u)=gu$ for some $g\in\mathcal
 C{\downarrow}$.

 (3) $\ra$ (2):~It follows from~(3) that
 all band preserving linear operators
 are order bounded in the
 universally complete vector lattice
 $\mathcal R{\downarrow}$. Thus, $\mathbb V^{(\mathbb B)}\models\mathbb
 R^{\scriptscriptstyle\wedge}=\mathcal R$ \big(cp.~2.1.7\,(2)\big);
 and so $\mathbb V^{(\mathbb B)}\models\mathcal
 C=\mathbb C^{\scriptscriptstyle\wedge}$.

 (3)~$\ra$~(4):~This follows from~3.3.2 and~3.3.4.

 (3)~$\ra$~(5): A~band preserving endomorphism
 $T:\mathcal C{\downarrow}\to\mathcal C{\downarrow}$ admits
 the representation $T=T_1+iT_2$, where $T_1$ and~$T_2$ are
 band preserving linear operators in the universally complete
 vector lattice $\mathcal R{\downarrow}$ (cp.~3.1.3).
 By~(3), $T_1$ and~$T_2$ are order bounded; consequently,
 $T_lx=c_lx$ ($x\in\mathcal R{\downarrow}$) for some constants
 $c_1,c_2\in\mathcal R{\downarrow}$. Hence, $Tz=c\cdot z$
 $(z\in\mathcal C{\downarrow})$, where $c\assign c_1+ic_2$.
 Multiplicativity of~$T$ implies $c^2=c$; therefore, the equalities
 $c^2_1-c_2^2=c_1$ and $2c_1c_2=c_2$ are valid.
 If $\pi \assign [c_2]$ is
 the~projection in~$\mathcal R{\downarrow}$ onto the
 band~$\{c_2\}^{\perp\perp}$ then from the second equality we
 derive $\pi c_1=(1/2)\pi (\mathbbm{1})$, while the first equality
 implies $-\pi(c_2^2)=(1/4)\pi(\mathbbm{1})$. The last is possible only
 for $\pi=0$; hence, $c_2=0$ and $0\leq c^2_1=c_1$. But we also
 have $0\leq(\mathbbm{1}-c_1)^2=\mathbbm{1}-c_1$; consequently,
 $c_1\leq\mathbbm{1}$. Now, we see that the operator $x\mapsto
 T_1x=c_1x$ is a~band projection in~$\mathcal R{\downarrow}$ and,
 in view of~$T_2=0$, its canonical extension to~$\mathcal
 C{\downarrow}$ coincides with~$T$.

 (5) $\ra$ (6):~This is obvious.

 The implications~(4) $\ra$ (2) and~(6) $\ra$ (2)
 follow from~3.3.5.

 (4) $\ra$ (2):~If the equality $\mathcal C={\mathbb
 C}^{\scriptscriptstyle\wedge}$ is violated inside~${\mathbb
 V}^{(\mathbb B)}$ then $b\assign [\![\mathcal C={\mathbb C}^{\scriptscriptstyle\wedge}]\!]<\mathbbm{1}$. But then
 $b^\ast=[\![\mathcal C\ne{\mathbb C}^{\scriptscriptstyle\wedge}]\!]
 \ne{\mathbb 0}$. The inequality $\mathcal C\ne{\mathbb
 C}^{\scriptscriptstyle\wedge}$ is valid in the Boolean valued
 model ${\mathbb V}^{(\mathbb B_0)}$ over the Boolean algebra~$\mathbb
 B_0\assign [\mathbbm 0,b^\ast]$. By~3.3.5, there is a~nonzero
 derivation~$D$ on the band~$b^\ast\mathcal C{\downarrow}$. The
 unique extension $D\oplus 0$ of the operator~$D$ coinciding with
 zero on the band $b\mathcal C{\downarrow}$ is a~nonzero derivation
 on~$\mathcal C{\downarrow}$,~too.

 (6)~$\ra$~(2): Similarly, using 3.3.5,
 for the same $b\in \mathbb B$ we can find a~nontrivial
 automorphism $A^*$ of the band $b^\ast\mathcal C{\downarrow}$.
 If $A$ is the identity mapping in the band
 $b\mathcal C{\downarrow}$ then $A^*\oplus A$ is a~nontrivial
 automorphism of~$\mathcal C{\downarrow}$.~\endproof

 \subsec{3.3.7.} \corollary{}For a~universally complete real
 vector lattice~$G$ with a~fixed structure of an~$f$-algebra, the
 following are equivalent:

 \subsubsec{(1)} $\mathbb B\assign \goth P(G)$ is
 a~$\sigma$-distributive Boolean algebra;

 \subsubsec{(2)} There are no nontrivial derivations on the complex
 $f$-algebra $G_{\mathbb C}$;

 \subsubsec{(3)} There are no nontrivial band preserving automorphisms of
 the complex $f$-algebra $G_{\mathbb C}$. \Endproc

 \section{PART 4. VARIATIONS ON THE THEME}

 In this part we consider briefly the band preserving
 phenomenon in some natural environments (the endomorphisms
 of lattice ordered modules, bilinear operators on vector
 lattices, and  derivations in $AW^\ast$-algebras) and
 state some problems that may be viewed as versions of the
 Wickstead problem.

 \section{4.1. The Wickstead Problem in Lattice Ordered Modules}

 In this section we state a kind of the Wickstead problem for
 lattice ordered modules.

 \subsec{4.1.1.}~Let $K$ be a~lattice ordered ring, and let
 $X$  be a lattice ordered module over~$K$. The Wickstead
 problem  for lattice ordered modules can be stated as
 follows:

 \smallskip
 \textrm{WP(A)}:~\textsl{When are all
 band preserving $K$-linear endomorphisms of a lattice
 ordered $K$-module $X$ order bounded?}
 \smallskip

 Little is known about this problem. Boolean valued analysis
 provides a transfer principle which might send \textrm{WP} to
 \textrm{WP(A)}. Below we describe the class of lattice
 ordered modules for which this transfer works perfectly.

 \subsec{4.1.2.}~A subset $S$ of $K$ is {\it dense\/}
 provided that $S^\perp=\{0\}$; i.e., the equality $k\cdot S=\{0\}$
 implies $k=0$ for all $k\in K$. A ring $K$ is {\it rationally
 complete\/} whenever, to each dense ideal $J\subset K$ and each
 group homomorphism $h:J\to K$ such that $h(kx)=kh(x)$
 for all $k\in K$ and $x\in J$, there is an element $r$ in
 $K$ satisfying $h(x)=rx$ for all $x\in J$. A ring $K$ is
 rationally complete if and only if $K$ is
 selfinjective (cp.~\cite[Theorem 8.2.7\,(3)]{IBA}).

 \subsec{4.1.3.}~If $\mathcal{K}$ is an ordered field inside
 $\mathbbm{V}^{(\mathbbm{B})}$ then $\mathcal{K}{\downarrow}$
 is a rationally complete semiprime $f$-ring, and there is
 an isomorphism $\chi$ of $\mathbbm{B}$ onto the Boolean
 algebra $\mathfrak{B}(\mathcal{K}{\downarrow})$ of the annihilator
 ideals (coinciding in the case under consideration with the
 Boolean algebra of all bands) of $\mathcal{K}{\downarrow}$
 such that
 $$
 b\le[\![x=0]\!]\,\iff\, x\in \chi(b^*)\quad (x\in
 K,\ b\in B)
 $$
 (cp.~\cite[Theorem 8.3.1]{IBA}). Conversely, assume that
 $K$ is a rationally complete semiprime $f$-ring and
 $\mathbbm{B}$ stands for the Boolean algebra $\mathfrak{B}(K)$
 of all annihilator ideals (bands) of $K$. Then there is
 an element $\mathcal{K}\in\mathbbm{V}^{(\mathbbm{B})}$, called
 the  {\it Boolean valued representation of\/} $K$, such
 that $[\![\,\mathcal{K}$ is an ordered field\,$]\!]=\mathbbm{1}$
 and the lattice ordered rings $K$ and $\mathcal{K}{\downarrow}$
 are isomorphic (cp. \cite[Theorem 8.3.2]{IBA}).

 \subsec{4.1.4.}~A $K$-module $X$ is  {\it separated\/}
 provided that for every dense ideal $J\subset K$
  the identity  $Jx=\{0\}$
 implies $x=0$. Recall that a $K$-module $X$ is
 {\it injective\/} whenever, given a $K$-module $Y$,
 a $K$-submodule
 $Y_0\subset Y$, and a $K$-homomorphism $h_0:Y_0\to X$,
 there exists a $K$-homomorphism $h:Y\to X$ extending
 $h_0$. The Baer criterion says that a $K$-module $X$ is
 injective if and only if for each ideal $J\subset K$ and
 each $K$\hbox{-}homomorphism $h:J\to X$ there exists $x\in X$ with
 $h(a)=xa$ for all $a\in J$~(cp.~\cite{Lam}).

 \subsec{4.1.5.}~Let $\mathcal{X}$ be a~vector lattice over an ordered
 field $\mathcal{K}$ inside $\mathbbm{V}^{(\mathbbm{B})}$,
 and let $\chi:\mathbbm{B}\to\mathfrak{B}(\mathcal{K}{\downarrow})$ be
 a~Boolean isomorphism from 4.1.3. Then
 $\mathcal{X}\!\!\downarrow$ is a~separated unital
 injective lattice ordered module over
 $\mathcal{R}{\downarrow}$ satisfying
 $$
 b\le[\![x=0]\!]\,\iff\,\chi(b)x=\{0\}
 \quad (x\in\mathcal{X}\!\!\downarrow,\ b\in\mathbbm{B}).
 $$
 Conversely, let $K$ be
 a~rationally complete semiprime $f$-ring,
 $\mathbbm{B}\assign \mathfrak{B}(K)$, and let $\mathcal{K}$ be
 the Boolean valued representation of~$K$.
 Assume that $X$ is a~unital separated injective
 lattice ordered $K$-module. Then there exists some
 $\mathcal{X}\in\mathbbm{V}^{(\mathbbm{B})}$ such that
 $[\![\,\mathcal{X}$~is a~vector lattice over the ordered
 field $\mathcal{K}\,]\!]=\mathbbm{1}$ and there are
 algebraic and order isomorphisms
 $\jmath:K\to\mathcal{K}\!\!\downarrow$ and
 $\imath:X\to\mathcal{X}\!\!\downarrow$ such that
 $$
 \imath(ax)=\jmath(a)\imath(x)\quad (a\in K,\ x\in X)
 $$
 (cp.~\cite[Theorems 8.3.12 and 8.3.13]{IBA}). Thus, the
 Boolean transfer principle is applicable to unital separated
 injective lattice ordered modules over rationally complete
 semiprime $f$-rings. Consider an example.

 \subsec{4.1.6.}~Let $\mathbf{B}$ be a complete Boolean
 algebra and let $\mathbbm{B}$ be a complete subalgebra
 of $\mathbf B$. We say that $\mathbf{B}$ is
 $\mathbbm{B}$-$\sigma$-{\it distributive\/} if for every
 sequence $(b_{n})_{n\in\mathbbm{N}}$ in $\mathbf{B}$ we
 have
 $$
 \bigvee _{\varepsilon\in\mathbbm{B}^{\mathbbm{N}}}
 \bigwedge_{n\in\mathbbm{N}}\varepsilon (n)b_{n}=\mathbbm{1},
 $$
 where $\varepsilon (n)b_{n}\assign \bigl(\varepsilon(n)\wedge b_{n}\bigr)\vee\bigl(\varepsilon(n)^\ast\wedge b_n^\ast\bigr)$
 and~$b^\ast$ is the~complement of~$b\in\mathbf{B}$. Clearly, the
 $\{\mathbbm{0},\mathbbm{1}\}$-$\sigma$-distributivity of
 $\mathbf{B}$ means that $\mathbf{B}$ is $\sigma$-distributive
 \big(cp.~1.3.1\,(3)\big).

 There exists a $\mathcal{B}\in\mathbbm{V}^{(\mathbbm{B})}$
 such that $[\![\,\mathcal{B}$ is a~complete Boolean
 algebra$\,]\!]=\mathbbm{1}$ and $\mathcal{B}\!\!\downarrow$ is
 a complete Boolean algebra isomorphic to $\mathbf{B}$
 (cp.~\cite[Theorem 4.7.11]{IBA}). Moreover, $\mathbf{B}$
 is $\mathbbm{B}$-$\sigma$-distributive if and only if
 $\mathcal{B}$ is $\sigma$-distributive inside
 $\mathbbm{V}^{(\mathbbm{B})}$.  We~now interpret Theorem 1.3.7
 inside $\mathbbm{V}^{(\mathbbm{B})}$ to obtain:

 \subsec{4.1.7.}~\theorem{}Let $X$ be a universally complete
 vector lattice with a fixed order unity   $\mathbbm{1}$
 and let $K$ be an order closed sublattice containing $\mathbbm{1}$.
 Put $\mathbf{B}\assign \mathfrak{E}(X)\assign \mathfrak{E}(\mathbbm{1}_X)$
 and $\mathbb{B}\assign \mathfrak{E}(K)\assign \mathfrak{E}(\mathbbm{1}_K)$.
 Then $K$ is a rationally complete $f$-algebra, $X$ is an
 injective lattice ordered $K$-module, and the following
 are equivalent:

 \subsubsec{(1)}~$\mathbf{B}$ is $\mathbbm{B}$-$\sigma$-distributive;

 \subsubsec{(2)}~Every~element $x\in X_{+}$ is locally
 $K$-constant, i.e., $x=\sup_{\xi \in \Xi }a_{\xi}\pi _{\xi }\mathbbm{1}$
 for some family $(a_{\xi })_{\xi\in\Xi}$ of elements of $K$
 and a~family $(\pi _{\xi })_{\xi \in \Xi }$ of pairwise
 disjoint band projections in $X$;

 \subsubsec{(3)}~Every band preserving $K$-linear endomorphism of $X$ is order bounded.
 \Endproc

 \section{4.2. The Wickstead Problem for Bilinear Operators}

 In this section we present the main results of~\cite{K18}.

 \subsec{4.2.1.}~Let $E$ be a vector lattice.
 A bilinear operator $b:E\times E\ra E$ is
 {\it separately band preserving\/} provided that the mappings
 $b(\cdot,e):x\mapsto b(x,e)$ and $b(e,\cdot):x\mapsto b(e,x)$
 $(x\in E)$ are band preserving for all
 $e\in E$  or, which is the same, provided that $b(L\times E)\subset L$
 and $b(E\times L)\subset L$ for every band $L$ in $E$.

 \Proclaim{4.2.2.} Assume that $E$ is a vector lattice and
 $b:E\times E\ra E$ is a bilinear operator. Then
 the following  are equivalent:

 \subsubsec{(1)}~$b$ is separately band preserving;

 \subsubsec{(2)}~$b(x,y)\in\{x\}^{\perp\perp}\cap
 \{y\}^{\perp\perp}$ for all $x,y\in E$;

 \subsubsec{(3)}~$b(x,y)\perp z$ for all $z\in E$
 provided that $x\perp z$ or $y\perp z$.\\
 If $E$ has the principal projection property, then (1)--(3)
 are equivalent to:

 \subsubsec{(4)}~$\pi b(x,y)=b(\pi x,\pi y)$ for
 every $\pi\in\mathfrak{P}(E)$ and all $x,y\in E$;

 \subsubsec{(5)}~$\pi b(x,y)=b(\pi x,y)=b(x,\pi y)$ for
 every $\pi\in\mathfrak{P}(E)$ and all $x,y\in E$.
 \Endproc

 \beginproof~We omit the routine arguments which are similar to
 \cite[Theorem 8.2]{AB}.~\endproof

 \subsec{4.2.3.}~Let $E$ and $F$ be vector lattices. A bilinear operator $b:E\times E\to F$
 is  {\it orthosymmetric\/} provided that $|x|\wedge |y|=0$
 implies $b(x,y)=0$ for arbitrary $x,y\in E$
 (cp.~\cite{BuR2}). The difference of two positive orthosymmetric
 bilinear operators is {\it orthoregular} (cp.~\cite{BuK, K15}).
 Recall also that a bilinear operator $b$ is
 {\it symmetric\/} or {\it antisymmetric\/} provided that
 $b(x,y)=b(y,x)$ or $b(x,y)=-b(y,x)$ for all $x,y\in E$.

 The following important property of orthosymmetric
 bilinear operators was established
 in~\cite[Corollary 2]{BuR2}:

 \Theorem{}If $E$ and $F$ are vector lattices then every orthosymmetric
 positive bilinear operator from $E\times E$ into $F$ is
 symmetric.
 \Endproc

 \subsec{4.2.4.}~It is evident from 4.2.2 that a separately
 band preserving bilinear operator is ortho\-sym\-metric.
 Hence, all orthoregular separately band preserving operators
 are symmetric by 4.2.3. At the same time an order bounded
 separately band preserving bilinear operator $b$ is of the
 form $b=\pi\odot$ with $\pi$ an orthomorphism on
 $E^{{\scriptscriptstyle\odot}}$ and $\odot$ is the canonical
 bimorphism from $E\times E$ to $E^{{\scriptscriptstyle\odot}}$
 (cp.~\cite[Section 2]{BuK} and \cite{BuR4}). This brings up the following question:

 \smallskip
 \textrm{WP(B)}:~{\sl Under what conditions are all
 separately band preserving bilinear operators in a~vector
 lattice  symmetric? Order bounded?}
 \smallskip

 In the case of a~universally complete vector lattice the
 answer is similar to the linear case and is presented below
 in 4.2.5. The general case was not yet examined.

 \subsec{4.2.5.} \theorem{}Let $G$ be a universally complete
 vector lattice and
 let $\mathbbm{B}\assign \mathfrak{P}(G)$ denote the complete Boolean
 algebra of all bands in $G$. Then the following
 are equivalent:

 \subsubsec{(1)}~$\mathbbm{B}$ is $\sigma$-distributive;

 \subsubsec{(2)}~There is no nonzero separately band
 preserving antisymmetric bilinear operator in $G$;

 \subsubsec{(3)}~All separately band preserving bilinear
 operators in $G$ are symmetric;

 \subsubsec{(4)}~All separately band preserving bilinear
 operators in $G$ are order bounded.
 \Endproc

\beginproof~The only nontrivial implication is
 (2)~$\ra$~(1).

 We may assume that $G=\mathcal{R}{\downarrow}$. Suppose that $\mathbbm{B}$ is not
 $\sigma$-distributive. Then
 $\mathbbm{R}^{\scriptscriptstyle\wedge}\ne\mathcal{R}$
 by WP(1)~$\iff$~WP(2) (cp. Section~2.3) and
 a~separately band preserving antisymmetric
 bilinear operator can be constructed on using the bilinear
 version of 2.1.6\,(1).
 Indeed, inside $\mathbbm{V}^{(\mathbbm{B})}$, a Hamel bases $\mathcal{E}$ for $\mathcal{R}$ over
 $\mathbbm{R}^{\scriptscriptstyle\wedge}$ contains at least
 two different elements $e_1\ne e_2$. Define a function
 $\beta_0:\mathcal{E}\times\mathcal{E}\to\mathcal{R}$ so
 that $1=\beta_0(e_1,e_2)=-\beta_0(e_2,e_1)$, and
 $\beta(e_1^\prime,e_2^\prime)=0$ for all other pairs
 $(e_1^\prime,e_2^\prime)\in\mathcal{E}\times\mathcal{E}$
 \big(in~particular, $0=\beta_0(e_1,e_1)=\beta_0(e_2,e_2)$\big).
 Then $\beta_0$ can be extended to an
 $\mathbbm{R}^{\scriptscriptstyle\wedge}$-bilinear function
 $\beta:\mathcal{R}\times\mathcal{R}\ra\mathcal{R}$.
 The descent $b$ of $\beta$ is a separately band preserving
 bilinear operator in $G$ by 4.2.6, the bilinear version
 of 2.1.5. Moreover, $b$ is
 nonzero and antisymmetric, since $\beta$ is nonzero and
 antisymmetric by construction. This contradiction proves
 that $\mathbbm{R}^{\scriptscriptstyle\wedge}=\mathcal{R}$
 and $\mathbbm{B}$~is $\sigma$-distributive.~\endproof

 \subsec{4.2.6.}~Let $B\!L_N(G)$ stand for the set
 of all separately band preserving
 bilinear operators in $G=\mathcal{R}{\downarrow}$. Clearly, $B\!L_N(G)$ becomes
 a~faithful unitary module over  $G$ provided that we define
 $gT$ as $gT:x\mapsto g\cdot Tx$ for all $x\in G$.
 Denote by
 $B\!L_{\mathbbm{R}^{\scriptscriptstyle\wedge}}(\mathcal{R})$
 the element of~$\mathbbm{V}^{(\mathbbm{B})}$ that depicts the
 space of all $\mathbbm{R}^{\scriptscriptstyle\wedge}$-bilinear
 mappings from $\mathcal{R}\times\mathcal{R}$ into~%
 $\mathcal{R}$. Then
 $B\!L_{\mathbbm{R}^{\scriptscriptstyle\wedge}}(\mathcal{R})$
 is a~vector space over
 $\mathbbm{R}^{\scriptscriptstyle\wedge}$
 inside $\mathbbm{V}^{(\mathbbm{B})}$, and
 $B\!L_{\mathbbm{R}^{\scriptscriptstyle\wedge}}(\mathcal{R}){\downarrow}$
 is a~faithful unitary module over~$G$.

 \Proclaim{}The modules $B\!L_N(G)$ and $B\!L_{\mathbbm{R}^{\scriptscriptstyle\wedge}}
 (\mathcal{R}){\downarrow}$ are isomorphic by sending each
 band preserving bilinear operator to its ascent.
 \Endproc

 \beginproof~See 2.1.5.~\endproof

 \smallskip
 \proclaim{4.2.7.}There exists a nonatomic universally
 complete vector lattice in which all separately band
 preserving bilinear operators are symmetric and order bounded.
 \Endproc

 \beginproof~It follows from 4.2.5 and 1.3.8.~\endproof

\goodbreak
 \section{4.3. The Noncommutative Wickstead Problem}

 The relevant information on the theory of Baer
 $\ast$-algebras and $AW^\ast$-algebras  can be found
 in~\cite{Berb, Chil, DOP}.

 \subsec{4.3.1.}~A \textit{Baer $\ast$-algebra\/} is a complex
 involutive algebra $A$ provided that, for each nonempty
 $M\subset A$, there is a projection, i.e., a~hermitian
 idempotent, $p$ satisfying $M^\perp=pA$, where
 $M^\perp\assign \{y\in A:(\forall\,x\in M)\,xy=0\}$ is the
 right annihilator of~$M$. Clearly,
 this amounts to the condition that each left annihilator
 has the form ${}^\perp M=Aq$ for an appropriate projection
 $q$. To each left annihilator $L$ in a~Baer $\ast$-algebra
 there is a unique projection $q_L\in A$ such that $x=xq_L$
 for all $x\in L$ and $q_L y=0$ whenever
 $y\in L^\perp$. The mapping $L\mapsto q_L$ is an isomorphism
 between the poset of left annihilators and the poset of
 all projections. Thus, the poset $\mathfrak{P}(A)$ of all projections
 in a Baer $\ast$-algebra is an~order complete lattice.
 \big(Clearly, the formula $q\leq p\iff q=qp=pq$,
 sometimes pronounced  as  ``$p$ contains $q$,'' specifies some order on the set of
 projections $\mathfrak{P}(A)$.\big)

  An element $z$ in $A$ is  \textit{central\/}
 provided that $z$ commutes with every member of~$A$; i.e.,
 $(\forall\,x\in A)\,xz=zx$. The \textit{center\/}
 of a Baer $\ast$-algebra  $A$ is the set $\mathcal{Z}(A)$
 comprising central elements. Clearly, $\mathcal{Z}(A)$
 is a commutative Baer $\ast$-subalgebra of~$A$, with
 $\lambda\mathbbm{1}\in\mathcal{Z}(A)$ for all
 $\lambda\in\mathbbm{C}$. A \textit{central projection\/}
 of $A$ is a projection belonging to $\mathcal{Z}(A)$. Put
 $\mathfrak{P}_c(A)\assign \mathfrak{P}(A)\cap\mathcal{Z}(A)$.

 \subsec{4.3.2.}~A {\it derivation\/} on a Baer $\ast$-algebra $A$
 is a linear operator $d:A\to A$ satisfying $d(xy)=d(x)y+xd(y)$
for all $x,y\in A$. A derivation $d$ is  {\it inner\/}
 provided that $d(x)=ax-xa$ $(x\in A)$ for some $a\in A$. Clearly, an
 inner derivation vanishes on $\mathcal{Z}(A)$ and is
 $\mathcal{Z}(A)$-linear, i.e., $d(ex)=ed(x)$ for all $x\in A$
 and $e\in\mathcal{Z}(A)$.

 Consider a derivation $d:A\to A$ on a Baer
 $\ast$-algebra $A$. If $p\in A$ is a central projection
 then $d(p)=d(p^2)=2pd(p)$. Multiplying this identity
 by $p$ we have $pd(p)=2pd(p)$ so that $d(p)=pd(p)=0$.
 Consequently, every derivation vanishes on the linear span
 of $\mathfrak{P}_c(A)$, the set of all central projections.
 In particular, $d(ex)=ed(x)$
 whenever $x\in A$ and $e$ is a linear combination of central
 projections. Even if the linear span of central projections
 is dense in a sense  in $\mathcal{Z}(A)$, the derivation $d$
 may fail to be $\mathcal{Z}(A)$-linear.

 This brings up the natural question: \textsl{Under what
 conditions is every derivation
 $Z$-linear on a Baer $\ast$-algebra $A$  provided that $Z$ is a Baer $\ast$-subalgebra
 of $\mathcal{Z}(A)$?}

 \subsec{4.3.3.}~An~{\it $AW^\ast$-algebra\/} is a
 $C^\ast$-algebra with unity $\mathbbm{1}$ which is also a~Baer
 $\ast$\hbox{-}algebra. More explicitly, an~$AW^\ast$-algebra is
 a~$C^\ast$-algebra whose every right annihilator has the
 form~$pA$, with $p$ a~projection. Clearly, $\mathcal{Z}(A)$
 is a commutative $AW^*$\hbox{-}subalgebra of~$A$. If $\mathcal{Z}(A)=
 \{\lambda\mathbbm{1}:\lambda\in\mathbbm{C}\}$ then the
 $AW^*$-algebra $A$ is  an~$AW^*$\hbox{-}{\it factor}.

 \Proclaim{} A~$C^\ast$-algebra $A$ is an~$AW^\ast$-algebra
 if and only if the following hold:

 \subsubsec{(1)} Each orthogonal family in $\goth P(A)$
 has a supremum;

 \subsubsec{(2)} Each maximal commutative $\ast$-subalgebra
 $A_0\subset A$ is a Dedekind complete $f$-algebra (or,
 equivalently, coincides with the least norm closed
 $\ast$-subalgebra containing all  projections of~$A_0$).
 \Endproc

 \subsec{4.3.4.}~Given an $AW^\ast$-algebra $A$,
 define the two sets $C(A)$ and $S(A)$ of measurable and locally
 measurable operators, respectively. Both are Baer
 $\ast$-algebras, cp. \cite{Chil}. Suppose that $\Lambda$ is an
 $AW^\ast$-subalgebra in $\mathcal{Z}(A)$, and $\Phi$ is a
 $\Lambda$ valued trace on $A_+$. Then we may define another
 Baer $\ast$-algebra, $L(A,\Phi)$, of $\Phi$-measurable
 operators. The center $\mathcal{Z}(A)$ is a vector lattice
 with a strong unity, while the centers of
 $C(A)$, $S(A)$, and $L(A,\Phi)$ coincide with the universal
 completion of~$\mathcal{Z}(A)$.
 If $d$ is a derivation  on $C(A)$, $S(A)$, or
 $L(A,\Phi)$ then $d(px)=pd(x)$ $\bigl(p\in\mathfrak{P}_c(A)\bigr)$ so
 that $d$ can be considered as band preserving in a sense
 (cp. 1.1.1\,(4) and 3.1.3).

 \smallskip
 \textrm{WP(C)}: \textsl{When are all derivations on $C(A)$, $S(A)$,
 or $L(A,\Phi)$  inner?}

 \subsec{4.3.5.}~The classification of $AW^\ast$-algebras
 into types is determined from the structure of their lattices
 of projections $\mathfrak{P}(A)$ \cite{DOP, Sak}. We recall  only the definition
 of type~I $AW^\ast$-algebra. A~projection $\pi\in A$ is
  {\it abelian\/} if  $\pi A\pi$ is
 a~commutative algebra. An algebra $A$ has {\it type\/}~I provided that
 each  nonzero projection in $A$ contains a~nonzero abelian
 projection.

 A~$C^\ast$-algebra $A$ is  ${\mathbbm B}$-{\it embeddable}
 provided that there is a~type~I $AW^\ast$-algebra $N$ and a~$\ast$-monomorphism
 $\imath: A\ra N$ such that ${\mathbb B}={\mathfrak P}_c(N)$
 and $\imath (A)=\imath(A)^{\prime\prime}$, where
 $\imath (A)^{\prime\prime}$ is the bicommutant
 of~$\imath (A)$ in~$N$. Note that in this event $A$ is
 an~$AW^\ast$-algebra and ${\mathbbm B}$ is a~complete
 subalgebra of~${\mathfrak P}_c(A)$.

 \subsec{4.3.6.} \theorem{}Let $A$ be a type I
 $AW^\ast$-algebra, let $\Lambda$ be an $AW^\ast$\hbox{-}subalgebra of
 $\mathcal{Z}(A)$, and let $\Phi$ be a $\Lambda$ valued faithful
 normal semifinite trace on $A$. If the complete Boolean algebra
 $\mathbbm{B}\assign \mathfrak{P}(\Lambda)$ is $\sigma$-distributive
 and $A$ is $\mathbbm{B}$-embeddable, then every derivation
 on $L(A,\Phi)$ is inner.
 \Endproc

 \beginproof~We briefly sketch the proof.
 Let $\mathcal{A}\in\mathbbm{V}^{(\mathbbm{B})}$
 be the Boolean valued representation of~$A$. Then $\mathcal{A}$
 is a von Neumann algebra inside $\mathbbm{V}^{(\mathbbm{B})}$.
 Since  the Boolean valued interpretation preserves classification
 into types, $\mathcal{A}$ is of type I. Let $\varphi$ stand for
 the Boolean valued representation of~$\Phi$. Then $\varphi$
 is a $\mathcal{C}$ valued faithful normal semifinite
 trace on $\mathcal{A}$ and the descent of $L(\mathcal{A},\varphi)$
 is $\ast$-$\Lambda$-isomorphic to $L(A,\Phi)$, cp.~\cite{KCh}.
 Suppose that $d$ is a derivation on $L(A,\Phi)$ and
 $\delta$ is the Boolean valued representation~of~$d$. Then $\delta$
 is a $\mathcal{C}$ valued $\mathbbm{C}^{\scriptscriptstyle\wedge}$-linear
 derivation on $L(\mathcal{A},\varphi)$. Since $\mathbbm{B}$ is
 $\sigma$-distributive, $\mathcal{C}=
 {\mathbbm C}^{\scriptscriptstyle\wedge}$ inside
 $\mathbbm{V}^{(\mathbbm{B})}$ and $\delta$ is
 $\mathcal{C}$-linear. But it is well known that any derivation
 on a type I von Neumann algebra is inner, cp. \cite{AAKud}.
 Therefore, $d$ is also inner.~\endproof

 \normalsize{

 \section{PART 5. COMMENTS}
 \centerline{\scshape 5.1. Comments on Part 1}
\medskip

 \subsec{5.1.1.}~The theory of orthomorphisms stems from
 Nakano~\cite{Nak2}. Ortho\-mor\-phisms have been studied by
 many authors under various names (cp. \cite{AB}):
 {\it dilatators\/}
 (Nakano~\cite{Nak2}), {\it essentially positive
 operators\/} (Birkhoff~\cite{Bir}), {\it polar
 preserving endomorphisms\/} (Conrad and
 Diem~\cite{CD}), {\it multiplication operators\/}
 (Buck~\cite{Buc} and~Wickstead~\cite{Wic1}),
 and {\it stabilisateurs\/} (Meyer~\cite{Mey1}). The
 main stages of this development as well as the various
 aspects of the theory of orthomorphisms are reflected in
 the books by Bigard, Keimel, and Wolfenstein
 \cite{BKW}, Aliprantis and Burkinshaw~\cite{AB},
Zaanen \cite[Chapter 20]{Z} etc.;  also see the
 survey papers by Bukhvalov \cite[Section 2.2]{Buh1}
 and Gutman \cite[Chapter 6]{Gut5}.

 \subsec{5.1.2.}~%Theorem 3.3.2\,(1) stems from Nakano~\cite{Nak2, Nak3}. ???????????????????????????????????
 Order continuity of an extended orthomorphism (cp.~1.1.4)
 was established independently by Bigard and Keimel
 \cite{BK} and Conrad and Diem \cite{CD}
 using functional representation. A direct proof was found
 by Luxemburg and Schep \cite{LS}.
 Commutativity of every Archimedean $f$\hbox{-}algebra was proved by
Birkhoff and Pierce \cite{BirP};  this paper
 also introduced  the concept of~$f$\hbox{-}algebra.
The lattice ordered algebras were surveyed by Boulabiar, Buskes,
and Triki \cite{BBT1, BBT2}.
 The fact that $\Orth(D,E)$ is a vector lattice under
the pointwise algebraic and latticial  operations was also obtained in~%
 \cite{BK} and~\cite{CD}. Extensive is the
 bibliography on the theory of orthomorphisms; and so we
 indicate a~portion of it: \cite{AVK, AW, Ber3, BK, DuM,
 Gut6, Gut1, HP3, HW, Lux2, LS, MW, Pag1, Pag2, Wic1, Wic3, Z2}.

 \subsec{5.1.3.}~The terms ``local linear independence'' and ``local
 Hamel basis'' (coined in \cite{MW})
 appeared in~\cite{AVK}  as
 $d$-\textit{independence} and $d$-\textit{basis}. Using
 these concepts Abramovich and Kitover \cite{AK1}
 gave complete description for a band preserving
 projection $P$ on a Dedekind complete vector lattice $E$.
 The order bounded part $\pi P$ of $P$ (cp.~1.1.2) is a
 band projection, whereas the unbounded part
 $P_0\assign P|_{E_0}$, with $E_{0}\assign \pi^\perp(E)$, is uniquely
 determined from the following conditions: (1)~every
 principal band in $E_0$ is laterally complete;
(2)~$P_0^{-1}(0)$ is componentwise closed; i.e.,
 $\mathfrak{E}(u)\subset P_0^{-1}(0)$ for all
 $0\leq u\in P_0^{-1}(0)$; (3)~$L\cap P_0^{-1}(0)$ is laterally
 complete for each principal band $L$ in $E_0$. Cp.~\cite{AK3}
 for  applications of this concept.

 \subsec{5.1.4.}~The notions of $d$-independence and
 $d$-basis can be introduced in an arbitrary vector lattice
 (cp.~\cite{AK2}).
 A collection $(x_\gamma)_{\gamma\in\Gamma}$ of elements
 in a vector lattice $E$ is  $d$-\textit{independent}
 provided that for each band $B$ in $E$, each finite subset
 $\{\gamma_1,\dots,\gamma_n\}$ of $\Gamma$, and  each family of nonzero
 scalars $c_1,\dots,c_n$ the condition
 $\sum_{i=1}^nc_i x_{\gamma_i}\perp B$ implies that
 $x_{\gamma_i}\perp B$ for $i=1,\dots,n$. A~$d$-independent
 system $(x_\gamma)_{\gamma\in\Gamma}$ is
 a~$d$-{\it basis\/} provided that for each $x\in E$ there is a full
  system $(B_\alpha)_{\alpha\in\textrm{A}}$ of pairwise
 disjoint bands in~$E$ and a system of elements
 $(y_\alpha)_{\alpha\in\textrm{A}}$ in $E$ such that each
 $y_\alpha$ is a linear combination of elements in
 $(x_\gamma)_{\gamma\in\Gamma}$ and $(x-y_\alpha)\perp B_\alpha$
 for all $\alpha\in\textrm{A}$.

 \subsec{5.1.5.}~Theorem 1.3.7 can be considered as an~exhaustive
 answer to the Wickstead problem about the order boundedness of
 all band preserving operators. However, a~new
 notion of locally one-dimensional vector lattice crept
 into the~answer. The novelty of this notion led to
 the~conjecture that it coincides with that of
 discrete ($=$~atomic) vector lattice. In~1981, Abramovich, Veksler,
 and Koldunov~\cite[Theorem~2.1]{AVK1} gave a~proof for
 existence of an~order unbounded band preserving operator in
 every nondiscrete universally complete vector lattice,
 thus seemingly corroborating the~conjecture that a~locally
 one-dimensional  vector lattice is discrete
 (also cp.~\cite[Section~5]{Ab1}). However, the~proof
 was  erroneous. Later in~1985, McPolin and
 Wickstead~\cite[Section~3]{MW} gave an~example of
 a~nondiscrete locally one-dimensional vector lattice,
 confuting the~conjecture. However, there was an~error
 in the~example. Finally,  Wickstead~\cite{AW} stated the~conjecture
 as an~open  question in~1993.

 \subsec{5.1.6.}~In the case of a~universally complete
 vector lattice, a~band preserving order  unbounded operator can
 be constructed on using~$\mathbbm{V}^{(\mathbbm{B})}$. Moreover,
 inside an~appropriate $\mathbbm{V}^{(\mathbbm{B})}$ this
 problem reduces to existence of a~discontinuous solution
 $\varphi:\mathcal{R}\to\mathcal{R}$ to
 the Cauchy functional equation $\varphi(s+t)=\varphi(s)+\varphi(t)$
 $(s,t\in\mathcal{R})$ with an additional property
 $\varphi(\lambda s)=\lambda\varphi(s)$
 $(\lambda\in\mathbbm{R}^{\scriptscriptstyle\wedge}$,
 $s\in\nobreak\mathcal{R})$.  Let~$E$ be a~universally complete vector lattice
 such that ${\mathbb V}^{(\mathbb B)}\models\mathbbm{R}^{\scriptscriptstyle\wedge}\neq\mathcal{R}$
 (cp.~Section 2.3) with $\mathbbm{B}\assign \goth B(E)$.
 Then $\mathcal{R}$ is an~infinite-dimensional
 vector space over~$\mathbbm{R}^{\scriptscriptstyle\wedge}$ inside~%
 $\mathbbm{V}^{(\mathbbm{B})}$.  By
 the Kuratowski--Zorn Lemma, there exists an~$\mathbbm{R}^{\scriptscriptstyle\wedge}$-linear
 but not $\mathcal{R}$-linear function $\varphi:\mathcal{R}\to\mathcal{R}$
 inside~$\mathbbm{V}^{(\mathbbm{B})}$.  The operator~%
 $\Phi_0\assign  \varphi{\downarrow}:\mathcal{R}{\downarrow}\to\mathcal{R}{\downarrow}$
 is linear, band preserving, but order  unbounded. If
 $\iota$ is an~isomorphism of~$E$ onto~$\mathcal{R}{\downarrow}$
 then $\Phi\assign \iota^{-1} \circ \Phi_0\circ\iota$ is
 an~order unbounded band preserving operator in~$E$.

 \section{5.2. Comments on Part 2}

 We see that the claim of Theorem \textrm{WP} reduces to
 simple properties of reals and cardinals. However, even the
 reader who mastered the technique (of ascending and descending)
 of Boolean valued analysis might find the above demonstration
 bulky as compared with the standard proof in the
 articles by Abramovich, Veksler, and Koldunov~\cite{AVK1},
 McPolin and Wickstead \cite{MW}, and~Gutman~\cite{Gut6}.
 However, the aim of the exposition in~Part 2
 was not to simplify the  available proof but rather
 demonstrate that the Boolean approach to the problem
 reveals many new interconnections. A few clarifications are
 now in order.

 \subsec{5.2.1.}~Since the space of ${\mathbb R}^{\scriptscriptstyle\wedge}$-linear
 functions in~$\mathcal R$ admits a complete description
 that uses a~Hamel basis \big(cp.~2.1.7\,(2)\big); therefore,
 $\End_N(\mathcal R{\downarrow})$ may be described completely
by means of a~(strict) local Hamel basis. However, this approach
will evoke some problems of unicity.

 \subsec{5.2.2.}~The dimension $\delta(\mathcal R)$ of the vector
 space $\mathcal R$ over~${\mathbb R}^{\scriptscriptstyle\wedge}$
 is a~cardinal inside~${\mathbb V}^{(\mathbb B)}$.
The object $\delta(\mathcal R)$ carries   important information on
the interconnection of the Boolean algebra $\mathbb B$ and the
reals~$\mathbb R$. By the properties of Boolean valued ordinals,
we obtain the representation
 $\delta(\mathcal R)=\mix_\xi b_\xi\alpha_\xi^{\scriptscriptstyle\wedge}$,
 where $(b_\xi)$ is a~partition of unity in~$\mathbb B$ and
  $(\alpha_\xi)$ is a~family of standard cardinals.
 This representation is an instance of a ``decomposition series''
 of~$\mathbb B$ such that the principal ideals
 $[\mathbb 0,b_\xi]$  are ``$\alpha_\xi$-homogeneous'' in a~sense.

 \subsec{5.2.3.}~If we replace the class of band preserving
 linear operators with the class of band preserving additive operators
 then the equivalence WP(1)$\,\iff\,$WP(4) fails to
 hold in Theorem WP. Moreover, in each nonzero universally complete
 vector lattice there exist order unbounded band preserving additive
 operators. This reflects the fact that there is no
 Boolean valued model satisfying
 ${\mathbb V}^{(\mathbb B)}\!\models\!\mathcal R={\mathbb Q}^{\scriptscriptstyle\wedge}$.

 \subsec{5.2.4.}~The property of $\lambda$ in~2.3.4
 is usually referred to as  {\it absolute definability}.
 Gordon \cite{Gord}
 called a continuous function absolutely definable if it possesses
 an analogous property. For instance, the functions
  $e^x$, $\log x$, $\sin x$, and $\cos x$ are absolutely definable. In particular, these
  functions reside  in every Boolean valued universe, presenting the mappings
  from   $\mathcal R$ to $\mathcal R$ that are continuations of
  the corresponding functions
  $\exp^{\scriptscriptstyle\wedge}(\cdot)$,
 $\log^{\scriptscriptstyle\wedge}(\cdot)$,
 $\sin^{\scriptscriptstyle\wedge}(\cdot)$, and
 $\cos^{\scriptscriptstyle\wedge}(\cdot)$
 from $\mathbb R^{\scriptscriptstyle\wedge}$
 into~$\mathbb R^{\scriptscriptstyle\wedge}$. Practically all functions
 admitting a constructive definition are absolutely definable.

 \subsec{5.2.5.}~Consider a band preserving operator
 $S:\mathcal R{\downarrow}\to\mathcal R{\downarrow}$
satisfying the Cauchy exponential equation:
 $S(x+y)=S(x)S(y)$ for all $x,y\in\mathcal R{\downarrow}$.
 If, moreover, $S$ enjoys the condition $S(\lambda x)=S(x)^\lambda$
 for all $0<\lambda\in\mathbb R$ and
 $x\in\mathcal R{\downarrow}$; then we call
 $S$ an {\it exponential operator}.
 Say that $S$ is order bounded if
 $S$ takes order bounded sets into order bounded sets.
 If $\sigma$ is the ascent
 of~$S$ then $\sigma$ is exponential
 inside~$\mathbb V^{(\mathbb B)}$. Therefore, in the class of
 functions bounded above on some nondegenerate interval
 we see that $\sigma =0$ or $\sigma(x)=e^{cx}$ for all
 $x\in\mathcal R$ and some $c\in\mathcal R$.
 This implies that WP(1)--WP(7)
 of Theorem WP amount to the following:

 \smallskip
 \subsubsec{WP(8)}~{\sl Each band preserving exponential
 operator $S$ on~$\mathbb B(\mathbb R)\assign \mathcal R{\downarrow}$
 is order bounded (and thus, $S$ may be presented as
 $S(x)=e^{cx}$ for all  $x\in\mathcal R{\downarrow}$
 and some $c\in\mathcal R{\downarrow}$ or $S$
 is identically zero).}

 \subsec{5.2.6.}~An analogous situation takes place
 if $S$ satisfies the {\it Cauchy logarithmic  equation\/}
 $S(xy)=S(x)+S(y)$ for all $0\ll x,y\in\mathcal R{\downarrow}$ and
 enjoys the condition $S(x^\lambda)=\lambda S(x)$ for all
  $\lambda\in\mathbb R$ and
 $x\gg0$. (The record $0\ll x$ means that
 $0\leq x$ and $x^{\perp\perp}=\mathcal R{\downarrow}$.) We call
 an $S$ of this sort a~{\it logarithmic operator}. We may now
 formulate another equivalent claim as follows:

 \smallskip
 \subsubsec{WP(9)}~{\sl Every band preserving logarithmic operator
 $S$ on~$\{x\in\mathbb B(\mathbb R)\assign \mathcal R{\downarrow}:x\gg0\}$ is order bounded
 (and, consequently, $S$ may presented as $S(x)=c\log x$
 for all $0\ll x\in\mathcal R{\downarrow}$ with some $c\in\mathcal
 R{\downarrow}$).}

 \subsec{5.2.7.}~Instead of using continued fraction expansions
 in Section 2.3 we may involve binary expansions. In this
 event we have to construct a bijection
 of~$\mathcal P(\omega)$ onto some set of reals
 and apply A3.9\,(3) in place of~A3.9\,(2).

\section{5.3. Comments on Part 3}

 Part~3 may be considered as an evidence of the productivity
 of combining  algebraic and logical methods  in
 operator theory.

 \subsec{5.3.1.} Using the same arguments as in~3.3.5 and~3.3.6,
 from~3.2.6, we can infer that if
 $\mathbb R^{\scriptscriptstyle\wedge}\ne\mathcal R$ then
 there are nontrivial derivations on the real $f$-algebra
 $\mathcal R{\downarrow}$. Thus, in the class of universally
 complete real vector lattices with a~fixed structure of
 an~$f$-algebra we have WP(1)$\,\iff\,$WP(5); i.e.,
 the absence of nontrivial derivations is equivalent to the
 $\sigma$-distributivity of the~base of the algebra under consideration.
 At the same time there are no nontrivial band preserving
 automorphisms of the $f$-algebra $\mathcal R{\downarrow}$,
 regardless of the properties of its base.

 \subsec{5.3.2.} It is well known that if $Q$ is a~compact
 space then there are no nontrivial derivations on the algebra
 $C(Q,{\mathbb C})$ of complex valued continuous functions
 on~$Q$; for example, see~\cite[Chapter 19, Theorem~21]{AD}.
 At the same time, we see from~3.3.6\,(1),\,(4) that if
 $Q$ is an~extremally disconnected  compact space and the
 Boolean algebra of the clopen sets of~$Q$ is not $\sigma$-distributive
 then there is a~nontrivial derivation on~$C_\infty(Q,{\mathbb C})$.

 \subsec{5.3.3.} Let $(\Omega,\Sigma,\mu)$ be a~measure
 space with the direct sum property (cp.~\cite[1.1.7 and 1.1.8]{DOP}).
 {\sl The Boolean algebra $\mathbb{B}\assign \mathbb{B}(\Omega,\Sigma,\mu)$
 of measurable sets modulo negligible sets is
 $\sigma$-distributive if and only if $\mathbb{B}$ is atomic
 \big(and thus isomorphic to the boolean $\mathcal{P}(A)$
 of a nonempty set $A$\big)}. Indeed, suppose that~$\mathbb{B}$ is not atomic.
 By choosing a nonzero atomless coset $b_0\in\mathbb{B}$ of finite measure,
 taking an instance $B_0\in b_0$, and replacing $(\Omega,\Sigma,\mu)$
 with $(B_0,\Sigma_0,\mu\vert{}_{\Sigma_0})$, where $\Sigma_0=\{B\cap B_0:B\in\Sigma\}$,
 we may assume that~$\mu$ is finite and~$\mathbb{B}$ is atomless.
 Define a~strictly positive
 countably additive function $\nu:\nobreak\mathbb{B}\to\mathbb{R}$
 by $\nu(b)=\mu(B)$ where $b\in\mathbb{B}$ is the coset of
 $B\in\Sigma$.
 Since any finite atomless measure admits
 halving, by induction it is easy to construct a sequence
 of finite partitions $P_m\assign \{b_{1}^{m}, b_{2}^{m},\ldots,b_{2^{m}}^{m}\}$
 of $\mathbbm{1}\in\mathbb{B}$ with $\mathbbm{1}=b^1_1\vee b^1_2$,
 $\nu(b^1_1)=\nu(b^1_2)$, and
 $b_{j}^{m}=b_{2j-1}^{m+1}\lor b_{2j}^{m+1}$,
 $\nu(b^{m+1}_{2j-1})=\nu(b^{m+1}_{2j})$, for all
 $m\in\mathbb{N}$ and $j\in\{1,2,\ldots,2^{m}\}$.
 Since $\nu(b^m_j)\to 0$ as $m\to\infty$ for each $j$,
 there is no partition refined from $(P_m)_{m\in\mathbb{N}}$.
 It remains to refer to 1.3.4\,(1),\,(3).

 \subsec{5.3.4.} Let $L^0_{\mathbbm{C}}(\Omega,\Sigma,\mu)$ be the space of
 all (cosets of) measurable complex valued
 functions, and let $L^\infty_{\mathbbm{C}}(\Omega,\Sigma,\mu)$
 be the space of essentially bounded measurable complex valued
 functions. Then the space $L^\infty_{\mathbbm{C}}(\Omega,\Sigma,\mu)$
 is isomorphic to some $C(Q,\mathbbm{C})$; consequently,
 there are no nontrivial derivations on it. If the Boolean
 algebra $\mathbb{B}(\Omega,\Sigma,\mu)$ of measurable sets
 modulo sets of measure zero is not atomic (and therefore is not
 $\sigma$-distributive, cp.~5.3.3); then, by 3.3.6\,(4),
 there exist nontrivial derivations
 on~$L^0_{\mathbb C}(\Omega,\Sigma,\mu)$
 (cp.~\cite{BCS, Kus3, Kus2}). The same is true about the spaces
 $L^\infty(\Omega,\Sigma,\mu)$ and $L^0(\Omega,\Sigma,\mu)$
 of real valued measurable functions. Moreover:

 \subsec{5.3.5.} A derivation (an automorphism) $S$ on $G$ is
 \textit{essentially nontrivial\/} provided that $\pi S=0$ ($\pi S=\pi I_G$)
 implies $\pi=0$ for every band projection $\pi\in\mathfrak{P}(G)$.
 If $(\Omega,\Sigma,\mu)$ is an~atomless measure space
 with the direct sum property then
 (cp. \cite{Kus3})

 \subsubsec{(1)}~There is an essentially nontrivial
 derivation on  $L^0_{\mathbb{R}}(\Omega,\Sigma,\mu)$;

 \subsubsec{(2)}~There is an essentially nontrivial
 derivation on $L^0_{\mathbb{C}}(\Omega,\Sigma,\mu)$;

 \subsubsec{(3)}~The identity operator is the unique automorphism of
 $L^0_{\mathbb{R}}(\Omega,\Sigma,\mu)$;

 \subsubsec{(4)}~There is an essentially nontrivial band
 preserving automorphism of $L^0_{\mathbb{C}}(\Omega,\Sigma,\mu)$.

 Also there exists an essentially nontrivial separately
 band preserving antisymmetric bilinear operator in
 $L^0_{\mathbb{R}}(\Omega,\Sigma,\mu)$, cp. \cite{Kus4}.

 \subsec{5.3.6.} Two arbitrary transcendence bases for a~field
 over a~subfield have the same cardinality called the
 {\it transcendence degree} (cp.~\cite[Chapter~II, Theorem~25]{ZS}).
 Let $\tau(\mathcal C)$ be the transcendence degree of
 $\mathcal C$ over~${\mathbb{C}}^{\scriptscriptstyle\wedge}$
 inside~$\mathbbm{V}^{(\mathbbm{B})}$. The Boolean valued cardinal
 $\tau(\mathcal C)$ carries some information on the connection
 between the Boolean algebra $\mathbb{B}$ and the complexes
 $\mathcal{C}$. Each Boolean valued cardinal is a~mixing of
 standard cardinals; i.e., the representation
 $\tau(\mathcal C)=
 \operatorname{mix}_\xi b_\xi\alpha_\xi^{\scriptscriptstyle\wedge}$
 holds, where $(b_\xi)$ is a~partition of unity of~$\mathbb B$ and $(\alpha_\xi)$ is some
 family of cardinals \big(cp.~A36\,(3) and A3.8\,(1)\big). Moreover, for
 $\mathbbm{B}_\xi\assign [\mathbbm{0},b_\xi]$ we have
 $
 \mathbbm{V}^{(\mathbbm{B}_\xi)}\models\tau(\mathcal C)=
 \alpha_\xi^{\scriptscriptstyle\wedge}.
 $
 In this connection, it would be interesting to characterize
 the complete Boolean algebras~$\mathbb B$ such that
 $\tau(\mathcal C)=\alpha^{\scriptscriptstyle\wedge}$ inside~$\mathbbm{V}^{(\mathbbm{B})}$ for
 some cardinal~$\alpha$.

 \subsec{5.3.7.} Given  $\mathcal E\subset G$, denote
 by~$\langle\mathcal E\rangle$ the set of elements of the form
 $e_1^{n_1}\cdot\ldots\cdot e_k^{n_k}$, where $e_1,\ldots, e_k\in
 \mathcal E$ and $k,n_1,\ldots,n_k\in\mathbb N$. A~set $\mathcal
 E\subset G$ is  {\it locally algebraically independent\/} provided that
 $\langle\mathcal E\rangle$ is locally linearly independent in the sense
 of 2.2.2. This
 notion, presenting the external interpretation of the internal notion
 of algebraic independence (or transcendence), seems to turn out useful in
 studying the descents of fields \cite[Section 8.3]{IBA} or
 general regular rings~\cite{Go}.

 \subsec{5.3.8.} Consider a~band preserving operator $S:\mathcal
 C{\downarrow}\to\mathcal C{\downarrow}$ satisfying the Cauchy
 functional equation $S(u+v)=S(u)S(v)$ for all
 $u,v\in\mathcal C{\downarrow}$. If, in addition, $S$ satisfies the
 condition $S(\lambda u)=S(u)^\lambda$ for arbitrary
 $\lambda\in\mathbb C$ and $u\in\mathcal C{\downarrow}$ then we say
 that~$S$ is {\it exponential}. Say that $S$ is order bounded if
 $S$ takes order bounded sets into order bounded sets. If $\sigma$
 is the ascent of~$S$ then $\sigma$ is exponential inside~$\mathbb
 V^{(\mathbb B)}$; therefore, in the class of functions bounded from
 above on a~nonzero interval, we have either $\sigma =0$ or
 $\sigma(x)=e^{cx}$ $(x\in\mathcal C)$ for some $c\in\mathcal
 C$~\cite[Chapter~5, Theorem~5]{AD}. Hence, we conclude that
 conditions~WP(1)--WP(7) of Theorem~WP are also equivalent to the
 following: {\sl every band preserving exponential operator in
 $\mathbb B(\mathbb C)\assign \mathcal C{\downarrow}$ is order bounded} ({\sl
 and consequently has the form $S=0$ or $S(x)=e^{cx}$,
 $x\in\mathcal C{\downarrow}$, for some $c\in\mathcal
 C{\downarrow})$.}

\section {Appendix. Boolean Valued Analysis}
\centerline{\scshape A1. Boolean  Valued Universes}
\medskip

We start with recalling some auxiliary facts about the
construction and  treatment of Boolean valued models.

\subsec{A1.1.}
% Universes and truth values}
Let ${\mathbbm B}$ be a~complete Boolean algebra. Given an ordinal
$\alpha$, put
$$
{\mathbbm V}_{\alpha}^{({\mathbbm B})} \assign  \bigl\{\,x\,:\,x \mbox{
is\ a\ function}\ \wedge\ (\exists\,\beta)\bigl(\beta<\alpha\ \wedge\
\dom (x)
\subset {\mathbbm V}_{\beta}^{({\mathbbm B})}\ \wedge\ \im
(x)\subset {\mathbbm B}\bigr)\bigr\}.
$$

 After this recursive definition the {\it Boolean valued
 universe\/} ${\mathbb V}^{({\mathbbm B})}$ or, in other words, the
 {\it class of ${\mathbbm B}$ valued sets\/} is introduced by
 $$
 {\mathbb V}^{({\mathbbm B})}\assign \bigcup\limits_{\alpha\in\On}
 {\mathbb V}_{\alpha}^{({\mathbbm B})},
 $$
 with $\On$ standing for the class of all ordinals.

 In case of the two element Boolean algebra $\mathbbm 2\assign \{\mathbbm
 0, \mathbbm 1\}$ this procedure yields a version of the classical
 {\it von Neumann universe\/} ${\mathbbm V}$
 (cp.~\cite[Theorem 4.2.8]{IBA}).

 Let $\varphi$ be an arbitrary formula of $\ZFC$,
 Zermelo--Fraenkel set theory with  choice. The {\it Boolean truth
 value\/} $[\![\varphi]\!]\in {\mathbbm B}$ is introduced by
 induction on the complexity of~$\varphi$ by naturally
 interpreting the propositional connectives and quantifiers in the
 Boolean algebra ${\mathbbm B}$
 \big(for instance, $[\![\varphi_1\lor\varphi_2]\!]\assign[\![\varphi_1]\!]\lor[\![\varphi_2]\!]$\big)
 and taking into consideration the
 way in which a formula is built up from atomic formulas. The
 Boolean truth values of the {\it atomic formulas\/} $x\in y$ and
 $x=y$ \big(with $x,y$ assumed to be elements of ${\mathbbm V}^{({\mathbbm B})}$\big) are defined by
 means of the following recursion schema:
 $$
 \gathered {} [\![x\in y]\!]= \bigvee\limits_{t\in\dom(y)}\!\!
 \bigl(y(t)\wedge [\![t=x]\!]\bigr),
 \\%
 [0.5\jot]
 [\![x=y]\!]=\bigvee\limits_{t\in\dom(x)}\!\!\bigl(x(t)\Rightarrow [\![t\in
 y]\!]\bigr)\wedge \bigvee\limits_{t\in\dom(y)}\!\!\bigl(y(t)\Rightarrow [\![t\in
 x]\!]\bigr).
 \endgathered
 $$
 The sign $\Rightarrow$ symbolizes the implication in ${\mathbbm
 B}$; i.e., $(a\Rightarrow b)\assign (a^\ast\vee b)$, where $a^\ast$ is as
 usual the {\it complement\/} of~$a$.
 The universe ${\mathbbm V}^{({\mathbbm B})}$ with the  Boolean
 truth value of a~formula is a~model of set theory in the sense
 that the following statement is fulfilled:

 \subsec{A1.2.} \proclaim{Transfer Principle.}For every theorem
 $\varphi$ of~ZFC, we have
 $[\![\varphi]\!]=\mathbbm 1$ (also in~ZFC);
 i.e., $\varphi$ is true inside the Boolean valued universe
 ${\mathbbm V}^{({\mathbbm B})}$.
 \endproc
 \smallskip

We enter into the next agreement: If $\varphi (x)$ is a~formula of $\ZFC$
then, on assuming $x$ to be an element of
${\mathbbm V}^{({\mathbbm B})}$, the phrase ``$x$ satisfies $\varphi$ inside
${\mathbbm V}^{({\mathbbm B})}$'' or, briefly, ``$\varphi (x)$ is
true inside ${\mathbbm V}^{({\mathbbm B})}$'' means that
$[\![\varphi(x)]\!]=\mathbbm 1$. This is sometimes written as
${\mathbbm V}^{({\mathbbm B})}\models \varphi (x)$.

Given $x\in {\mathbbm V}^{({\mathbbm B})}$ and $b\in {\mathbbm
B}$, define the function $ bx:z\mapsto b\land x (z)$ $\bigl(z\in\dom
(x)\bigr)$. Here we presume that $b\varnothing \assign \varnothing$ for all
$b\in {\mathbbm B}$.

There is a~natural equivalence relation $x\sim y\iff
[\![x=y]\!]=\mathbbm 1$ in  the class ${\mathbbm V}^{({\mathbbm
B})}$. Choosing a~representative of the smallest rank in each
equivalence class or, more exactly, using the so-called
``Frege--Russell--Scott trick,'' we obtain a~{\it separated
Boolean valued universe\/} $\overline {\mathbbm V}{}^{
({\mathbbm B})}$ for which $ x=y\iff
[\![x=y]\!]=\mathbbm 1. $

The Boolean truth value of a~formula~$\varphi$
remains unaltered if we replace in~$\varphi$ each element of~${\mathbbm
V}^{({\mathbbm B})}$ by one of its equivalents. In this connection
from now on we take ${\mathbbm V}^{({\mathbbm B})}\assign \overline
{\mathbbm V}{}^{({\mathbbm B})}$ without further specification.

Observe that in $\overline {\mathbbm V}{}^{({\mathbbm B})}$ the
element $bx$ is defined correctly for $x\in\overline {\mathbbm
V}{}^{({\mathbbm B})}$ and $b\in {\mathbbm B}$, since
 $[\![x_1=\nobreak x_2]\!]=\nobreak{\mathbbm 1}$ implies $[\![bx_1=b x_2]\!]={\mathbbm 1}$.
%For a~similar reason, we often write ${\mathbbm 0}\assign \varnothing$,
%and in particular ${\mathbbm 0}\varnothing =\varnothing={\mathbbm
%0} x$ for $x\in {\mathbbm V}^{({\mathbbm B})}$.

\subsec{A1.3.} \proclaim{Mixing Principle.} Let
$(b_{\xi})_{\xi\in\Xi}$ be a~{\it partition of unity\/} in~%
${\mathbbm B}$, i.e., $\sup_{\xi\in\Xi} b_{\xi}=\mathbbm 1$ and $\xi\neq\eta\ra b_{\xi}\wedge
b_{\eta}=\mathbbm 0$. To each family $(x_{\xi})_{\xi\in\Xi}$
in~${\mathbbm V}^{({\mathbbm B})}$ there exists a~unique element
$x$ in the separated universe such that $[\![x=x_{\xi}]\!]\ge
b_{\xi}$ for all $\xi\in\Xi$.
\endproc
\smallskip

This element $x$ is called the {\it mixing\/} of
$(x_{\xi})_{\xi\in\Xi}$ by~$(b_{\xi})_{\xi\in\Xi}$ and is denoted
by $\mix_{\xi\in\Xi} b_{\xi} x_{\xi}$.

\subsec{A1.4.} \proclaim{Maximum Principle.} Let $\varphi(x)$ be
a~formula of~ZFC. Then (in~ZFC) there is a ${\mathbbm B}$~valued set
$x_0$ satisfying $[\![(\exists\,x)\varphi (x)]\!]=[\![\varphi
(x_0)]\!]. $
\endproc

%\newpage

\section{A2.  Escher Rules}

Boolean valued analysis  consists primarily in comparison of the
instances of a~mathematical object or idea in two  Boolean valued
models. This is impossible to achieve without some dialog  between
the universes~${\mathbbm V}$ and~${\mathbbm V}^{({\mathbbm B})}$.
In other words, we need a~smooth mathematical toolkit for
revealing interplay between the interpretations of one and the
same fact in the two models~${\mathbbm V}$ and~${\mathbbm
V}^{({\mathbbm B})}$. The relevant {\it ascending-and-descending
technique\/} rests on the functors of canonical embedding,
descent, and ascent.

\subsec{A2.1.}
%Embedding}
We start with the canonical embedding of the von Neumann universe
$\mathbbm V$.

Given ${x\in {\mathbbm V}}$, we denote
by~$x^{\scriptscriptstyle\wedge}$ the {\it standard name\/} of~$x$
in~${\mathbbm V}^{({\mathbbm B})}$; i.e., the element defined by
the following recursion schema: $
\varnothing^{\scriptscriptstyle\wedge}\assign \varnothing$,\quad $\dom
(x^{\scriptscriptstyle\wedge})\assign  \{y^{\scriptscriptstyle\wedge} :
y\in x \}$,\quad $\im
(x^{\scriptscriptstyle\wedge})\assign \{{\mathbbm 1}\}$. Observe some
properties of the mapping $x\mapsto x^{\scriptscriptstyle\wedge}$
we need in the sequel.

{\subsubsec(1)} For an~arbitrary formula $\varphi(y)$ of~$\ZFC$ we
have (in ZFC) for each ${x\in {\mathbbm V}}$
$$
\gathered {}[\![(\exists\,y\in
x^{\scriptscriptstyle\wedge})\,\varphi(y)]\!]=
\bigvee\limits_{z\in x}
[\![\varphi(z^{\scriptscriptstyle\wedge})]\!],
\\
[\![(\forall\,y\in x^{\scriptscriptstyle\wedge})\,\varphi(y)]\!]
=\bigwedge\limits_{z\in
x}[\![\varphi(z^{\scriptscriptstyle\wedge})]\!].
\endgathered
$$

{\subsubsec(2)} If $x,y\in{\mathbbm V}$ then, by transfinite
induction, we establish $x\in y\iff {\mathbbm V}^{({\mathbbm
B})}\models x^{\scriptscriptstyle\wedge}\in
y^{\scriptscriptstyle\wedge}$, $x=\nobreak y\iff {\mathbbm
V}^{({\mathbbm B})}\models
x^{\scriptscriptstyle\wedge}=y^{\scriptscriptstyle\wedge}$. In
other words, the standard name can be considered as an~embedding
of~${\mathbbm V}$ into~${\mathbbm V}^{({\mathbbm B})}$. Moreover,
it is beyond a~doubt that the standard name sends ${\mathbbm V}$
onto ${\mathbbm V}^{(\mathbb 2)}$, which fact is demonstrated by
the next proposition:

\smallskip
{\subsubsec(3)} The following holds: $ (\forall\, u\in {\mathbbm
V}^{(\mathbb 2)})\,(\exists!\,x\in {\mathbbm V})\ {\mathbbm
V}^{({\mathbbm B})}\models u=x^{\scriptscriptstyle\wedge}. $
\smallskip

A~formula is  {\it bounded\/} or {\it restricted\/} provided that
each bound variable in it is restricted by a~bounded quantifier;
i.e., a~quantifier ranging over a~particular set. The latter means
that each bound variable~$x$ is restricted by a~quantifier of the
form~$(\forall\,x\in y)$ or $(\exists\,x\in y)$.

\subsec{A2.2.} \proclaim{Restricted Transfer Principle.}Let $\varphi(x_1,\dots,x_n)$
be a bounded formula of~ZFC. Then (in ZFC) for every
collection~$x_1,\dots,x_n\in {\mathbbm V}$ we have \
$\varphi(x_1,\dots,x_n)\iff {\mathbbm V}^{({\mathbbm
B})}\models
\varphi(x_1^{\scriptscriptstyle\wedge},\dots,x_n^{\scriptscriptstyle\wedge})$.
\endproc
\smallskip

Henceforth, working in the separated universe~$\overline{{\mathbbm
V}}{}^{({\mathbbm B})}$, we agree to preserve the
symbol~$x^{\scriptscriptstyle\wedge}$ for the distinguished
element of the class corresponding to~$x$.

Observe for example that
 the restricted transfer principle yields:
$$
\gathered \text{``$\Phi$ is a correspondence from $x$ into $y$''}
\ \iff
\\%
[-4pt]
{\mathbbm V}^{({\mathbbm B})}\models
\text{``$\Phi^{\scriptscriptstyle\wedge}$ is a correspondence from
$x^{\scriptscriptstyle\wedge}$ into
$y^{\scriptscriptstyle\wedge}$'';}
\\
\text{``$f: x\to y$'' $\ \iff\ $ ${\mathbbm V}^{({\mathbbm
B})}\models$ ``$f^{\scriptscriptstyle\wedge}:
x^{\scriptscriptstyle\wedge}\to y^{\scriptscriptstyle\wedge}$''}
\endgathered
$$
\noindent \big(moreover,
$f(a)^{\scriptscriptstyle\wedge}=f^{\scriptscriptstyle\wedge}(a^{\scriptscriptstyle\wedge})$
for all $a\in x$\big). Thus, the standard name can be considered as a
covariant functor from the category of sets (or correspondences)
inside~${\mathbbm V}$ to an~appropriate subcategory of~${\mathbbm
V}^{(\mathbb 2)}$ in the separated universe~${\mathbbm
V}^{({\mathbbm B})}$.

\subsec{A2.3.}~A~set $X$ is  {\it finite\/} provided that $X$
coincides with the image of a~function on a~finite ordinal. In
symbols, this is expressed as $\Fin (X)$; hence,
$$
\Fin(X)\assign (\exists\,n)(\exists\,f)\bigl(n\in\omega\ \wedge\ f \mbox{ is a
function}\ \wedge\ \dom(f)=n\ \wedge\ \im(f)=X\bigr)
$$
(as usual $\omega\assign \{0,1,2,\dots\}$). Obviously, the above formula
is  not bounded. Nevertheless there is a~simple transformation
rule for the  class of finite sets under the canonical embedding.
Denote by ${\mathcal P}_{\Fin}(X)$ the class of all finite subsets
of $X$; i.e., ${\mathcal P}_{\Fin}(X)\assign \{Y\in{\mathcal P}(X)
:\Fin (Y)\}. $ For an arbitrary set $X$ the following holds: $
{\mathbbm V} ^{({\mathbbm B})}\models\mathcal P _{\Fin}
(X)^{\scriptscriptstyle\wedge} = \mathcal P _{\Fin}
(X^{\scriptscriptstyle\wedge}). $

\subsec{A2.4.}
%. Descent}
Given an~arbitrary element~$x$ of the (separated) Boolean valued
universe~${\mathbbm V}^{({\mathbbm B})}$, we define the {\it
descent\/} $x{\downarrow}$ of~$x$ as $ x{\downarrow}\assign \{y\in
{\mathbbm V}^{({\mathbbm B})}: [\![y\in x]\!]={\mathbbm 1} \}$. We
list the simplest properties of descending:

\smallskip
{\subsubsec(1)} The class~$x{\downarrow}$ is a~set, i.e.,
$x{\downarrow}\in {\mathbbm V}$ for all~$x\in {\mathbbm
V}^{({\mathbbm B})}$. If $[\![x\ne \varnothing]\!]={\mathbbm 1}$
then $x{\downarrow}$ is nonempty.

{\subsubsec(2)} Let $\varphi(x)$ be a formula of~ZFC. Then (in ZFC)
for every $z\in {\mathbbm V}^{({\mathbbm B})}$ such that
$[\![z\ne \varnothing]\!]={\mathbbm 1}$ we have
$$
\gathered {} [\![(\forall\,x\in z)\,\varphi(x)]\!]
=\bigwedge\limits_{ x\in z{\downarrow}}[\![\varphi(x)]\!],
\\
[\![(\exists\,x\in z)\,\varphi(x)]\!] =\bigvee\limits_{x\in
z{\downarrow}}[\![\varphi(x)]\!].
\endgathered
$$
Moreover, there exists $x_0\in z{\downarrow}$ such that
$[\![(\exists\,x\in z)\,\varphi(x)]\!]=[\![\varphi(x_0)]\!]$.

{\subsubsec(3)} Let $\Phi$ be a~correspondence from~$X$ into~$Y$ in
${\mathbbm V}^{({\mathbbm B})}$. Thus, $\Phi$, $X$, and~$Y$ are
elements of~${\mathbbm V}^{({\mathbbm B})}$ and, moreover,
$[\![\Phi\subset X\times Y]\!]={\mathbbm 1}$. There is a~unique
correspondence~$\Phi{\downarrow}$ from $X{\downarrow}$
into~$Y{\downarrow}$ such that $
\Phi{\downarrow}(A{\downarrow})=\Phi(A){\downarrow} $ for every
nonempty subset~$A$ of~$X$ inside~${\mathbbm V}^{({\mathbbm B})}$.
The correspondence~$\Phi{\downarrow}$ is called the {\it descent\/} of $\Phi$.

{\subsubsec(4)} The descent of the composite of correspondences
inside~${\mathbbm V}^{({\mathbbm B})}$ is the composite of their
descents:
$(\Psi\circ\Phi){\downarrow}=\Psi{\downarrow}\circ\Phi{\downarrow}.
$

{\subsubsec(5)} If $\Phi$ is a~correspondence inside~${\mathbbm
V}^{({\mathbbm B})}$ then $
(\Phi^{-1}){\downarrow}=(\Phi{\downarrow})^{-1}. $

 {\subsubsec(6)} Let $\Id_X$ be the identity mapping
 inside~${\mathbbm V}^{({\mathbbm B})}$ of a~set~$X\in {\mathbbm
 V}^{({\mathbbm B})}$. Then
 $({\Id}_X){\downarrow}={\Id}_{X{\downarrow}}$.

 {\subsubsec(7)} Suppose that $X,Y,f\in {\mathbbm V}^{({\mathbbm
 B})}$ are such that $[\![f:X\to Y]\!]={\mathbbm 1}$, i.e., $f$ is
 a~mapping from~$X$ into~$Y$ inside~${\mathbbm V}^{({\mathbbm B})}$.
 Then $f{\downarrow}$ is a~unique mapping from~$X{\downarrow}$
 into~$Y{\downarrow}$ satisfying
 $[\![f{\downarrow}(x)=f(x)]\!]={\mathbbm 1}$ for all $x\in
 X{\downarrow}$. The descent of a mapping is
 \textit{extensional\/}:
 $[\![x_1=x_2]\!]\leq[\![f{\downarrow}(x_1)=f{\downarrow}(x_2)]\!]$
 for all $x_1,x_2\in X{\downarrow}$ \big(cp.~A2.5\,(4)\big).

\smallskip
By virtue of~\hbox{(1)--(7)}, we can consider the descent
operation as a~functor from the category of~${\mathbbm B}$~valued
sets and mappings (correspondences) to the category of the
standard sets and mappings (correspondences) (i.e., those in the sense
of~${\mathbbm V}$).
\smallskip

{\subsubsec(8)} Given $x_1,\dots,x_n\in {\mathbbm V}^{({\mathbbm
B})}$, denote by $(x_1,\dots,x_n)^{\mathbbm B}$ the corresponding
ordered $n$-tuple inside~${\mathbbm V}^{({\mathbbm B})}$. Assume
that $P$ is an~$n$-ary relation on~$X$ inside~${\mathbbm
V}^{({\mathbbm B})}$; i.e., $X,P\in {\mathbbm V}^{({\mathbbm B})}$
and $[\![P\subset X^n]\!]=\nobreak{\mathbbm 1}$.
Then there exists an~$n$-ary
relation~$P'$ on~$X{\downarrow}$ such that $(x_1,\dots,x_n)\in
P'\iff [\![(x_1,\dots,x_n)^{\mathbbm B}\in\nobreak
P]\!]=\nobreak{\mathbbm 1}$. Slightly abusing notation, we denote~$P'$ by
the occupied symbol~$P{\downarrow}$ and call~$P{\downarrow}$ the
{\it descent\/} of~$P$.

\subsec{A2.5.}
%. Ascent}
Let $x\in {\mathbbm V}$ and $x\subset {\mathbbm V}^{({\mathbbm
B})}$; i.e., let $x$ be some set composed of ${\mathbbm B}$ valued
sets or, in other words, $x\in \mathcal P({\mathbbm V}^{({\mathbbm
B})})$. Put $\varnothing{\uparrow}\assign \varnothing$ and $ \dom
(x{\uparrow})\assign x$, $\im (x{\uparrow})\assign \{{\mathbbm 1}\} $ if
$x\neq \varnothing$. The element~$x{\uparrow}$ \big(of the separated
universe~${\mathbbm V}^{({\mathbbm B})}$, i.e., the distinguished
representative of the class~$\{y\in {\mathbbm V}^{({\mathbbm B})}:
[\![y=\nobreak x{\uparrow}]\!]=\nobreak{\mathbbm 1}
 \}$\big)
is  the {\it ascent\/} of~$x$.

\smallskip
{\subsubsec(1)} Let $\varphi(y)$ be a~formula of~ZFC. Then (in~ZFC)
for all $x\in \mathcal P({\mathbbm V}^{({\mathbbm
B})})$ we have
$$
\gathered {}[\![(\forall\,y\in x{\uparrow})\,\varphi(y)]\!]=
\bigwedge_{y\in x} [\![\varphi (y)]\!],
\\
[\![(\exists\,y\in x{\uparrow})\,\varphi(y)]\!]= \bigvee_{y\in x}
[\![\varphi (y)]\!].
\endgathered
$$

Introducing the ascent of a~correspondence~$\Phi\subset X\times
Y$, we have to bear in mind a~possible distinction between the
domain of departure,~$X$, and the domain, $\dom (\Phi)\assign \{x\in
X:\Phi(x)\ne \varnothing \}$. This circumstance is immaterial for
the sequel; therefore,  speaking of ascents, we always imply total
correspondences; i.e., $\dom (\Phi)=X$.

{\subsubsec(2)} Let $X,Y,\Phi\in {\mathbbm V}^{({\mathbbm B})}$
and let $\Phi$ be a~correspondence from~$X$ into~$Y$. There exists
a~(unique) correspondence~$\Phi{\uparrow}$ from~$X{\uparrow}$
into~$Y{\uparrow}$ inside~${\mathbbm V}^{({\mathbbm B})}$, such that
$ \Phi{\uparrow}(A{\uparrow})=\Phi(A){\uparrow} $ is valid for
every subset~$A$ of~$\dom (\Phi)$, if and only if $\Phi$ is {\it
extensional\/}; i.e., satisfies the condition $ y_1\in \Phi(x_1)\to
[\![x_1=\nobreak x_2]\!]\le \bigvee\nolimits_{y_2\in \Phi(x_2)}
[\![y_1=y_2]\!] $ for $x_1,x_2\in \dom (\Phi)$. In this event,
$\Phi{\uparrow}=\Phi'{\uparrow}$, where $\Phi'\assign \{(x,y)^{\mathbbm
B}: (x,y)\in \Phi \}$. The element $\Phi{\uparrow}$ is  the {\it
ascent\/} of the initial correspondence~$\Phi$.

{\subsubsec(3)} The composite of extensional correspondences is
extensional. Moreover, the ascent of a~composite is equal to the
composite of the ascents inside~${\mathbbm V}^{({\mathbbm B})}$:
On assuming that $\dom (\Psi)\supset \im (\Phi)$ we have $
{\mathbbm V}^{({\mathbbm
B})}\vDash(\Psi\circ\Phi){\uparrow}=\Psi{\uparrow}\circ\Phi{\uparrow}.
$

Note that if $\Phi$ and $\Phi^{-1}$ are extensional then
$(\Phi{\uparrow})^{-1}=(\Phi^{-1}){\uparrow}$. However, in
general, the extensionality of~$\Phi$ in no way guarantees the
extensionality of~$\Phi^{-1}$.

{\subsubsec(4)} It is worth mentioning that if an~extensional
correspondence~$f$ is a~function from~$X$ into~$Y$ then the
ascent~$f{\uparrow}$ of~$f$ is a~function from~$X{\uparrow}$
into~$Y{\uparrow}$. Moreover, the extensionality property can be
stated as follows: $[\![x_1=x_2]\!]\le [\![f(x_1)=f(x_2)]\!]$  for
all $x_1,x_2\in X$.

\subsec{A2.6.} Given a~set $X\subset {\mathbbm V}^{({\mathbbm
B})}$, we denote by~$\mix(X)$ the set of all mixings of the
form~$\mix_\xi b_{\xi}x_{\xi}$, where $(x_{\xi})\subset X$ and
$(b_{\xi})$ is an~arbitrary partition of unity. The following
propositions are referred to as the {\it arrow cancellation
rules\/} or {\it ascending-and-descending rules}. There are many
good reasons to call them simply
the {\it Escher rules}.%~\cite{Hof}.

\smallskip
{\subsubsec(1)} Let $X$ and $X'$ be subsets of~${\mathbbm
V}^{({\mathbbm B})}$ and let $f:X\to X'$ be an~extensional
mapping. Suppose also that $Y,Y',g\in {\mathbbm V}^{({\mathbbm B})}$
are such that $[\![\,Y\ne \varnothing]\!]=[\![\,g:Y\to
Y']\!]={\mathbbm 1}$. Then
$X{\uparrow}{\downarrow}=\mix(X)$, $Y{\downarrow}{\uparrow}=Y$,
$f{\uparrow}{\downarrow}=f$ on $X$, and $g{\downarrow}{\uparrow}=g$.

{\subsubsec(2)}~If $X$ is a subset of ${\mathbbm V}^{({\mathbbm B})}$ then
${\mathbbm V}^{({\mathbbm B})}\models
{\mathcal P}_{\Fin}(X{\uparrow})=\{{\theta{\uparrow}} :
\theta\in{\mathcal P}_{\Fin}(X)\}{\uparrow}.
 $
\smallskip

%\section{The technique of ascending and descending}
Suppose that $X\in {\mathbbm V}$, $X\ne\varnothing$; i.e., $X$ is
a~nonempty set. Let the letter~$\iota$ denote the standard name
embedding $x\mapsto x^{\scriptscriptstyle\wedge}$ $(x\in X)$. Then
$\iota(X){\uparrow}=X^{\scriptscriptstyle\wedge}$ and
$X=\iota^{-1}(X^{\scriptscriptstyle\wedge}{\downarrow})$. Using
the above relations, we may extend the descent and ascent
operations to the case in which $\Phi$ is a~correspondence
from~$X$ into~$Y{\downarrow}$ and $[\![\,\Psi$ is a~correspondence
from~$X^{\scriptscriptstyle\wedge}$ into~$Y\,]\!]={\mathbbm 1}$, where
$Y\in {\mathbbm V}^{({\mathbbm B})}$. Namely, we put
$\Phi\upwardarrow\assign (\Phi\circ\iota^{-1}){\uparrow}$ and
$\Psi\downwardarrow\assign \Psi{\downarrow}\circ\iota$. In this case,
$\Phi\upwardarrow$ is the {\it modified ascent\/} of~$\Phi$ and
$\Psi\downwardarrow$ is  the {\it modified descent\/} of~$\Psi$.
(If the context excludes ambiguity then we briefly speak of
ascents and descents using simple arrows.) It is easy to see that
$\Phi\upwardarrow$ is a~unique correspondence inside~${\mathbbm
V}^{({\mathbbm B})}$ satisfying the relation $
[\![\Phi\upwardarrow(x^{\scriptscriptstyle\wedge})=\Phi(x){\uparrow}]\!]={\mathbbm
1}$ $(x\in X)$. Similarly, $\Psi\downwardarrow$ is a~unique
correspondence from~$X$ into~$Y{\downarrow}$ satisfying the equality
$
\Psi\downwardarrow(x)=\Psi(x^{\scriptscriptstyle\wedge}){\downarrow}$
$(x\in X)$. If $\Phi\assign f$ and~$\Psi\assign g$ are functions then these
relations take the form $
[\![f\upwardarrow(x^{\scriptscriptstyle\wedge})=f(x)]\!]={\mathbbm
1}$ and $g\downwardarrow(x)=g(x^{\scriptscriptstyle\wedge})$ for
all $x\in X$.

\subsec{A2.7.}
%Functional Representation of  Boolean Valued Universes}
Various function spaces reside in functional analysis, and so
the~problem is natural of replacing an~abstract Boolean valued
system by some function-space analog, a~model whose elements are
functions and in which the~basic logical operations are calculated
``pointwise.'' An~example of such a~model is given by
the~class~${\mathbbm V}^Q$ of all functions defined on a~fixed
nonempty set~$Q$ and acting into~${\mathbbm V}$. The truth values
on~${\mathbbm V}^Q$ are various subsets of~$Q$: The~truth value
$[\![\varphi(x_1,\dots,x_n)]\!]$ of a formula $\varphi(x_1,\dots,x_n)$ (at
functions $x_1,\dots,x_n\in {\mathbbm V}^Q$) is calculated as
follows:
$$
[\![\varphi(x_1,\dots,x_n)]\!]= \big\{q\in Q :
 \varphi\big(x_1(q),\dots,x_n(q)\big)\big\}.
$$

Gutman and Losenkov solved  the~above problem by the concept of
continuous polyverse which is a~continuous bundle of models of set
theory. It~is shown that the~class of continuous sections of
a~continuous polyverse is a~Boolean valued system satisfying all
basic principles of Boolean valued analysis and, conversely, each
Boolean valued algebraic system can be represented as the~class of
sections of a~suitable continuous polyverse. More details reside
in~\cite[Chapter~6]{IBA}.

\subsec{A2.8.}
% Analysis of Algebraic Systems}
Every Boolean valued universe has  the collection of mathematical
objects in full supply:  available in plenty are all sets with
extra structure (groups, rings, algebras, normed spaces, etc.).
Applying the descent functor to such {\it internal\/} algebraic
systems of a~Boolean valued model, we distinguish some bizarre
entities or recognize old acquaintances, which leads  to revealing
the new facts of their life and structure.

This technique of research, known as {\it direct Boolean valued
interpretation}, allows us to produce new theorems or, to be more
exact, to extend the semantical content of the available theorems
by means of  slavish translation.  The information we so acquire
might fail to be vital, valuable, or intriguing, in which case the
direct Boolean valued interpretation  turns out into a leisurely
game.

It thus stands to reason to raise  the following questions: What
structures significant for mathematical practice are obtainable by
the Boolean valued interpretation of the most typical algebraic
systems?  What transfer principles hold true in this process?
Clearly, the answers should imply specific objects whose
particular  features enable us to deal with their Boolean valued
representation which, if understood duly, is impossible to
implement for arbitrary algebraic systems.

An~{\it abstract\/ ${\mathbbm B}$-set\/} or~{\it set with ${\mathbbm
B}$-structure\/} is a~pair $(X,d)$, where $X\in {\mathbbm V}$,
$X\ne \varnothing$, and $d$ is a~mapping from~$X\times X$
into~${\mathbbm B}$ such that \ $d (x,y)={{\mathbbm 0}}\iff
x=y$; \ $d (x,y)=d (y,x)$; \ $d (x,y)\le d (x,z)\vee d (z,y)$ for all
$x,y,z\in X$.

To obtain an easy~example of an abstract ${\mathbbm B}$-set, given
$\varnothing\ne X\subset {\mathbbm V}^{({\mathbbm B})}$ put
$$
d (x,y)\assign [\![x\ne y]\!]= [\![x=y]\!]^\ast \quad\text{for $x,\ y\in X$.}
$$
Another easy example is a~nonempty $X$ with the {\it discrete
${\mathbbm B}$-metric} $d$; i.e., $d (x,y)={\mathbbm 1}$ if $x\ne
y$ and $d (x,y)={\mathbbm 0}$ if $x=y$.

Let $(X,d)$ be some abstract ${\mathbbm B}$-set. There exist
an~element ${\mathcal X}\in {\mathbbm V}^{({\mathbbm B})}$ and
an~injection~$\iota:X\to X'\assign {\mathcal X}{\downarrow}$ such that
$d (x,y)=[\![\iota x\ne \iota y]\!]$ for all $x, y\in X$ and
each~$x'\in X'$ admits the representation $x'=\mix_{\xi\in
\Xi}b_{\xi}\iota x_{\xi}$, where $(x_{\xi})_{\xi\in \Xi}\subset
X$ and $(b_{\xi})_{\xi\in \Xi}$ is a~partition of unity
in~${\mathbbm B}$.
We see that an abstract ${\mathbbm B}$-set $X$  embeds in
the~Boolean valued universe ${\mathbbm V}^{({\mathbbm B})}$ so
that the Boolean distance between the members of $X$ becomes the
Boolean truth value of the negation of their equality.  The
corresponding element ${\mathcal X}\in{\mathbbm V} ^{({\mathbbm B})}$ is, by
definition, the {\it Boolean valued representation\/} of~$X$.

If $X$ is a~discrete abstract ${\mathbbm B}$-set then ${\mathcal
X}=X^{\scriptscriptstyle\wedge}$ and $\iota
x=x^{\scriptscriptstyle\wedge}$ for all $x\in X$. If $X\subset
{\mathbbm V}^{({\mathbbm B})}$ then $\iota{\uparrow}$ is
an~injection of~$X{\uparrow}$ into~${\mathcal X}$
(inside~${\mathbbm V}^{({\mathbbm B})}$). A~mapping~$f$ from
a~${\mathbbm B}$-set $(X,d)$ into a~${\mathbbm B}$-set $(X',d')$ is
said to be {\it contractive\/} if $d (x,y)\ge d'\bigl(f(x),f(y)\bigr)$ for
all $x,y\in X$.

In case a~${\mathbbm B}$-set $X$ has some  a~priori structure we
may try to furnish the Boolean valued representation of $X$ with
an analogous structure, so as to apply the technique of ascending
and descending to the study of the original structure of~$X$.
Consequently, the above questions
 may be treated as instances of the~unique problem
of  searching a~well-qualified Boolean valued representation of
a~${\mathbbm B}$-set  with some additional structure.

We call these objects {\it algebraic ${\mathbbm B}$-systems\/}.
Located at the epicenter of Boolean valued analysis, the notion of
an algebraic ${\mathbbm B}$-system refers to a~nonempty ${\mathbbm
B}$-set endowed with a~few contractive operations  and ${\mathbbm
B}$-predicates, the latter  meaning  ${\mathbbm B}$ valued
contractive mappings.

The Boolean valued representation of an algebraic ${\mathbbm
B}$-system appears to be a~standard two valued algebraic system of
the same type.  This means that an appropriate completion of each
algebraic ${\mathbbm B}$-system coincides with the descent of
some~two valued algebraic system inside ${\mathbbm V} ^{({\mathbbm
B})}$.

On the other hand, each two valued algebraic system  may be
transformed into an algebraic ${\mathbbm B}$-system on
distinguishing a~complete Boolean algebra of congruences of the
original system.  In this event, the task is in order of finding
the formulas  holding true in direct or reverse transition from
a~${\mathbbm B}$-system to a~two valued system.  In other words,
we  have to seek and reveal here  some versions of transfer in the
form of identity preservation, a principle  of long standing in
vector lattice theory.

%\newpage
\section{A3. Boolean Valued Numbers, Ordinals, and Cardinals}

Boolean valued analysis stems from the fact that each internal
field of reals of a~Boolean valued model descends into
a~universally complete vector lattice. Thus, a~remarkable
opportunity opens up to expand and enrich the treasure-trove of
mathematical knowledge by translating information about the reals
to the language of other noble families of functional analysis. We
will elaborate  upon the matter in  this section.

\subsec{A3.1.} Recall a~few definitions. Two elements $x$ and $y$
of a~vector lattice $E$ are {\it disjoint\/} (in~symbols $x\perp
y$) provided that $|x| \wedge  |y| =0$. A~{\it band\/} of $E$ is
defined as the {\it disjoint complement\/} $M^\perp\assign  \{x\in E
:\, (\forall\,y\in M)\, x\perp y\}$ of a~nonempty set $M\subset E$.

The inclusion-ordered set ${\mathfrak B}(E)$ of all bands in $E$
is a~complete Boolean algebra with the Boolean operations:
$$
L\wedge K=L\cap K,\quad L\vee K=(L\cup K)^{\perp\perp},\quad L^*
=L^\perp\quad \bigl(L,K\in{\mathfrak B}(E)\bigr).
$$
The Boolean algebra ${\mathfrak B}(E)$ is often referred to as the
{\it base\/} of~$E$.

A~{\it band projection\/} in $E$ is a~linear idempotent operator
in $\pi:E\to E$ satisfying the inequalities $0\leq\pi x\leq x$ for
all $0\leq x\in E$. The set ${\mathfrak P}(E)$ of all band
projections ordered by
$\pi\le\rho\iff\pi\circ\rho=\pi$ is a~Boolean
algebra with the Boolean operations:
$$
\pi\wedge \rho =\pi\circ\rho,\quad \pi\vee\rho =\pi +\rho
-\pi\circ\rho,\quad \pi^* =I_E-\pi\quad \bigl(\pi,\rho\in\mathfrak P
(E)\bigr).
$$

Let $u\in E_+$ and $e\wedge (u-e)= 0$ for some $0\leq e\in E$.
Then $e$ is  a~{\it fragment\/} or {\it component\/} of~$u$. The
set ${\mathfrak E}(u)$ of all fragments of $u$ with the order
induced by~$E$ is a~Boolean algebra where the lattice operations
are taken from~$E$ and the Boolean complement has the form $e^*
\assign  u - e$.

\subsec{A3.2.} A Dedekind complete vector lattice is also called
a~{\it Kantorovich space\/} or $K$-{\it space\/}, for short. A
Dedekind complete vector lattice $E$ is {\it universally
complete\/} if every family of pairwise disjoint elements of~$E$
is order bounded.

\subsubsec{(1)}\proclaim{}Let $E$ be an arbitrary
\hbox{$K$-}space. Then the correspondence $\pi\mapsto\pi(E)$
determines an~isomorphism of the Boolean algebras~${\mathfrak
P}(E)$ and~${\mathfrak B}(E)$. If there is an~order unity
${\mathbbm 1}$ in~$E$ then the mappings $\pi\mapsto\pi{\mathbbm
1}$ from~${\mathfrak P}(E)$ into~${\mathfrak E}({\mathbbm 1})$ and $e\mapsto
\{e\}^{\perp\perp}$ from~${\mathfrak E}({\mathbbm 1})$ into~${\mathfrak
B}(E)$ are isomorphisms of Boolean algebras too.
\endproc

\subsubsec{(2)} \proclaim{}Each universally complete vector
lattice $E$ with order unity~${\mathbbm 1}$ can be uniquely
endowed with multiplication so as to make~$E$ into a~faithful
$f$-algebra and ${\mathbbm 1}$ into a~ring unity. In this
$f$-algebra each band projection $\pi\in{\mathfrak P}(E)$ is the
operator of multiplication by~$\pi({\mathbbm 1})$.
\endproc

\subsec{A3.3.} By a~{\it field of reals\/} we mean every algebraic
system that satisfies the axioms of an~Archimedean  ordered field
(with distinct zero and unity) and enjoys the axiom of
completeness. The same object can be defined as a~one-dimensional
$K$-space.

Recall the well-known assertion of $\ZFC$: {\sl There exists
a~field of reals ${\mathbbm R}$ that is unique up to isomorphism.}

Successively applying the transfer and maximum principles, we find
an~element ${\mathcal R}\in{\mathbbm V}^{({\mathbbm B})}$ for
which $[\![\,{\mathcal R}$ is a~field of reals$\,]\!]={\mathbbm
1}$.  Moreover, if an~arbitrary ${\mathcal R}\,'\in {\mathbbm
V}^{({\mathbbm B})}$ satisfies the condition $[\![\,{\mathcal R}\,'$~%
is a~field of  reals$\,]\!]={\mathbbm 1}$ then $[\![\,$the ordered
fields $\mathcal R$ and ${\mathcal R}\,'$ are
isomorphic\,$]\!]=\nobreak{\mathbbm 1}$.  In other words, there
exists an internal field~of  reals ${\mathcal R}\in{\mathbbm
V}^{({\mathbbm B})}$ which is unique up to isomorphism.

By the same reasons there exists an internal field~of complex
numbers ${\mathcal C}\in{\mathbbm V}^{({\mathbbm B})}$ which is
unique up to isomorphism. Moreover, ${\mathbbm V}^{({\mathbbm
B})}\models{\mathcal C}={\mathcal R}\oplus i{\mathcal R}$. We call
${\mathcal R}$ and ${\mathcal C}$ the {\it internal  reals} and
{\it internal complexes\/} in ${\mathbbm V}^{({\mathbbm B})}$.

\subsec{A3.4.} Consider another well-known assertion of $\ZFC$:
{\sl If\/ ${\mathbbm P}$ is an Archimedean ordered field then there
is an~isomorphic embedding~$h$ of the field\/~${\mathbbm P}$
into~${\mathbbm R}$ such that the image~$h({\mathbbm P})$ is
a~subfield of\/~${\mathbbm R}$ containing the subfield of rational
numbers.  In particular, $h({\mathbbm P})$ is dense in~${\mathbbm
R}$.}

Note also that $\varphi (x)$, presenting the conjunction of the
axioms of an~Archimedean ordered field~$x$, is bounded; therefore,
$[\![\,\varphi ({\mathbbm R}^{\scriptscriptstyle\wedge})\,]\!]
={\mathbbm 1}$, i.e., $[\![\,{\mathbbm
R}^{\scriptscriptstyle\wedge}$ is an~Archimedean ordered
field$\,]\!]={\mathbbm 1}$. ``Pulling'' the above assertion through the
transfer principle, we conclude that $[\![\,{\mathbbm
R}^{\scriptscriptstyle\wedge}$ is isomorphic to a~dense subfield
of~$\mathcal R\,]\!]={\mathbbm 1}$. We further assume that
${\mathbbm R}^{\scriptscriptstyle\wedge}$ is a~dense subfield of
${\mathcal R}$ and ${\mathbbm C}^{\scriptscriptstyle\wedge}$ is
a~dense subfield of ${\mathcal C}$. It is easy to see that the
elements~$0^{\scriptscriptstyle\wedge}$
and~$1^{\scriptscriptstyle\wedge}$ are the zero and unity
of~${\mathcal R}$.

Observe that the equalities ${\mathcal R}={\mathbbm
R}^{\scriptscriptstyle\wedge}$ and ${\mathcal C}={\mathbbm
C}^{\scriptscriptstyle\wedge}$ are not valid in general. Indeed,
the axiom of completeness for~${\mathbbm R}$ is not a~bounded
formula and so it may fail for ${\mathbbm
R}^{\scriptscriptstyle\wedge}$ inside~${\mathbbm V}^{({\mathbbm
B})}$. (The corresponding example is given in Section~1.3 of this paper.)

\subsec{A3.5.} Look now at the descent ${\mathcal R}{\downarrow}$
of the algebraic system~${\mathcal R}$. In other words,  consider
the descent of the underlying set of the system~${\mathcal R}$
together with the descended operations and order. For simplicity, we
denote the operations and order in~${\mathcal R}$ and~${\mathcal
R}{\downarrow}$ by the same symbols $+$, $\cdot\,$, and~$\le $. In
more detail, we introduce addition, multiplication, and order
in~${\mathcal R}{\downarrow}$ by the formulas
$$
\align z=x+y\ &\iff\ [\![\,z=x+y\,]\!]={\mathbbm 1},
\\
z=x\cdot y\ &\iff\ [\![\,z=x\cdot y\,]\!]={\mathbbm 1},
\\
x\le y\ &\iff\ [\![\,x\le y\,]\!]={\mathbbm 1} \quad (x, y,
z\in {\mathcal R}{\downarrow}).
\endalign
$$
Also, we may introduce multiplication by the usual reals
in~${\mathcal R}{\downarrow}$ by the rule
$$
y=\lambda x\ \iff\
[\![\,y=\lambda^{\scriptscriptstyle\wedge}x\,]\!]={\mathbbm
1}\quad (\lambda\in{\mathbbm R},\ x, y\in {\mathcal
R}{\downarrow}).
$$

The fundamental result of Boolean valued analysis is the Gordon
Theorem which  reads as follows: {\sl Each universally complete
vector lattice is an~interpretation of the reals in an~appropriate
Boolean valued model}. Formally:

\subsec{A3.6.} \proclaim{Gordon Theorem.} Let ${\mathcal R}$ be
the reals inside~${\mathbbm V}^{({\mathbbm B})}$. Then ${\mathcal
R}{\downarrow}$, with the descended operations and order, is
a~universally complete vector lattice with order unit~$1^{\scriptscriptstyle\wedge}$.
Moreover, there exists an~isomorphism~$\chi$ of\/~${\mathbbm B}$
onto~${\mathfrak P} ({\mathcal R}{\downarrow})$ such that
$$
\chi (b) x=\chi (b) y\ \iff\ b\le [\![\,x=y\,]\!], \quad
\chi (b) x\le \chi (b) y\ \iff\ b\le [\![\,x\le y\,]\!]
$$
for all $x, y\in {\mathcal R}{\downarrow}$ and $b\in {\mathbbm
B}$.
\endproc
\smallskip

The converse  is also true: {\sl Each Archimedean vector lattice
embeds in a~Boolean valued model, becoming a~vector sublattice of
the reals (viewed as such over some dense subfield of the reals)}.

\subsec{A3.7.} \proclaim{Theorem.} Let $E$ be an~Archimedean
vector lattice, let ${\mathcal R}$ be the reals inside~${\mathbbm
V}^{({\mathbbm B})}$, and let $\jmath$ be an isomorphism
of\/~${\mathbbm B}$ onto~${\mathfrak B}(E)$. Then there is
${\mathcal E}\in {\mathbbm V}^{({\mathbbm B})}$ satisfying the
following:

\subsubsec{(1)} $\mathcal E$ is a~vector sublattice of~${\mathcal
R}$ over~${\mathbbm R}^{\scriptscriptstyle\wedge}$ inside
${\mathbbm V}^{({\mathbbm B})}$;

\subsubsec{(2)} $E'\assign {\mathcal E}{\downarrow}$ is a~vector
sublattice of~${\mathcal R}{\downarrow}$ invariant under every
band projection $\chi (b)$ $(b\in {\mathbbm B})$ and such that
each set of pairwise disjoint elements in $E'$ has a~supremum;

\subsubsec{(3)} There is an~order continuous lattice
isomorphism~$\iota:E\to E'$ such that $\iota (E)$ is a~coinitial
sublattice of~${\mathcal R}{\downarrow}$;

\subsubsec{(4)} For every~$b\in {\mathbbm B}$ the band projection
in~${\mathcal R}{\downarrow}$ onto $\{\iota
(\jmath(b))\}^{\perp\perp}$ coincides with~$\chi (b)$.
\endproc

Note also that ${\mathcal E}$ and ${\mathcal R}$ coincide if and
only if $E$ is Dedekind complete.  Thus, each theorem about the
reals within Zermelo--Fraenkel set theory has an analog in an
arbitrary Dedekind complete vector lattice. Translation of
theorems is carried out by appropriate general functors of Boolean
valued analysis. In particular, the most important structural
properties of vector lattices such as the functional
representation, spectral theorem,  etc. are the ghosts of some
properties of the reals in an~appropriate Boolean valued model.

\goodbreak

\subsec{A3.6.}
 Let us dwell for a while on the properties of
 ordinals inside ${\mathbb V}^{(\mathbb B)}$.

  \subsubsec{(1)} Clearly, $\Ord(x)$ is a~bounded formula. Since
   $\lim(\alpha)\le\alpha$ for every ordinal~$\alpha$,
 the formula $\Ord(x)\wedge x=\lim(x)$ may be rewritten as
  $ \Ord(x)\wedge (\forall\,t\in x)(\exists\,s\in x)(t\in s)$.
  Hence, $\Ord(x)\wedge x=\nobreak\lim(x)$ is a bounded formula as well.
Finally, the record
 $$
 \Ord(x)\wedge x=\lim(x)\wedge (\forall\,t\in x)(t=\lim(t)\to t=0)
 $$
 convinces us that the ``least limit ordinal'' is
 a~bounded formula too. Hence
 $\alpha$ is the least limit ordinal
 if and only if ${\mathbb V}^{(\mathbb B)}$ $\models$
 ``$\alpha^{\scriptscriptstyle\wedge}$ is~the least limit ordinal.''
Since $\omega$ is the least limit ordinal,
  ${\mathbb V}^{(B)}\models$
 ``$\omega^{\scriptscriptstyle\wedge}$ is the least limit ordinal.''

 \subsubsec{(2)} It can be demonstrated that
 ${\mathbb V}^{(\mathbb B  )}\models$
 ``$\On^{\scriptscriptstyle\wedge}$ is the unique ordinal class
 that is not an ordinal'' (with $\On^{\scriptscriptstyle\wedge}$ defined in an appropriate way).
 Given $x\in{\mathbb V}^{(\mathbb B  )}$, we thus have
 $$
 [\![\Ord(x)]\!]=\bigvee_{\alpha\in\On}[\![x=\alpha^{\scriptscriptstyle\wedge} ]\!].
 $$

  \subsubsec{(3)} \proclaim{}
 Each ordinal inside ${\mathbb V}^{(\mathbb B  )}$ is  a~mixing of
 some set of standard ordinals.
 In other words, given $x\in{\mathbb V}^{(\mathbb B  )}$, we have
 ${\mathbb V}^{(\mathbb B  )}\models\Ord(x)$
 if and only if there
are an ordinal $\beta\in\On$ and a partition of unity
 $(b_\alpha)_{\alpha\in\beta}\subset \mathbb B$ such that
  $x=\mix_{\alpha\in\beta}b_\alpha\alpha
^{\scriptscriptstyle\wedge}$.
 \Endproc

  \subsubsec{(4)} This yields the convenient formulas for quantification over
  ordinals:
 $$
 \align {}
 [\![(\forall\,x)\bigl(\Ord(x)\to\psi (x)\bigr)]\!]&=
 \bigwedge\limits_{\alpha\in\On}\, [\![\psi (\alpha^{\scriptscriptstyle\wedge})]\!],
 \\
 [\![(\exists\,x)\bigl(\Ord(x)\wedge\psi (x)\bigr)]\!]&=
 \bigvee\limits_{\alpha\in\On}\, [\![\psi (\alpha^{\scriptscriptstyle\wedge})]\!].
 \endalign
 $$

 \subsec{A3.7.}~By transfer every Boolean valued model enjoys the
 classical principle of  cardinal comparability. In other words,
 there is a~${\mathbb V}^{(\mathbb B  )}$-class $\Cn$ whose elements
 are only cardinals. Let $\Card(\alpha)$ denote the formula that declares $\alpha$
 a cardinal.  Inside ${\mathbb V}^{(\mathbb B  )}$ we then see that
  $\alpha\in\Cn\,\iff\,\Card(\alpha)$. Clearly, the
  class of ordinals $\On^{\scriptscriptstyle\wedge}$ is
  similar to the class of infinite cardinals, and we denote
  the similarity from $\On^{\scriptscriptstyle\wedge}$ into $\Cn$
  by~$\alpha\mapsto\aleph_\alpha$.
 In particular, to each standard ordinal $\alpha\in\On$
 there is a unique infinite cardinal
 $\aleph_{\alpha^{\scriptscriptstyle\wedge}}$ inside
  ${\mathbb V}^{(\mathbb B  )}$. Indeed,
 $[\![\Ord(\alpha^{\scriptscriptstyle\wedge})]\!]=\mathbbm1$.

 Recall that it is customary to refer to
 the standard names of ordinals and cardinals as
  {\it standard ordinals\/} and
  {\it standard cardinals\/}
 inside~$\mathbb V^{(\mathbb B)}$.

 \smallskip
 \subsubsec{(1)}~\proclaim{}The standard name of the least
 infinite cardinal is the least infinite cardinal:
 $$
 {\mathbb V}^{(\mathbb B  )}\models
 (\omega_0)^{\scriptscriptstyle\wedge}=\aleph_0.
 $$
 \Endproc

  Inside ${\mathbb V}^{(\mathbb B )}$ there is a mapping
   $|{\cdot}|$ from the universal class $\mathbb U_{\mathbb B}$
   into the class $\Cn$ such that  $x$ and $|x|$ are equipollent
   for all $x$.

 \subsubsec{(2)}~\proclaim{}The standard names of equipollent sets
 are of the same cardinality:
 $$
 (\forall\, x\in\mathbb V)\,(\forall\, y\in\mathbb V)\,
 \bigl(|x|=|y|\ra[\![|x^{\scriptscriptstyle\wedge}|=
 |y^{\scriptscriptstyle\wedge}|]\!]=\mathbbm1\bigr).
 $$
 \Endproc

 \subsec{A3.8.}~%
 {\subsec\relax(1)}~\proclaim{}If the standard name of
 an ordinal $\alpha$ is a cardinal then
 $\alpha$ is a cardinal too:
 $$
 (\forall\,\alpha\in\On)\ \bigl({\mathbb V}^{(\mathbb B  )}\models\Card(\alpha^{\scriptscriptstyle\wedge})\bigr)
 \ra\Card(\alpha).
 $$
 \Endproc

 \subsubsec{(2)}~\proclaim{}The standard name of a finite
 cardinal is a finite cardinal too:
 $$
 (\forall\,\alpha\in\On)\ \bigl(\alpha<\omega\,\ra\, {\mathbb
 V}^{(\mathbb B  )}\models\Card(\alpha^{\scriptscriptstyle\wedge})
 \wedge\alpha^{\scriptscriptstyle\wedge}\in\aleph_0\bigr).
 $$
 \Endproc

 \subsec{A3.9.}~\proclaim{}Given $x\in{\mathbb V}^{(\mathbb B  )}$,
 we have ${\mathbb V}^{(\mathbb B)}\models\Card(x)$ if and only if
 there are nonempty set of cardinals $\Gamma$ and a partition of unity
 $(b_\alpha)_{\alpha\in\Gamma}\subset \mathbb B $ such that ${\mathbb
 V}^{(\mathbb B )}\models\Card(\gamma^{\scriptscriptstyle\wedge})$
 for all  $\gamma\in\Gamma$ and
 $x=\mix_{\gamma\in\Gamma} b_\gamma\gamma^{\scriptscriptstyle\wedge}$.
 In other words, each Boolean valued cardinal is a mixing of
 some set of standard cardinals. \Endproc

 \goodbreak

 \subsec{A3.10.}~It is worth noting that $\sigma$-distributive Boolean algebras
 are often referred to as $(\omega,\omega)$-distributive Boolean algebras.
 This term is related to a~more general notion, $(\alpha,\beta)$-distributivity,
 where $\alpha$ and $\beta$ are arbitrary cardinals.

 \smallskip
 \proclaim{}If\/ $\mathbb B$ is a~complete Boolean algebras then the~%
 following are equivalent:

 \subsubsec{(1)}~$\mathbb  B$ is $\sigma$-distributive;

 \subsubsec{(2)}~$\mathbb V^{(\mathbb B  )}\models(\aleph_0)^{\aleph_0}=
 (\omega^\omega)^{\scriptscriptstyle\wedge}$;

 \subsubsec{(3)}~$\mathbb V^{(\mathbb B  )}\models\mathcal P(\aleph_0)=
 \mathcal P(\omega)^{\scriptscriptstyle\wedge}$.
 \Endproc
 \smallskip

 The latter is a~result by Scott on $(\alpha,\beta)$-distributive Boolean algebras
 which was formulated in the~case $\alpha=\beta=\omega$ (cp.~\cite[2.14]{Bell}).

 \smallskip

More details and references are collected in~\cite{IBA}.
The~monographs \cite{Bell} and \cite{Jech} are also a~very good source
of facts concerning Boolean valued cardinals and, in~particular, continuum.

} % "\small" ends here

 \Ref
 \def\citeitem#1#2{\bibitem{#2}}

 \begin{enumerate}
 \itemsep=0pt\parskip=0pt
 {\normalsize

 \citeitem{Abramovich1983}{Ab1}    Abramovich Yu.\,A.~
 Multiplicative representation of disjointness preserving operators,
 \textit{Indag. Math.~\!(N.S.)}, \textbf{45}(3) (1983), 265--279.

 \citeitem{AbramovichVK1979}{AVK} Abramovich Yu.\,A., Veksler A.\,I., and Koldunov A.\,V.~
 On~dis\-joint\-ness preserving operators,
 {\it  Dokl. Akad. Nauk SSSR}, {\bf 289}(5) (1979), 1033--1036.

 \citeitem{AbramovichVK1981}{AVK1} Abramovich Yu.\,A., Veksler A.\,I., and Koldunov A.\,V.~
 Disjointness-preserving operators, their continuity, and multiplicative representation,
 In: {\it Linear Operators and Their Applications}, Leningrad Ped. Inst., Leningrad (1981)
 [in Russian].

 \citeitem{AbramovichKi1999}{AK1} Abramovich, Y.\,A. and Kitover, A.\,K.~
 d-Independence and d-bases in Vector lattices, {\it Rev. Romaine de Math. Pures et appliquees}, {\bf 44} (1999), 667-682.

 \citeitem{AbramovichKi2000}{AK3}  Abramovich Y.\,A. and Kitover A.\,K.~
 {\it Inverses of Disjointness Preserving Operators}, Memoirs Amer. Math. Soc., {\bf 679} (2000).

 \citeitem{AbramovichKi2003}{AK2}  Abramovich, Y.\,A. and Kitover~A.\,K.~ % AGu
 d-Independence and d-bases, {\it Positivity}, \textbf{7}(1) (2003), 95--97.

 \citeitem{AbramovichWi1993}{AW}   Abramovich~Yu.\,A. and Wickstead~A.\,W.~
 The regularity of order bounded operators into $C(K)$. {\rm II}, {\it Quart. J. Math. Oxford}, Ser.~2,  {\bf44} (1993), 257--270.

 \citeitem{AczelDhombres}{AD}      Acz\'el J. and Dhombres J.~
 \textit{Functional Equations in Several Variables}, Cambridge Univ. Press, Cambridge etc. (1989).

 \citeitem{AlbeverioAyuKud2007}{AAKud}  Albeverio S., Ayupov Sh.\,A., and Kudaybergenov K.\,K.~
 Derivations on the~algebras of measurable operators affiliated
 with a~type I von Neumann algebra, \textit{to appear}.

 \citeitem{AliprantisBur1985}{AB}   Aliprantis C.\,D. and Burkinshaw O.~
 \textit{Positive Operators}, Acad. Press Inc., London (1985).

 \citeitem{AyupovKud2007}{AKud}     Ayupov Sh.\,A. and Kudaybergenov K.\,K.~
 Derivations and automorphisms of algebras of bounded operators on Banach--Kantorovich spaces,
 \textit{to appear}.

 \citeitem{AzouziBoulBu2006}{ABB}   Azouzi~Y., Boulabiar~K., and Buskes~G.~ % AGu (initials -> end)
 The de Schipper formula and squares of Riesz spaces, \textit{Indag. Math.~\!(N.S.)}, \textbf{17}(4) (2006), 265--279.

 \citeitem{Bell1985}{Bell}       Bell~J.\,L.~
 \textit{Boolean Valued Models and Independence Proofs in Set Theory},
 Clarendon Press, New York etc. (1985).% xx+165~pp.

 \citeitem{BerChilSu2004}{BCS}     Ber~A.\,F., Chilin V.\,I., and Sukochev F.\,A.~
 Derivations in regular commutative algebras, {\it Math.
 Notes}, \textbf{75} (2004) 418--419.

 \citeitem{Berberian1972}{Berb}     Berberian~S.\,K.~
 {\it Baer $\ast$-Rings}, Springer-Verlag, Berlin (1972).%---xii+296~pp.

 \citeitem{Bernau1981}{Ber3}       Bernau~C.\,B.~
 Orthomorphisms of Archimedean vector lattices, {\it Math. Proc. Cambridge Phil. Soc.}, {\bf 89} (1981), 119--126.

 \citeitem{BeukersHuiPag1983}{BHP}  Beukers F., Huijsmans~C.\,B., and  de~Pagter~B.~
 Unital embedding and complexification of $f$-algebras, {\it Math. Z.}, {\bf 183} (1983), 131--144.

 \citeitem{BigardKeimel}{BK}       Bigard~A. and Keimel~K.~
 Sur les endomorphismes conservant les polaires d'un groupe r\'eticul\'e archim\'edien,
 {\it Bull. Soc. Math. France}, {\bf97} (1969),  381--398.

 \citeitem{BigardKeimelWolf}{BKW}  Bigard~A., Keimel~K., and Wolfenstein~S.~
 \textit{Groupes et Anneaux R\'eticul\'es}, Springer-Verlag, Berlin etc. (1977) (Lecture Notes in Math., {\bf 608}.)%---xi+334~p.

 \citeitem{Birkhoff}{Bir}          Birkhoff~G.~
 \textit{Lattice Theory},  3rd ed., Amer. Math. Soc. Colloq. Publ.,
 Providence, {\bf 25}  (1967).

 \citeitem{BirkhoffPierce}{BirP}   Birkhoff~G. and Pierce~R.\,S.~
 Lattice-ordered rings, {\it An. Acad. Brasil Ci\^enc}, {\bf28} (1956), 41--68.

 \citeitem{BoulabiarBuskesTr1}{BBT1}Boulabiar~K., Buskes~G., and Triki~A.~
 Recent results in lattice ordered algebras, In:
 \textit{Contemp. Math.}, {\bf328}, Providence, RI, (2003), 99--133.

 \citeitem{BoulabiarBuskesTr2}{BBT2} Boulabiar~K., Buskes~G., and~Triki~A.~
 Results in $f$-algebras, \textit{Positivity} (2007), 73--96.

 \citeitem{Bourbaki}{Bou1}         Bourbaki~N.~
 \textit{Algebra (Polynomials and Fields, Ordered Groups)},
 Hermann,  Paris (1967) [in~French].

 \citeitem{Buck1961}{Buc}          Buck~R.\,C.~
 Multiplication operators, {\it Pacific J. Math.}, {\bf 11} (1961), 95--103.

 \citeitem{Bukhvalov1988}{Buh1}    Bukhvalov~A.\,V.~
 Order bounded operators in vector lattices and
 spaces of measurable functions, In: {\it Mathematical Analysis}, % AGu (no [in Russian])
 VINITI, Moscow (1988),
 3--63 (Itogi Nauki i Tekhniki, {\bf26}) [in Russian].

 \citeitem{BuskesKusraev2007}{BuK}  Buskes G. and Kusraev A.\,G.~
 Representation and extension of orthoregular bilinear operators, {\it Vladikavkaz Math. J.}, \textbf{9}(1) (2007), 16--29.

 \citeitem{BuskesRedf2007}{BRe}   Buskes G. and Redfield~R.\,H.~
 Square means and geometric means in lattice-ordered groups and vector
 lattices, \textit{to appear}.

 \citeitem{BuskesRooij2000}{BuR2}   Buskes G. and van~Rooij~A.~
 Almost $f$-algebras: commutativity and
 the~Cauchy--Schwarz inequality, {\it Positivity}, \textbf{4}(3) (2000), 227--231.

 \citeitem{BuskesRooij2001}{BuR4}   Buskes~G. and van~Rooij~A.~ % AGu (initials -> end)
 Squares of Riesz spaces, \textit{Rocky Mountain J. Math.,} \textbf{31}(1) (2001), 45--56.

 \citeitem{Chilin1985}{Chil}       Chilin V.\,I.~
 Partially ordered
 involutive Baer algebras, in: \textit{Contemporary Problems of Mathematics.
 Newest Advances}, VINITI, Moscow, {\bf27} (1985), % AGu (year)
 pp.~99--128  [in Russian]. % AGu ([in Russian] -> end)

 \citeitem{ConradDiem}{CD}         Conrad~P.\,F. and Diem~J.\,E.~
 The ring of polar preserving endomorphisms of an~abelian lattice-ordered group,
 {\it Illinois J. Math.},  {\bf15} (1971), 222--240.

 \citeitem{DuhouxMeyer1982}{DuM}   Duhoux~M. and Meyer~M.~
 A new proof of the~lattice structure of orthomorphisms,
 {\it J. London Math. Soc.}, {\bf25}(2) (1982), 375--378.

 \citeitem{Gooderl1979}{Go}        Goodearl K.\,R.~
 \textit{Von Neumann Regular Rings}, Pitman, London (1979).

 \citeitem{Gordon1977}{Gor1}       Gordon~E.\,I.~
 Real numbers in Boolean valued models of set theory and $K$-spaces,
 {\it  Dokl. Akad. Nauk SSSR}, {\bf237}(4) (1977), 773--775.

 \citeitem{Gordon1981}{Gor2}       Gordon~E.\,I.~
 $K$-spaces in Boolean valued models of set theory,
 {\it Dokl. Akad. Nauk SSSR}, {\bf258}(4) (1981), 777--780.

 \citeitem{Gordon1982}{Gor3}       Gordon~E.\,I.~
 To the~theorems of identity preservation  in~$K$-spaces,
 {\it Sibirsk. Mat. Zh.}, {\bf23}(5) (1982), 55--65.

 \citeitem{Gordon1997}{Gord}   Gordon~E.\,I.~
 \textit{Nonstandard
 Methods in Commutative Harmonic Analysis,} Providence, RI: Amer. Math. Soc., (1997).

 \citeitem{Gutman1995}{Gut5} Gutman A.\,E.~
 Banach bundles in the~theory of lattice-normed spaces,
 In: {\it Linear Operators Compatible with Order},
 Sobolev Institute Press, Novosibirsk (1995), 63--211 [in Russian]. % AGu ([in Russian -> end])

 \citeitem{Gutman1995}{Gut6}       Gutman A.\,E.~
 Locally one-dimensional $K$-spaces and $\sigma$-distributi\-ve
 Boole\-an algebras, {\it Siberian Adv. Math.}, \textbf{5}(2) (1995), 99--121.

 \citeitem{Gutman1996}{Gut1}       Gutman~A.\,E.~
 Disjointness preserving operators,
 In: {\it Vector Lattices and Integral Operators} (Ed.: Kutateladze~S.\,S.), % AGu (inirtials -> end)
 Kluwer, Dordrecht etc. (1996), 361--454.

 \citeitem{HuijsmansPag1982}{HP3}  Huijsmans~C.\,B. and  de~Pagter~B.~
 Ideal theory in $f$-algebras, {\it Trans. Amer. Math. Soc.}, {\bf 269} (1982), 225--245.

 \citeitem{HuijsmansWic1992}{HW}  Huijsmans~C.\,B. and Wickstead A.\,W.~
 The inverse of band preserving and disjointness preserving operators
 Indag. Math. 179-183.  {\bf 3}(2) (1992), 179--183.

 \citeitem{Jech2002}{Jech}  Jech~T.\,J.~
 \textit{Set Theory. The~Third Millennium Edition}, Springer, Berlin etc. (2002).

 \citeitem{KorolChilin}{KCh}   Korol$'$~A.\,M. and Chilin~V.\,I.~
 Measurable operators in a~Boolean-valued model of set theory,
 \textit{Dokl. Akad. Nauk UzSSR,} (3) (1989), 7--9 [in Russian].

 \citeitem{Kusraev2000}{DOP}       Kusraev~A.\,G.~ % AGu (initials -> end)
 \textit{Dominated Operators}, Kluwer, Dordrecht (2000).

 \citeitem{Kusraev2004}{Kus1}      Kusraev A.\,G.~
 On band preserving operators, {\it  Vladikavkaz Math. J.}, {\bf6}(3) (2004), 47--58
 [in~Russian].

 \citeitem{Kusraev2005}{Kus3}      Kusraev A.\,G.~
 Automorphisms and derivations in the~%
 algebra of complex measurable functions, {\it Vladikavkaz Math. J.}, {\bf7}(3) (2005), 45--49 [in Russian].

 \citeitem{Kusraev2006}{K15}        Kusraev A.\,G.~
 On the~structure of orthosymmetric bilinear operators
 in vector lattices, {\it Dokl. RAS}, \textbf{408}(1) (2006), 25--27.

 \citeitem{Kusraev2006}{Kus2}      Kusraev~A.\,G.~
 Automorphisms and derivations in extended complex $f$-algebras,
 {\it Siberian Math. J.}, {\bf~47}(1) (2006), 97--107.

 \citeitem{Kusraev2006}{Kus4}      Kusraev~A.\,G.~
 Analysis, algebra, and logics in operator theory,
 In: \textit{Complex Analysis, Operator Theory, and Mathematical Modeling}
 (Eds.: Korobe\u{\i}nik Yu.\,F. and Kusraev A.\,G.), Vladikavkaz Scientific Center, Vladikavkaz (2006), 171--204 [in Russian].

 \citeitem{Kusraev2007}{K18}      Kusraev A.\,G.~
 When are all separately band preserving bilinear operators symmetric?
 {\it Vladikavkaz Math. J.}, \textbf{9}(2) (2007), 22--25.

 \citeitem{KusraevKut1999}{BVA} Kusraev~A.\,G. and Kutateladze~S.\,S.~
 {\it Boolean Valued Analysis},  Nauka, Novosibirsk (1999), Kluwer, Dordrecht (1999).

 \citeitem{KusraevKut2005}{IBA}    Kusraev~A.\,G. and Kutateladze~S.\,S.~
 \textit{Introduction to Boolean Valued Analysis}, Nauka, Moscow (2005) [in Russian].

 \citeitem{Lambek1966}{Lam}        Lambek\,J.~
 \textit{Lectures on Rings and Modules}, Blaisdell, Toronto (1966).

 \citeitem{Luxemburg1979}{Lux2}    Luxemburg~W.\,A.\,J.~
 {\it Some Aspects of the~Theory of Riesz Spaces}
 (Univ. Arkansas Lect. Notes Math. {\bf 4}), Fayetteville, Arkansas (1979).

 \citeitem{LuxemburgSchep}{LS}     Luxemburg~W.\,A.\,J. and Schep~A.~
 A Radon-Nikod\'ym type theorem for positive operators and a~dual, {\it Indag. Math.}, {\bf40} (1978), 357--375.

 \citeitem{McPolinWick1985}{MW}    McPolin~P.\,T.\,N. and Wickstead~A.\,W.~
 The order boundedness of band preserving operators on uniformly complete vector lattices,
 {\it Math. Proc. Cambridge Philos. Soc.},  {\bf97}(3) (1985), 481--487.

 \citeitem{Meyer1976}{Mey1}        Meyer~M.~
 Le stabilisateur d'un espace vectoriel r\'eticul\'e,
 {\it C. R. Acad. Sci. Paris}, S\'er.~A, {\bf283} (1976), 249--250.

 \citeitem{Nakano1941}{Nak2}       Nakano~H.~
 Teilweise geordnete algebra, {\it Japan J. Math.}, {\bf 17} (1950), 425--511.

 \citeitem{Nakano1953}{Nak3}       Nakano~H.~
 Product spaces of semi-ordered linear spaces,
 {\it J. Fac. Sci. Hokkaido Univ.}, Ser.~I., {\bf12} (1953), 163--210.

 \citeitem{Pagter1984}{Pag1}       de~Pagter~B.~
 A note on disjointness preserving operators,
 {\it Proc.  Amer. Math. Soc.}, {\bf 90}(4) (1984), 543--549.

 \citeitem{Pagter1984}{Pag2}       de~Pagter~B.~
 The space of extended orthomorphisms in Riesz space,
 {\it Pacific J. Math.}, {\bf 112}(1) (1984), 193--210.

 \citeitem{Sakai}{Sak}        Sakai\,S.~
 \textit{$C^*$-algebras and $W^*$-algebras}, Springer-Ver\-lag, Berlin etc. (1971).%---256~pp.

 \citeitem{Schaefer1974}{Sch}      Schaefer~H.\,H.~
 \textit{Banach Lattices and Positive Operators}, Springer-Verlag, Berlin etc. (1974).%---376~pp.

 \citeitem{Sikorski1964}{Sik}      Sikorski R.~
 \textit{Boolean Algebras}, Springer-Verlag, Berlin etc. (1964).

 \citeitem{Waerden1977}{Waer}     Van der Waerden~B.\,L.~
 \textit{Algebra}, Springer, Berlin etc. (1964). [English translation: New York
 (1949)].

 \citeitem{Wickstead1977}{Wic1}    Wickstead~A.\,W.~
 Representation and duality of multiplication
 operators on Archimedean Riesz spaces, {\it Compositio Math.}, {\bf35}(3) (1977), 225--238.

 \citeitem{Wickstead1980}{Wic3}    Wickstead~A.\,W.~
 Extensions of orthomorphisms, {\it J. Austral Math. Soc.,} Ser A , \textbf{29} (1980), 87--98.

 \citeitem{Wickstead1987}{Wic2}    Wickstead~A.\,W.~
 The injective hull of an~Archimedean $f$-algebra,
 {\it Compositio Math.}, {\bf62}(4) (1987), 329--342.

 \citeitem{Zaanen1975}{Z2}         Zaanen~A.\,C.~
 Examples of orthomorphisms,
 {\it J. Approx. Theory},  {\bf13}(2) (1975), 192--204.

 \citeitem{Zaanen1983}{Z}          Zaanen~A.\,C.~
 \textit{Riesz Spaces}, {\bf2}, North-Holland, Amsterdam etc. (1983).%---720~pp.

 \citeitem{ZariskiSamuel1963}{ZS}  Zariski O. and Samuel~P.~
 \textit{Commutative Algebra},
 Springer-Verlag, Berlin (1991).
 }

 \end{enumerate}

 %\end{document}

 %\newpage\thispagestyle{empty}
 \normalsize
 %\vspace*{40mm}

 \vspace*{10mm}

  \noindent
 Alexander E.~Gutman\\
 Sobolev Institute of Mathematics\\
 Siberian Division of the RAS\\
 Novosibirsk, 630090, RUSSIA\\
 E-mail: gutman@math.nsc.ru

 \vspace*{5mm}

 \noindent
 Anatoly G.~Kusraev\\
 Institute of Applied Mathematics and Informatics\\
 Vladikavkaz Science Center of the RAS\\
 Vladikavkaz, 362040, RUSSIA\\
 E-mail: kusraev@smath.ru

 \vspace*{5mm}

 \noindent
 Sem\"en S.~Kutateladze\\
 Sobolev Institute of Mathematics\\
 Siberian Division of the RAS\\
 Novosibirsk, 630090, RUSSIA\\
  E-mail: sskut@math.nsc.ru

 \end{document}